\documentclass[11pt,a4paper]{article}

\usepackage{epsf,epsfig,amsfonts,amsgen,amsmath,amstext,amsbsy,amsopn,amsthm
}
\usepackage{hyperref}
\hypersetup{
    colorlinks = true,     
    linkcolor = blue,      
    urlcolor = blue,       
    citecolor = blue,      
}
\usepackage{cleveref}
\usepackage{amssymb}
\usepackage{graphics}
\usepackage{latexsym,bm}
\usepackage{bbding}
\usepackage{indentfirst}
\usepackage{graphicx}
\usepackage{geometry}
\usepackage{mathrsfs}
\usepackage{exscale}
\usepackage{relsize}
\usepackage{tikz}
\usepackage{float}
\usepackage{mathtools,booktabs,placeins,needspace,xcolor}
\newtheorem{theorem}{Theorem}[section]

\newtheorem{lemma}{Lemma}
\newtheorem{false statement}{False statement}

\theoremstyle{definition}

\newtheorem{remark}{Remark}
\newtheorem{conjecture}{Conjecture}
\newtheorem{corollary}{Corollary}
\newtheorem{problem}{Problem}

\usetikzlibrary{calc,shapes.geometric}
\allowdisplaybreaks[2]
\AddToHook{env/proof/begin}{%
  \setlength{\abovedisplayskip}{10pt plus 3pt minus 6pt}%
  \setlength{\belowdisplayskip}{10pt plus 3pt minus 6pt}%
}
\tikzset{
 vertex/.style={circle,fill=black,inner sep=1.7pt,outer sep=1pt},
 block/.style={ellipse,draw,fill=white,minimum width=20mm,minimum height=10mm,inner sep=2pt,align=center},
 papergraph/.style={graphline,font=\small},
 graphline/.style={line width=.55pt},
 arbitrary/.style={densely dashed,line width=.55pt}
}
\newcommand{\one}{\mathbf1}

\begin{document}

\title{\bfseries Nordhaus-Gaddum type inequalities for Laplacian eigenvalues of connected graphs and its complement}

\author{{\Large Qi Chen,\,\, Ji-Ming Guo\thanks{Corresponding author. Email addresses: \href{mailto:Y30231282@mail.ecust.edu.cn}{Y30231282@mail.ecust.edu.cn} (Qi Chen), \href{mailto:jimingguo@hotmail.com}{jimingguo@hotmail.com} (Ji-Ming Guo), \href{mailto:walkerwzw@163.com}{walkerwzw@163.com} (Zhiwen Wang).},\,\, Zhiwen Wang}\\[2mm]
\small School of Mathematics, East China University of Science and Technology, \\
\small Shanghai, 200237, China}

\date{}

\maketitle
\begin{abstract}
For a graph $G$ of order $n$, let $\alpha(G)$ denote its algebraic connectivity and $\bar{G}$ its complement.
Einollahzadeh and Karkhaneei [Journal of Combinatorial Theory, Series B, 2021] and Afshari [Journal of Graph Theory, 2023] used different methods to prove that $\alpha(G)+\alpha(\bar{G})\geq1$ for $n\geq2$.
In this paper, we first generalize the above result and settle a conjecture of Barrett et al. [Linear Algebra and its Applications, 2022].
If $n\geq4$ and both $G$ and $\bar{G}$ are connected, then
\begin{align}
    \alpha(G)+\alpha(\bar{G})\geq n+a_{n}-c_{n}\quad(>1.38\text{ for }n\geq6), \notag
\end{align}
where $a_{n}$ and $c_{n}$ are the smallest and largest roots of $x^{3}-(n+2)x^{2}+(3n-2)x-n$, respectively.
Equality holds if and only if $G$ or $\bar{G}$ is isomorphic to a graph $\mathcal D(H)$ obtained from a graph $H$ of order $n-3$ by adding a disjoint path $rsq$ and joining $r$ to every vertex of $H$.

Next, we prove that if $n\geq4$ and both $G$ and $\bar{G}$ are connected, then
\begin{align}
     \mu_{1}(G)\mu_{1}(\bar{G})\leq c_{n}(n-a_{n}), \notag
\end{align}
where $\mu_{1}(G)$ denotes the Laplacian spectral radius of $G$.
Equality holds if and only if $G$ or $\bar{G}$ is isomorphic to a graph $\mathcal D(H)$, where $H$ is an arbitrary graph of order $n-3$.
Our result sharpens the upper bound $n(n-1)$ proved for arbitrary graphs by Chen et al. [The Electronic Journal of Combinatorics, 2025].

Finally, as an application of our sum bound, we prove that every nonstar tree $T$ of order $n\geq4$ with maximum degree $\Delta$ satisfies
\begin{align}
    \alpha(\bar{T})\geq n+a_{n}-c_{n}-
    \min\left\{\frac{\Delta-(3\sqrt{2}-4)}{n-1},\,a_{2n-2-\Delta}\right\}. \notag
\end{align}
Equality holds if and only if $T\cong\mathcal D((n-3)K_{1})$.
This result refines the lower bound $\alpha(\bar{T})\geq n-c_{n}$ of Jiang et al. [Discrete Dynamics in Nature and Society, 2015].

		\par\vspace{2mm}
		\par\vspace{1mm}
\end{abstract}

\section{Introduction}
All graphs considered in this paper are finite and simple.
Let $G=(V(G),E(G))$ be a graph with vertex set $V(G)=\{v_{1},\ldots,v_{n}\}$ and edge set $E(G)$, and $\bar{G}$ be its complement.
Write $|G|$ for the order of $G$.
For $v\in V(G)$, denote its neighborhood and degree by $N_{G}(v)$ and $d_{G}(v)$, respectively, and put $N_{G}[v]=N_{G}(v)\cup\{v\}$.
A vertex of degree one is called a pendant vertex.
Write $d_{G}(u,v)$ for the distance between vertices $u,v$, and $d(G)$ for the diameter of $G$.
Let $K_{m}$ denote the complete graph of order $m$.
For vertex-disjoint graphs $G$ and $H$, write $G\cup H$ for their disjoint union.
Write $kF$ for the disjoint union of $k$ copies of a graph $F$.
The join $G\vee H$ is obtained from $G\cup H$ by adding all edges between $V(G)$ and $V(H)$.

Let $A(G)$ be the adjacency matrix of $G$, and $D(G)=\operatorname{diag}(d_{G}(v_{1}),\allowbreak\ldots,d_{G}(v_{n}))$ be its degree matrix.
The Laplacian matrix of $G$ is $L(G)=D(G)-A(G)$.
It is positive semidefinite, and its eigenvalues are written in nonincreasing order as $\mu_{1}(G)\geq\cdots\geq\mu_{n-1}(G)\geq\mu_{n}(G)=0$.
The Laplacian spectral radius of $G$ is the largest eigenvalue $\mu_{1}(G)$ of $L(G)$.
The second smallest eigenvalue $\alpha(G)=\mu_{n-1}(G)$ is called the algebraic connectivity of $G$~\cite{Fiedler1973}.
A unit eigenvector of $L(G)$ corresponding to $\alpha(G)$ is called a Fiedler vector.

Nordhaus and Gaddum~\cite{NG1956} established upper and lower bounds on the sum and product of the chromatic numbers of a graph and its complement.
For a graph parameter $p$, bounds on $p(G)+p(\bar{G})$ or $p(G)p(\bar{G})$ are called Nordhaus--Gaddum type inequalities.

Put $S(G)=\alpha(G)+\alpha(\bar{G})$.
Zhai et al.~\cite{ZSH2011}, as well as You and Liu~\cite{YL2012}, conjectured that $S(G)\geq1$ for every graph of order at least two.
The conjecture was proved for trees~\cite{Fan2008}, unicyclic graphs~\cite{BTF2009}, bicyclic graphs~\cite{FLT2010}, tricyclic graphs~\cite{CW2009}, cactus graphs~\cite{Liu2010}, certain quasi-tree graphs~\cite{XM2011}, connected graphs of diameter other than three~\cite{ZSH2011}, bipartite graphs~\cite{ATR2014}, and $K_{3}$-free graphs~\cite{CD2016}.
Einollahzadeh and Karkhaneei~\cite{EK2021} proved the conjecture in full and characterized all equality cases.
Afshari~\cite{Afshari2023} gave another proof of the inequality.

\begin{figure}[H]
\centering
\begin{tikzpicture}[papergraph,scale=.92,transform shape]
\begin{pgfinterruptboundingbox}
\node[block,minimum width=12mm,minimum height=19mm,inner sep=0pt] (H) at (0,0) {};
\node at (-.28,0) {$H$};
\node[font=\scriptsize,inner sep=0pt] at (.16,0) {$\vdots$};
\node[font=\scriptsize] at (0,-1.15) {$|H|=n-3$};
\node[font=\scriptsize] at (0,-1.52) {$|E(H)|\geq0$};
\node[vertex] (h1) at (0.16,0.46) {};
\node[vertex] (h2) at (0.16,-0.46) {};
\node[vertex,label=above:$r$] (r) at (1.65,0) {};
\node[vertex,label=above:$s$] (s) at (2.9,0) {};
\node[vertex,label=above:$q$] (q) at (4.15,0) {};
\draw (h1)--(r);
\draw (h2)--(r);
\draw (r)--(s);
\draw (s)--(q);
\draw (0,0) ellipse [x radius=6mm,y radius=9.5mm];
\foreach \v in {h1,h2,r,s,q} {\fill (\v) circle[radius=1.7pt];}
\end{pgfinterruptboundingbox}
\path[use as bounding box] (-1.0766,-1.74) rectangle (4.3112,1.0267);
\end{tikzpicture}
\caption{$\mathcal D(H)$.}
\label{fig:dandelion}
\end{figure}
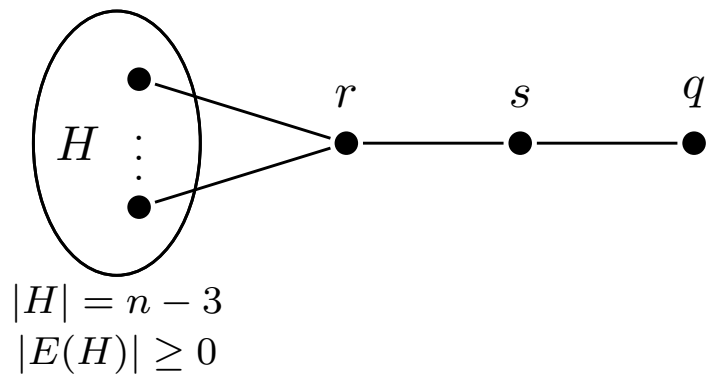

For a graph $H$ of order $n-3$ ($n\geq4$), let $\mathcal D(H)$ be obtained from $H$ by adding a path $rsq$ and joining $r$ to every vertex of $H$ (see Figure~\ref{fig:dandelion}).

\FloatBarrier

Barrett et al.~\cite[Conjecture~2]{BEHK2022} posed the following conjecture.
\begin{conjecture}\label{conj:BEHK2022}
For each $n\geq4$ and any graph $H$ of order $n-3$, the graph $\mathcal D(H)$ and its complement minimize $\alpha(G)+\alpha(\bar{G})$ among all graphs $G$ of order $n$ with $d(G)=d(\bar{G})=3$.
\end{conjecture}
The following theorem implies Conjecture~\ref{conj:BEHK2022}.

\begin{theorem}\label{thm:main}

Let $G$ be a graph of order $n\geq4$.
If both $G$ and $\bar{G}$ are connected, then
\begin{align}
    \alpha(G)+\alpha(\bar{G})\geq n+a_{n}-c_{n}\quad(>1.38\text{ for }n\geq6), \notag
\end{align}
where $a_{n}$ and $c_{n}$ are the smallest and largest roots, respectively, of the polynomial $p_{n}(x)=x^{3}-(n+2)x^{2}+(3n-2)x-n$.
Equality holds if and only if $G$ or $\bar{G}$ is isomorphic to a graph $\mathcal D(H)$, where $H$ is an arbitrary graph of order $n-3$.

\end{theorem}

Ashraf and Tayfeh-Rezaie~\cite{ATR2014} proved $\mu_{1}(G)\mu_{1}(\bar{G})\leq n(n-1)$ for bipartite graphs and conjectured the same bound for arbitrary graphs.
Chen et al.~\cite{CGLW2025} proved the conjecture, with equality if and only if $G$ or $\bar{G}$ is the join of $K_{1}$ and a disconnected graph $H$ of order $n-1$.

For connected graphs with connected complements, we sharpen this bound as follows.

\begin{theorem}\label{thm:laplacian-product}
Let $G$ be a graph of order $n\geq4$.
If both $G$ and $\bar{G}$ are connected, then
\begin{align}
    \mu_{1}(G)\mu_{1}(\bar{G})\leq c_{n}(n-a_{n}),   \notag
\end{align}
where $a_{n}$ and $c_{n}$ are as in Theorem~\ref{thm:main}.
Equality holds if and only if $G$ or $\bar{G}$ is isomorphic to a graph $\mathcal D(H)$, where $H$ is an arbitrary graph of order $n-3$.
\end{theorem}

Write $\Delta(G)$ for the maximum degree of $G$.
Jiang et al.~\cite{JYC2015} showed that every nonstar tree $T$ of order $n\geq4$ satisfies $\alpha(\bar{T})\geq n-c_{n}$, with equality precisely when $T\cong\mathcal D((n-3)K_{1})$.

As an application of Theorem~\ref{thm:main}, we obtain an explicit lower bound for $\alpha(\bar{T})$ in terms of $n$ and $\Delta(T)$.
This refines the lower bound $n-c_{n}$.

\begin{theorem}\label{thm:tree-complement}
Let $T$ be a nonstar tree of order $n\geq4$ with maximum degree $\Delta$.
Then
\begin{align}
    \alpha(\bar{T})\geq n+a_{n}-c_{n}-
    \min\left\{\frac{\Delta-(3\sqrt{2}-4)}{n-1},\,a_{2n-2-\Delta}\right\},            \notag
\end{align}
where $a_{n}$, $c_{n}$, and $a_{2n-2-\Delta}$ are as in Theorem~\ref{thm:main}.
Equality holds if and only if $T\cong\mathcal D((n-3)K_{1})$.
\end{theorem}

\section{Preliminaries}\label{sec:prelim}

For $X\subseteq V(G)$, let $G[X]$ be the subgraph of $G$ induced by $X$, and put $G-X=G[V(G)\setminus X]$.
We write $G-e$ ($G-v$) for the graph obtained from $G$ by deleting an edge $e$ (a vertex $v$).
Let $K_{r,s}$ and $P_{n}$ denote, respectively, the complete bipartite graph with parts of orders $r$ and $s$ and the path of order $n$.
Write $\nu(G)$ for the matching number of $G$.

We use $I$ and $\one$ for the identity matrix and all-ones vector of the required sizes, and $\|\mathbf{x}\|$ for the Euclidean norm.
For $\mathbf{x}=(x_{1},\ldots,x_{n})\in\mathbb R^{n}$, write $x_{v}$ or $x(v)$ for the component corresponding to $v\in V(G)$, so that $x_{v_{i}}=x(v_{i})=x_{i}$.
For $X\subseteq V(G)$, let $\mathbf{x}_{X}$ denote the subvector of $\mathbf{x}$ indexed by $X$.
The Dirichlet energy of $\mathbf{x}\in\mathbb R^{n}$ on $G$ is $\mathcal E_{G}(\mathbf{x})=\sum_{v_{i}v_{j}\in E(G)}(x_{i}-x_{j})^{2}=\mathbf{x}^{\mathsf T}L(G)\mathbf{x}$.

For disjoint sets $A,B\subseteq V(G)$, write $E_{G}(A,B)=\{uv\in E(G):u\in A,v\in B\}$ and let $G[A,B]$ be the bipartite graph on $A\cup B$ with edge set $E_{G}(A,B)$.
We use the same notation for a vertex in place of a singleton set.
An edge $uv\in E(G)$ is called dominating if $N_{G}(u)\cup N_{G}(v)=V(G)$.
A vertex cover of $G$ is a set $W\subseteq V(G)$ containing at least one endpoint of every edge of $G$.

For any graph $G$ of order $n\geq 2$, we have
\begin{equation}\label{eq:variational}
\alpha(G)
=\min_{\substack{\mathbf{x}\in\mathbb R^{n}\setminus\{\mathbf{0}\}\\
		\sum_{i=1}^{n} x_{i}=0}}
\frac{\mathcal E_{G}(\mathbf{x})}{\|\mathbf{x}\|^{2}}.
\end{equation}

\begin{lemma}[{\cite[Corollary~3.2]{Fiedler1973}}]\label{lem:subgraph}
If $H$ is a spanning subgraph of a graph $G$ of order $n\geq2$, then $\alpha(H)\leq\alpha(G)$.
\end{lemma}

\begin{lemma}[{\cite[Section~1.3.2]{BH2012}}]\label{lem:complement-spectrum}
	For any graph $G$ of order $n\geq2$, $\mu_{k}(\bar{G})=n-\mu_{n-k}(G)$ for $1\leq k\leq n-1$.
	In particular, $\alpha(\bar{G})=n-\mu_{1}(G)$.
\end{lemma}

The following lemma follows from~\cite[Lemma~4 and the proof of Theorem~1]{EK2021}.

\begin{lemma}\label{lem:fiedler-range}

Let $G$ and $\bar{G}$ be connected graphs, and $\mathbf{x}=(x_{1},\ldots,x_{n})$ and $\mathbf{y}=(y_{1},\ldots,y_{n})$ be their Fiedler vectors, respectively.
Suppose that $x_{p}=\max _{1 \leq i \leq n} x_{i}, \quad x_{q}=\min _{1 \leq i \leq n} x_{i}$.
If $x_{p}-x_{q}\geq\max_{u,v\in V(G)}|y_{u}-y_{v}|$, then $S(G)\geq\frac{1}{(x_{p}-x_{q})^{2}}$.

\end{lemma}

\begin{lemma}[{\cite[Lemma~5]{EK2021}}]\label{lem:diameter-bound}
	If $G$ is a connected graph of order $n\geq2$ with $d(G)\leq2$, then $\alpha(G)\geq1$.
\end{lemma}

\begin{lemma}\label{lem:pendant}
Let $G'$ be obtained from a graph $G$ of order $n-1$ ($n\geq3$) by adding a pendant vertex.
If $0<c\leq\alpha(G)$, then $\alpha(G')\geq\frac{c+2-\sqrt{(c+2)^{2}-\frac{4nc}{n-1}}}{2}$.
\end{lemma}

\begin{proof}
\begingroup
Put $a=\alpha(G')$.
Since $c\leq\alpha(G)\leq n-1$ by Lemma~\ref{lem:complement-spectrum}, we have $(c+2)^{2}-\frac{4nc}{n-1}=c^{2}+4\left(1-\frac{c}{n-1}\right)\geq c^{2}$.
Hence $\frac{c+2-\sqrt{(c+2)^{2}-\frac{4nc}{n-1}}}{2}\leq1$, so we may assume that $a<1$.
Applying the estimate in the proof of Lemma~3.4 in~\cite{Kirkland2010} to $G'$, we have
\begin{align}
    c\leq\alpha(G)
    &\leq a+\frac{a(n-2)\left(2-a-\frac{n-2}{n-1}\right)}
    {(2-a)^{2}-\frac{n-2}{n-1}+(n-2)(1-a)^{2}} \notag\\
    &=a+\frac{a(n-2)}{n-(n-1)a}
    =\frac{(n-1)a(2-a)}{n-(n-1)a}. \notag
\end{align}
Since $n-(n-1)a>0$, we have $c[n-(n-1)a]\leq(n-1)a(2-a)$, that is, $a^{2}-(c+2)a+\frac{nc}{n-1}\leq0$.
Since $\frac{c+2-\sqrt{(c+2)^{2}-\frac{4nc}{n-1}}}{2}$ is the smaller root of $a^{2}-(c+2)a+\frac{nc}{n-1}=0$, it follows that $a\geq\frac{c+2-\sqrt{(c+2)^{2}-\frac{4nc}{n-1}}}{2}$.
\endgroup
\end{proof}

\begin{lemma}[{\cite[Theorem~1]{Haynsworth1968}}]\label{lem:schur}
	Let $M=\left(\begin{smallmatrix}A&C\\C^{\mathsf T}&D\end{smallmatrix}\right)$ be a real symmetric matrix with $D$ positive definite.
	Put $K=A-CD^{-1}C^{\mathsf T}$.
	Then we have $n_-(M)=n_-(K)$ and $n_{0}(M)=n_{0}(K)$, where $n_-$ and $n_{0}$ denote the numbers of negative and zero eigenvalues, respectively, counted with multiplicity.
\end{lemma}

The following lemma is an immediate consequence of the Courant--Fischer theorem~\cite[Theorem~2.4.1]{BH2012}.

\begin{lemma}\label{lem:hyperplane}
	Let $G$ be a graph of order $n\geq2$ and $\tau>0$.
	If $\mathcal E_{G}(\mathbf{x})\geq\tau\|\mathbf{x}\|^{2}$ for all $\mathbf{x}$ in an $(n-1)$-dimensional subspace $H\subseteq\mathbb R^{n}$, then $\alpha(G)\geq\tau$.
\end{lemma}

\begin{lemma}\label{lem:normalization}
Let $G$ be a connected graph of order $n\geq2$ and $\mathbf{x}$ be a Fiedler vector.
Choose $p,q\in V(G)$ such that $x_{p}=\max_{v\in V(G)}x_{v}$ and $x_{q}=\min_{v\in V(G)}x_{v}$.
Put $\mathbf{z}=\frac{\mathbf{x}-x_{q}\one}{x_{p}-x_{q}}$.
Then the following hold.
\begin{enumerate}
\item[\textup{(i)}] We have $\alpha(G)\geq|N_{G}(p)\cap N_{G}(q)|$, and $\alpha(G)\geq2$ if $pq\in E(G)$.

\item[\textup{(ii)}] If $\bar{G}$ is connected and has a Fiedler vector $\mathbf{y}$ such that $x_{p}-x_{q}\geq\max_{u,v\in V(G)}|y_{u}-y_{v}|$, then
\begin{equation}\label{eq:energy}
\mathcal E_{G}(\mathbf{z})=\frac{\alpha(G)}{(x_{p}-x_{q})^{2}}\leq\alpha(G)S(G).
\end{equation}
\end{enumerate}
\end{lemma}

\begin{proof}
	\textup{(i)} Since $\sum_{i=1}^{n} x_{i}=0$ and $\|\mathbf{x}\|=1$, we have $x_{p}>0>x_{q}$.
	The definition of $\mathbf{z}$ gives $z_{p}=1$, $z_{q}=0$ and $0\leq z_{i}\leq1$ for $1\leq i\leq n$.
	Also, $n^{-1}\sum_{i=1}^{n} z_{i}=-\frac{x_{q}}{x_{p}-x_{q}}$, so $L(G)\one=\mathbf{0}$ and $L(G)\mathbf{x}=\alpha(G)\mathbf{x}$ give
	\begin{equation}\label{eq:vertex-equations}
		L(G)\mathbf{z}=\alpha(G)\left(\mathbf{z}-\frac1n\left(\sum_{i=1}^{n} z_{i}\right)\one\right).
	\end{equation}
	Since $z_{p}-z_{q}=1$, Eq.~\eqref{eq:vertex-equations} gives
	\begin{equation}\label{eq:extreme-sum}
		\alpha(G)
		=(L(G)\mathbf{z})_{p}-(L(G)\mathbf{z})_{q}=\sum_{v\in N_{G}(p)}(1-z_{v})+\sum_{v\in N_{G}(q)}z_{v}.
	\end{equation}
	Hence $\alpha(G)\geq\sum_{v\in N_{G}(p)\cap N_{G}(q)}((1-z_{v})+z_{v})=|N_{G}(p)\cap N_{G}(q)|$.
	If $pq\in E(G)$, then $\alpha(G)\geq(1-z_{q})+z_{p}=2$.
	
	\textup{(ii)} Since $\|\mathbf{x}\|=1$ and $L(G)\one=\mathbf{0}$, we have $\mathcal E_{G}(\mathbf{z})=\frac{\mathcal E_{G}(\mathbf{x})}{(x_{p}-x_{q})^{2}}=\frac{\alpha(G)}{(x_{p}-x_{q})^{2}}$.
	By Lemma~\ref{lem:fiedler-range}, $(x_{p}-x_{q})^{-2}\leq S(G)$, which gives Eq.~\eqref{eq:energy}.
\end{proof}

Under the hypotheses and notation of Lemma~\ref{lem:normalization}\textup{(ii)}, suppose that $d_{G}(p,q)=3$, and put $U=N_{G}(p)$, $V=N_{G}(q)$ and $F=G[U,V]$.
Then $U\cap V=\varnothing$ and $E(F)\ne\varnothing$.
Put $t(v)=1-z(v)$ for $v\in U$ and $t(v)=z(v)$ for $v\in V$.
Thus $t(v)\geq0$ for every $v\in U\cup V$.
Eq.~\eqref{eq:extreme-sum} gives
\begin{equation}\label{eq:endpoint-sums}
	\sum_{v\in U\cup V}t(v)=\alpha(G).
\end{equation}

Let $F_{0}$ be a subgraph of $G[U\cup V]$ with at least one edge between $U$ and $V$ and no isolated vertices.
Define graphs $F_{1}$ and $F_{2}$ such that $V(F_{1})=V(F_{2})=V(F_{0})$, $E(F_{1})=\{uv\in E(F_{0}):u\in U,\ v\in V\}$, and $E(F_{2})=E(F_{0})\setminus E(F_{1})$.
For a real vector $\mathbf{w}=(w(v))_{v\in V(F_{0})}$ indexed by the vertices of $F_{0}$, define
\begin{align}
    C_{F_{0}}(\mathbf{w})
    &=\sum_{v\in V(F_{0})}w(v)^{2}
    +\sum_{uv\in E(F_{1})}(1-w(u)-w(v))^{2}
    +\sum_{uv\in E(F_{2})}(w(u)-w(v))^{2}, \notag\\
    \kappa_{F_{0}}
    &=\min_{\mathbf{w}\in\mathbb R^{|F_{0}|}}C_{F_{0}}(\mathbf{w}). \notag
\end{align}
For $\mathbf{t}=(t(v))_{v\in V(F_{0})}$ with $t(v)$ as defined above, Eq.~\eqref{eq:energy} gives
\begin{equation}\label{eq:retained-energy}
\kappa_{F_{0}}\leq C_{F_{0}}(\mathbf{t})
\leq\mathcal E_{G}(\mathbf{z})\leq\alpha(G)S(G).
\end{equation}

\begin{lemma}\label{lem:energy-bound}
	Let $G$ be a connected graph.
	Let $\mathbf{x},\mathbf{y}$ be Fiedler vectors of $G,\bar{G}$, respectively, and choose $p,q\in V(G)$ such that $x_{p}=\max_{v} x_{v}$ and $x_{q}=\min_{v} x_{v}$.
	Suppose that $d_{G}(p,q)=3$ and $x_{p}-x_{q}\geq\max_{u,v}|y_{u}-y_{v}|$.
	Put $\mathbf{z}=\frac{\mathbf{x}-x_{q}\one}{x_{p}-x_{q}}$.
	Let $F_{0}$ and $\kappa_{F_{0}}$ be as defined above, and put $h=|F_{0}|$.
	Then the following hold.
	\begin{enumerate}
	\item[\textup{(i)}] $h-4\kappa_{F_{0}}>0$.
	\item[\textup{(ii)}] If $\alpha(G)\leq2\kappa_{F_{0}}$, then $\alpha(G)S(G)\mathrel{\geq}\mathcal E_{G}(\mathbf{z})\geq\kappa_{F_{0}}+\frac{(2\kappa_{F_{0}}-\alpha(G))^{2}}{h-4\kappa_{F_{0}}}$.
	\end{enumerate}

\end{lemma}

\begin{proof}
Since $d_{G}(p,q)=3$, $pq$ is a dominating edge of $\bar{G}$, so $\bar{G}$ is connected.
Use $F_{1},F_{2}$ and $C_{F_{0}}$ as defined above.
\textup{(i)} Let $m=|E(F_{1})|$, $\mathbf{d}=(d_{F_{1}}(v))_{v\in V(F_{0})}$, and $M=I_{h}+D(F_{1})+A(F_{1})+L(F_{2})$.
The vectors and matrices are indexed by the same ordering of $V(F_{0})$.
By the definitions of $M$ and $C_{F_{0}}$, $M\one=\one+2\mathbf{d}$ and $C_{F_{0}}(\mathbf{w})=\mathbf{w}^{\mathsf T}M\mathbf{w}-2\mathbf{d}^{\mathsf T}\mathbf{w}+m$.
For $\mathbf{w}\ne\mathbf{0}$, $\mathbf{w}^{\mathsf T}M\mathbf{w}=\|\mathbf{w}\|^{2}+\sum_{uv\in E(F_{1})}(w(u)+w(v))^{2}+\sum_{uv\in E(F_{2})}(w(u)-w(v))^{2}>0$, so $M$ is positive definite.

Put $\mathbf{w}_{0}=M^{-1}\mathbf{d}$.
Since $M\mathbf{w}_{0}=\mathbf{d}$ and $M$ is positive definite, we have
\begin{align}
C_{F_{0}}(\mathbf{w})
&=(\mathbf{w}_{0}-\mathbf{w})^{\mathsf T}M(\mathbf{w}_{0}-\mathbf{w})+m-\mathbf{d}^{\mathsf T}\mathbf{w}_{0}\geq m-\mathbf{d}^{\mathsf T}\mathbf{w}_{0},
\label{eq:energy-square}\\
\kappa_{F_{0}}
&=C_{F_{0}}(\mathbf{w}_{0})=m-\mathbf{d}^{\mathsf T}\mathbf{w}_{0},
\label{eq:energy-minimum}
\end{align}
with $C_{F_{0}}(\mathbf{w})=\kappa_{F_{0}}$ if and only if $\mathbf{w}=\mathbf{w}_{0}$.

Using $M\mathbf{w}_{0}=\mathbf{d}$, $M\one=\one+2\mathbf{d}$ and the symmetry of $M$, we obtain $2m=\one^{\mathsf T}\mathbf{d}=\one^{\mathsf T}M\mathbf{w}_{0}=(M\one)^{\mathsf T}\mathbf{w}_{0}=\one^{\mathsf T}\mathbf{w}_{0}+2\mathbf{d}^{\mathsf T}\mathbf{w}_{0}$.
Hence Eq.~\eqref{eq:energy-minimum} gives
\begin{equation}\label{eq:energy-sum}
\one^{\mathsf T}\mathbf{w}_{0}
=2m-2\mathbf{d}^{\mathsf T}\mathbf{w}_{0}
=2\kappa_{F_{0}}.
\end{equation}

From $M\one=\one+2\mathbf{d}$ and $\mathbf{w}_{0}=M^{-1}\mathbf{d}$, we have $M^{-1}\one=M^{-1}(M\one-2\mathbf{d})=\one-2\mathbf{w}_{0}$.
Since $F_{0}$ has an edge, $\one\ne\mathbf{0}$.
The positive definiteness of $M^{-1}$ and Eq.~\eqref{eq:energy-sum} now give
\begin{equation}\label{eq:energy-inverse}
0<\one^{\mathsf T}M^{-1}\one
=\one^{\mathsf T}(\one-2\mathbf{w}_{0})
=h-2\one^{\mathsf T}\mathbf{w}_{0}=h-4\kappa_{F_{0}}.
\end{equation}
Thus \textup{(i)} holds.

\textup{(ii)} Assume that $\alpha(G)\leq2\kappa_{F_{0}}$, and use $U,V$ and $\mathbf{t}$ as defined above.
Since $V(F_{0})\subseteq U\cup V$ and $t(v)\geq0$ for every $v\in U\cup V$, Eq.~\eqref{eq:endpoint-sums} gives $\one^{\mathsf T}\mathbf{t}\leq\sum_{v\in U\cup V}t(v)=\alpha(G)$.
Using Eq.~\eqref{eq:energy-sum}, we obtain
\begin{equation}\label{eq:energy-deficit}
0\leq2\kappa_{F_{0}}-\alpha(G)
\leq2\kappa_{F_{0}}-\one^{\mathsf T}\mathbf{t}
=\one^{\mathsf T}(\mathbf{w}_{0}-\mathbf{t}).
\end{equation}
Taking $\mathbf{w}=\mathbf{t}$ in Eq.~\eqref{eq:energy-square} and using Eq.~\eqref{eq:energy-minimum} gives
\begin{align}
    C_{F_{0}}(\mathbf{t})-\kappa_{F_{0}}
    =(\mathbf{w}_{0}-\mathbf{t})^{\mathsf T}M(\mathbf{w}_{0}-\mathbf{t}). \notag
\end{align}
Since $M$ is positive definite, the Cauchy--Schwarz inequality for the inner product $\langle\boldsymbol{\xi},\boldsymbol{\eta}\rangle_{M}=\boldsymbol{\xi}^{\mathsf T}M\boldsymbol{\eta}$ on $\mathbb R^{h}$ and Eqs.~\eqref{eq:energy-sum} and~\eqref{eq:energy-inverse} give
\begin{align}
    (2\kappa_{F_{0}}-\one^{\mathsf T}\mathbf{t})^{2}
    &=\langle M^{-1}\one,\mathbf{w}_{0}-\mathbf{t}\rangle_{M}^{2}\notag\\
    &\leq(\one^{\mathsf T}M^{-1}\one)
    (\mathbf{w}_{0}-\mathbf{t})^{\mathsf T}M(\mathbf{w}_{0}-\mathbf{t})   \notag\\
    &=(h-4\kappa_{F_{0}})
    \bigl(C_{F_{0}}(\mathbf{t})-\kappa_{F_{0}}\bigr).    \notag
\end{align}
Combining this inequality with \textup{(i)} and Eqs.~\eqref{eq:retained-energy} and~\eqref{eq:energy-deficit}, we obtain
\begin{align*}
\alpha(G)S(G)\mathrel{\geq}\mathcal E_{G}(\mathbf{z})
&\geq C_{F_{0}}(\mathbf{t})\geq\kappa_{F_{0}}
+\frac{(2\kappa_{F_{0}}-\one^{\mathsf T}\mathbf{t})^{2}}{h-4\kappa_{F_{0}}}\\
&\geq\kappa_{F_{0}}
+\frac{(2\kappa_{F_{0}}-\alpha(G))^{2}}{h-4\kappa_{F_{0}}}.\qedhere
\end{align*}
\end{proof}

\begin{remark}\label{rem:energy-components}
	\textup{(i)} Suppose that $E(F_{2})=\varnothing$.
	Let $H_{1},\ldots,H_{r}$ be pairwise vertex-disjoint subgraphs of $F_{0}$ such that $E(H_{1}\cup\cdots\cup H_{r})=E(F_{0})$.
	Then $\kappa_{F_{0}}=\sum_{j=1}^{r}\kappa_{H_{j}}$.
	For $K_{a,b}$ with $a,b\geq1$, the solution of $M\mathbf{w}_{0}=\mathbf{d}$ takes the values $\frac{b}{a+b+1}$ and $\frac{a}{a+b+1}$ on the parts of orders $a$ and $b$, respectively.
	By Eq.~\eqref{eq:energy-sum}, $\kappa_{K_{a,b}}=\frac{ab}{a+b+1}$, so $\kappa_{K_{1,k}}=\frac{k}{k+2}$, $\kappa_{K_{2,2}}=\frac45$, and $\kappa_{K_{2,3}}=1$.
	A similar calculation gives $\kappa_{P_{4}}=\frac57$.

	\textup{(ii)} For a star $F_{0}$ with center $v_{0}$, $E(F_{2})=\{v_{0}v_{1}\}$ and $m=|E(F_{1})|$, the solution of $M\mathbf{w}_{0}=\mathbf{d}$ satisfies $w_{0}(v_{0})=2w_{0}(v_{1})=\frac{m}{m+3}$ and $w_{0}(v)=\frac{3}{2(m+3)}$ for $v\in V(F_{0})\setminus\{v_{0},v_{1}\}$.
	By Eq.~\eqref{eq:energy-sum}, $\kappa_{F_{0}}=\frac{3m}{2(m+3)}$.
\end{remark}

\begin{lemma}[K\H{o}nig's theorem \cite{Konig1931}]\label{lem:konig}
	For any bipartite graph $G$, the minimum size of a vertex cover is $\nu(G)$.
\end{lemma}

\begin{lemma}\label{lem:bipartite-bound}
	Let $G$ be a graph of order $n$, and $p,q\in V(G)$ satisfy $d_{G}(p,q)=3$.
	Put $U=N_{G}(p)$, $V=N_{G}(q)$, $F=G[U,V]$, $m=|E(F)|$, $s=\nu(F)$ and $\sigma(F)=\max_{uv\in E(F)}(d_{F}(u)+d_{F}(v))$.
	Then
\begin{align}
    \alpha(\bar{G})\geq\max\left\{\frac{s+2-\sqrt{s^{2}+4s}}2,\frac{n}{n+m+\sigma(F)}\right\}. \notag
\end{align}
\end{lemma}

\begin{proof}
Since $d_{G}(p,q)=3$, we have $n\geq4$, $U\cap V=\varnothing$ and $m,s\geq1$.
Put $C=V(G)\setminus(\{p,q\}\cup U\cup V)$.
Then
\begin{equation}\label{eq:bipartite-neighborhoods}
N_{\bar{G}}(p)=\{q\}\cup V\cup C,\qquad
N_{\bar{G}}(q)=\{p\}\cup U\cup C.
\end{equation}

Put $\tau=\frac{s+2-\sqrt{s^{2}+4s}}{2}$, so $0<\tau<1$.
By Lemma~\ref{lem:konig}, $F$ has a vertex cover $W\subseteq U\cup V$ with $|W|=s$.
Thus $E_{G}(U\setminus W,V\setminus W)=\varnothing$.
Together with Eq.~\eqref{eq:bipartite-neighborhoods}, this shows that $K=\bar{G}-W$ is connected and $d(K)\leq2$.
By Lemma~\ref{lem:diameter-bound}, $\alpha(K)\geq1$.
Take any $\mathbf{y}\in\mathbb R^{n}$ with $\sum_{u\in V(K)}y(u)=0$, and put $a^{2}=\sum_{u\in V(K)}y(u)^{2}$.
By Lemma~\ref{lem:hyperplane}, it suffices to show that $\mathcal E_{\bar{G}}(\mathbf{y})\geq\tau\|\mathbf{y}\|^{2}$.
Since $1-\tau>0$, we have
\begin{align}
    (y(u)-y(q))^{2}-\tau y(u)^{2}
    &=(1-\tau)\left(y(u)-\frac{y(q)}{1-\tau}\right)^{2}-\frac{\tau}{1-\tau}y(q)^{2} \notag\\
    &\geq-\frac{\tau}{1-\tau}y(q)^{2}\qquad(u\in W\cap U), \notag\\
    (y(v)-y(p))^{2}-\tau y(v)^{2}
    &\geq-\frac{\tau}{1-\tau}y(p)^{2}\qquad(v\in W\cap V).   \notag
\end{align}
By Eq.~\eqref{eq:variational} and $\alpha(K)\geq1$, we have
\begin{align}
    \mathcal E_{K}(\mathbf{y}_{V(K)})\geq\alpha(K)a^{2}\geq a^{2}.    \notag
\end{align}
Combining the above three inequalities with $\|\mathbf{y}\|^{2}=a^{2}+\sum_{u\in W}y(u)^{2}$ (since $K=\bar{G}-W$), we obtain
\begin{align}
    \mathcal E_{\bar{G}}(\mathbf{y})-\tau\|\mathbf{y}\|^{2}
    &\geq\mathcal E_{K}(\mathbf{y}_{V(K)})-\tau a^{2}
    +\sum_{u\in W\cap U}\bigl((y(u)-y(q))^{2}-\tau y(u)^{2}\bigr)  \notag\\
    &\quad+\sum_{v\in W\cap V}\bigl((y(v)-y(p))^{2}-\tau y(v)^{2}\bigr) \notag\\
    &\geq(1-\tau)a^{2}-\frac{\tau}{1-\tau}
    \bigl(|W\cap U|y(q)^{2}+|W\cap V|y(p)^{2}\bigr).    \notag
\end{align}
Since $p,q\in V(K)$ and $|W\cap U|+|W\cap V|=s$, we have $|W\cap U|y(q)^{2}+|W\cap V|y(p)^{2}\leq sa^{2}$.
By the definition of $\tau$, $(1-\tau)^{2}=s\tau$.
Hence
\begin{align}
    \mathcal E_{\bar{G}}(\mathbf{y})-\tau\|\mathbf{y}\|^{2}
    \geq \frac{(1-\tau)^{2}-s\tau}{1-\tau}a^{2}=0. \notag
\end{align}
Thus $\alpha(\bar{G})\geq\tau$.

For the bound $\alpha(\bar{G})\geq \frac{n}{n+m+\sigma(F)}$, put $n_{1}=|U|\geq1$, $n_{2}=|V|\geq1$ and $c=|C|$.
Take any $\mathbf{y}\in\mathbb R^{n}$ with $n_{2}y(p)+n_{1}y(q)=0$.
By Lemma~\ref{lem:hyperplane}, it suffices to show that $n\|\mathbf{y}\|^{2}\leq(n+m+\sigma(F))\mathcal E_{\bar{G}}(\mathbf{y})$.
Put
\begin{equation}\label{eq:bipartite-endpoint-energy}
\begin{aligned}
E_{0}&=(y(p)-y(q))^{2}+\sum_{u\in U}(y(u)-y(q))^{2}+\sum_{v\in V}(y(v)-y(p))^{2}\\
&\quad+\sum_{u\in C}\bigl((y(u)-y(p))^{2}+(y(u)-y(q))^{2}\bigr).
\end{aligned}
\end{equation}
Let $S_{1}=\sum_{\substack{u\in U\\v\in V}}(y(u)-y(v))=n_{2}\sum_{u\in U}y(u)-n_{1}\sum_{v\in V}y(v)$ and $S_{2}=\sum_{\substack{u\in U\\v\in V}}(y(u)-y(v))^{2}$.
The Cauchy--Schwarz inequality gives $S_{1}^{2}\leq n_{1}n_{2}S_{2}$.
Combining the inequality above with Eq.~\eqref{eq:bipartite-endpoint-energy} and $n=n_{1}+n_{2}+c+2$, we obtain
\begin{align}
    E_{0}-\|\mathbf{y}\|^{2}
    &=(y(p)-y(q))^{2}-y(p)^{2}-y(q)^{2}
    \quad+\sum_{u\in U}\bigl((y(u)-y(q))^{2}-y(u)^{2}\bigr)\notag\\
    &\quad+\sum_{v\in V}\bigl((y(v)-y(p))^{2}-y(v)^{2}\bigr)  \notag\\
    &\quad+\sum_{u\in C}\bigl((y(u)-y(p))^{2}+(y(u)-y(q))^{2}-y(u)^{2}\bigr)  \notag\\
    &=n n_{1}n_{2}\theta^{2}+2\theta S_{1}
    +\sum_{u\in C}(y(u)-(n_{1}-n_{2})\theta)^{2}+c n_{1}n_{2}\theta^{2}  \notag\\
    &\geq n n_{1}n_{2}\left(\theta+\frac{S_{1}}{n n_{1}n_{2}}\right)^{2}
    -\frac{S_{1}^{2}}{n n_{1}n_{2}}\geq-\frac{S_{1}^{2}}{n n_{1}n_{2}}\geq-\frac{S_{2}}{n}, \notag
\end{align}
where $\theta=\frac{y(p)}{n_{1}}=-\frac{y(q)}{n_{2}}$.
Hence
\begin{equation}\label{eq:bipartite-allpairs-energy}
n\|\mathbf{y}\|^{2}\leq nE_{0}+S_{2}.
\end{equation}

For $uv\in E(F)$ with $u\in U$ and $v\in V$, we have $m,d_{F}(u),d_{F}(v)>0$.
The Cauchy--Schwarz inequality, together with $d_{F}(u)+d_{F}(v)\leq\sigma(F)$, gives
\begin{align}
    (y(u)-y(v))^{2}
    &=\left(\sqrt{d_{F}(u)}\,\frac{y(u)-y(q)}{\sqrt{d_{F}(u)}}+\sqrt{m}\,\frac{y(q)-y(p)}{\sqrt{m}}+\sqrt{d_{F}(v)}\,\frac{y(p)-y(v)}{\sqrt{d_{F}(v)}}\right)^{2} \notag\\
    &\leq(m+\sigma(F))
    \left(\frac{(y(u)-y(q))^{2}}{d_{F}(u)}
    +\frac{(y(q)-y(p))^{2}}{m}
    +\frac{(y(p)-y(v))^{2}}{d_{F}(v)}\right). \notag
\end{align}
Summing over $E(F)$ and using Eq.~\eqref{eq:bipartite-endpoint-energy}, we obtain
\begin{align}
    \sum_{uv\in E(F)}(y(u)-y(v))^{2}
    &\leq(m+\sigma(F))\Biggl[
    \sum_{\substack{u\in U\\d_{F}(u)>0}}(y(u)-y(q))^{2}+(y(p)-y(q))^{2} \notag\\
    &\qquad\qquad+\sum_{\substack{v\in V\\d_{F}(v)>0}}(y(p)-y(v))^{2}\Biggr]\leq(m+\sigma(F))E_{0}. \notag
\end{align}
Combining Eq.~\eqref{eq:bipartite-allpairs-energy} with the inequality above, we obtain
\allowdisplaybreaks[4]
\begin{align*}
n\|\mathbf{y}\|^{2}
&\leq nE_{0}+S_{2}=nE_{0}+\sum_{uv\in E(F)}(y(u)-y(v))^{2}
+\sum_{uv\in E_{\bar{G}}(U,V)}(y(u)-y(v))^{2}\\
&\leq nE_{0}+(m+\sigma(F))E_{0}
+\sum_{uv\in E_{\bar{G}}(U,V)}(y(u)-y(v))^{2}\\
&\leq(n+m+\sigma(F))\left(E_{0}+
\sum_{uv\in E_{\bar{G}}(U,V)}(y(u)-y(v))^{2}\right)\\
&\leq(n+m+\sigma(F))\mathcal E_{\bar{G}}(\mathbf{y}).
\end{align*}

\end{proof}

\begin{lemma}[{\cite[Theorem~2.1 and Remark~2.2]{Mathil2022}}]\label{lem:quotient}
	Let $V_{1},\ldots,V_{k}$ be a partition of $V(G)$ into nonempty sets such that $G[V_{i},V_{j}]$ is complete bipartite or empty whenever $i\ne j$.
	Put $n_{i}=|V_{i}|$ and $d_{i}=|N_{G}(v)\setminus V_{i}|$ for any $v\in V_{i}$.
	Define $Q(G)=(q_{ij})_{k\times k}$ by
\begin{align}
    q_{ij}=
    	\begin{cases}
    	  d_{i}, & i=j;\\
    	  -n_{j}, & i\ne j\text{ and }G[V_{i},V_{j}]\text{ is complete bipartite};\\
    	  0, & \text{otherwise}.
    	\end{cases} \notag
\end{align}
	Then the Laplacian spectrum of $G$ consists of the eigenvalues of $Q(G)$ together with $d_{i}+\mu_{j}(G[V_{i}])$ for $1\leq i\leq k$ and $1\leq j\leq n_{i}-1$, all counted with multiplicity.
\end{lemma}

\begin{lemma}[{\cite[Lemma~4.2]{BEHK2022}}]
\label{lem:dandelion-spectrum}
Let $H$ be a graph of order $n-3$, where $n\geq4$.
Then $\alpha(\mathcal D(H))=a_{n}$, $\mu_{1}(\mathcal D(H))=c_{n}$, where $a_{n}$ and $c_{n}$ are, respectively, the smallest and largest roots of $p_{n}(x)=x^{3}-(n+2)x^{2}+(3n-2)x-n=0$.
\end{lemma}

For the rest of the paper, put $\xi_{n}=\frac{3-\sqrt{5-\frac{4}{n-1}}}{2}$, so that $0<\xi_{n}<1$.
\begin{lemma}\label{cor:Mn}
Let $H$ be a graph of order $n-3$, where $n\geq4$.
Then
\begin{align}
    S(\mathcal D(H))=n+a_{n}-c_{n}<\min\left\{\frac{71}{50},1+\xi_{n}\right\}, \notag
\end{align}
where $a_{n}$ and $c_{n}$ are as in Lemma~\ref{lem:dandelion-spectrum}.
\end{lemma}

\begin{proof}
By Lemmas~\ref{lem:complement-spectrum} and~\ref{lem:dandelion-spectrum}, we have $S(\mathcal D(H))=n+a_{n}-c_{n}$.
The sequence $\xi_{n}$ is decreasing, and $\xi_{13}<\frac{21}{50}<\xi_{12}$.
For $4\leq n\leq12$, direct calculation gives $n+a_{n}-c_{n}<\frac{71}{50}<1+\xi_{n}$.

Now suppose that $n\geq13$.
It suffices to prove $n+a_{n}-c_{n}<1+\xi_{n}$.
Put $d=\frac{1}{(n-2)(n-3)}$.
Since $0<d\leq\frac{1}{6}$, we have
\begin{align}
    p_{n}(n-1+d)=-d^{2}\bigl((n-5)(n-2)+1-d\bigr)<0. \notag
\end{align}
Since $p_{n}(n)=n(n-3)>0$ and $c_{n}$ is the largest root of $p_{n}(x)$, we have $c_{n}>n-1+d$.

Since $0<\xi_{n}\leq\frac{1}{2}$, we have $\xi_{n}+d<1$.
Noting that $p_{n}(x)=(n-1-x)(3x-x^{2}-1)-1$ and $3\xi_{n}-\xi_{n}^{2}-1=\frac{1}{n-1}$, we obtain
\begin{align}
    p_{n}(\xi_{n}+d)
    &=d(n-1-\xi_{n}-d)(3-2\xi_{n}-d)-\frac{\xi_{n}+d}{n-1} \notag\\
    &>d(n-2)-\frac1{n-1}
    =\frac1{n-3}-\frac1{n-1}>0. \notag
\end{align}
Since $p_{n}(0)=-n<0$ and $a_{n}$ is the smallest root of $p_{n}(x)$, we have $a_{n}<\xi_{n}+d$.
Using $c_{n}>n-1+d$, we obtain $n+a_{n}-c_{n}<1+\xi_{n}$.

\end{proof}

For an injective function $\psi$ on $I$, let $\psi^{-1}\colon\psi(I)\to I$ denote its inverse, so that $\psi^{-1}(s)=\lambda$ if and only if $\psi(\lambda)=s$.

\begin{lemma}\label{lem:spectral-parameter}
	For $i=1,2$, let $\mathcal G_{i}$ be a class of $n$-vertex graphs such that each graph $G_{i}(s)\in\mathcal G_{i}$ with parameter $s\in[a,b]$ satisfies $\Phi_{i}(\mu_{k_{i}}(G_{i}(s)),s)=0$, where $a\leq b$ are nonnegative integers and $k_{i}\in\{1,\ldots,n\}$.
	Suppose there exist intervals $I_{i}$ such that $\Phi_{i}(\lambda,s)=0$ has a unique real solution $s=\psi_{i}(\lambda)$ for each $\lambda\in I_{i}$, $\mu_{k_{i}}(G_{i}(s))\in I_{i}$, and $\psi_{i}'(\lambda)>0$ on $I_{i}$.

	\begin{enumerate}
		\item[\textup{(i)}] If $\psi_{i}'(\lambda)$ is nondecreasing on $I_{i}$ for $i=1,2$, then
\begin{align}
    \mu_{k_{1}}(G_{1}(s))+\mu_{k_{2}}(G_{2}(a+b-s))
    		\geq\min_{r\in\{a,b\}}\bigl\{\mu_{k_{1}}(G_{1}(r))+\mu_{k_{2}}(G_{2}(a+b-r))\bigr\}. \notag
\end{align}
		If $\psi_{i}'(\lambda)$ is nonincreasing on $I_{i}$ for $i=1,2$, then
\begin{align}
    \mu_{k_{1}}(G_{1}(s))+\mu_{k_{2}}(G_{2}(a+b-s))
    		\leq\max_{r\in\{a,b\}}\bigl\{\mu_{k_{1}}(G_{1}(r))+\mu_{k_{2}}(G_{2}(a+b-r))\bigr\}. \notag
\end{align}

		\item[\textup{(ii)}] If $(n-\lambda)\psi_{i}'(\lambda)$ is nonincreasing on $I_{i}$ for $i=1,2$, then
\begin{align}
    &\bigl(n-\mu_{k_{1}}(G_{1}(s))\bigr)\bigl(n-\mu_{k_{2}}(G_{2}(a+b-s))\bigr) \notag\\
    &\quad\geq\min_{r\in\{a,b\}}\bigl\{\bigl(n-\mu_{k_{1}}(G_{1}(r))\bigr)\bigl(n-\mu_{k_{2}}(G_{2}(a+b-r))\bigr)\bigr\}. \notag
\end{align}
		If $(n-\lambda)\psi_{i}'(\lambda)$ is nondecreasing on $I_{i}$ for $i=1,2$, then
\begin{align}
    &\bigl(n-\mu_{k_{1}}(G_{1}(s))\bigr)\bigl(n-\mu_{k_{2}}(G_{2}(a+b-s))\bigr) \notag\\
    &\quad\leq\max_{r\in\{a,b\}}\bigl\{\bigl(n-\mu_{k_{1}}(G_{1}(r))\bigr)\bigl(n-\mu_{k_{2}}(G_{2}(a+b-r))\bigr)\bigr\}. \notag
\end{align}
	\end{enumerate}

	In either part, if the monotonicity is strict for at least one $i$, then the corresponding inequality is strict for $s\in(a,b)$.
\end{lemma}

\begin{proof}
	For $i=1,2$ and each $G_{i}(s)\in\mathcal G_{i}$, $\mu_{k_{i}}(G_{i}(s))=\psi_{i}^{-1}(s)$, since $\Phi_{i}(\lambda,s)=0$ gives $s=\psi_{i}(\lambda)$ for $\lambda=\mu_{k_{i}}(G_{i}(s))\in I_{i}$.

	\textup{(i)} Since $\psi_{i}'(\lambda)>0$ on $I_{i}$ for $i=1,2$, each $\psi_{i}^{-1}(s)$ is strictly increasing. Put
\begin{align}
    F(x)=\psi_{1}^{-1}(x)+\psi_{2}^{-1}(a+b-x)\qquad (x\in[a,b]). \notag
\end{align}
	Then
\begin{align}
    F'(x)=\frac{1}{\psi_{1}'(\psi_{1}^{-1}(x))}
    	-\frac{1}{\psi_{2}'(\psi_{2}^{-1}(a+b-x))}
    	\qquad (x\in(a,b)). \notag
\end{align}
	If both $\psi_{i}'(\lambda)$ are nondecreasing on $I_{i}$, then $F'(x)$ is nonincreasing.
	For $a<s<b$, if $F'(s)\geq0$, then $F'(x)\geq0$ for $a<x<s$, so $F(s)\geq F(a)$, with strict inequality if $F'(x)$ is strictly decreasing.
	If $F'(s)<0$, then $F'(x)<0$ for $s<x<b$, so $F(s)>F(b)$.
	Thus $F(s)\geq\min\{F(a),F(b)\}$.
	If both $\psi_{i}'(\lambda)$ are nonincreasing on $I_{i}$, then $F'(x)$ is nondecreasing, and the same argument applied to $-F$ gives $F(s)\leq\max\{F(a),F(b)\}$.

	\textup{(ii)} Put
\begin{align}
    P(x)=(n-\psi_{1}^{-1}(x))(n-\psi_{2}^{-1}(a+b-x))\qquad (x\in[a,b]). \notag
\end{align}
	Put $h_{i}(\lambda)=(n-\lambda)\psi_{i}'(\lambda)$ for $\lambda\in I_{i}$ and $i=1,2$.
	Then
\begin{align}
    P'(x)=
    	\frac{h_{1}(\psi_{1}^{-1}(x))-h_{2}(\psi_{2}^{-1}(a+b-x))}
    	{\psi_{1}'(\psi_{1}^{-1}(x))\psi_{2}'(\psi_{2}^{-1}(a+b-x))}
    	\qquad (x\in(a,b)). \notag
\end{align}
	The denominator is positive since $\psi_{i}'(\lambda)>0$ on $I_{i}$ for $i=1,2$.
	If both $h_{i}(\lambda)$ are nonincreasing on $I_{i}$, then the numerator is nonincreasing in $x$, so the sign argument in \textup{(i)} gives $P(s)\geq\min\{P(a),P(b)\}$.
	If both $h_{i}(\lambda)$ are nondecreasing on $I_{i}$, then the sign argument applied to $-P$ gives $P(s)\leq\max\{P(a),P(b)\}$.

	Strict monotonicity for at least one $i$ makes $F'(x)$ and the numerator of $P'(x)$ strictly monotone in \textup{(i)} and \textup{(ii)}, respectively.
\end{proof}

\begin{figure}[H]
\centering
\setlength{\abovecaptionskip}{4pt}
\setlength{\belowcaptionskip}{0pt}
\tikzset{familyvertex/.style={circle,fill=black,inner sep=1.35pt,outer sep=.8pt},familyblock/.style={ellipse,draw,minimum width=6.4mm,minimum height=12mm,inner sep=0pt}}
\begin{minipage}[b]{.32\linewidth}
\centering
\resizebox{\linewidth}{!}{%
\begin{tikzpicture}[font=\small,line width=.45pt,every label/.style={font=\scriptsize}]
\path[use as bounding box] (-.65,-.2) rectangle (4.60,2.1);
\node[familyblock] (partone) at (0,.85) {};
\node[familyblock] (parttwo) at (1.60,.85) {};
\node[font=\scriptsize,anchor=south,inner sep=0pt] at (0,1.55) {$|V_{1}|=m_{1}$};
\node[font=\scriptsize,anchor=south,inner sep=0pt] at (1.60,1.55) {$|V_{2}|=m_{2}$};
\node[font=\scriptsize,anchor=north,inner sep=0pt] at (2.90,.60) {$|V_{3}|=1$};
\node[font=\scriptsize,anchor=north,inner sep=0pt] at (4.08,.60) {$|V_{4}|=1$};
\node[font=\scriptsize] at (0,.85) {$\vdots$};
\node[font=\scriptsize] at (1.60,.85) {$\vdots$};
\node[familyvertex] (a1) at (0,1.19) {};
\node[familyvertex] (a2) at (0,.51) {};
\node[familyvertex] (b1) at (1.60,1.19) {};
\node[familyvertex] (b2) at (1.60,.51) {};
\node[familyvertex] (c) at (2.90,.85) {};
\node[familyvertex] (d) at (4.08,.85) {};
\foreach \a in {a1,a2} {\foreach \b in {b1,b2} {\draw (\a)--(\b);}}
\draw (b1)--(c)--(d);
\draw (b2)--(c);
\end{tikzpicture}%
}\par
{\small (a) $G\in\mathcal H(m_{1},m_{2},1,1)$}
\end{minipage}\hfill
\begin{minipage}[b]{.32\linewidth}
\centering
\resizebox{\linewidth}{!}{%
\begin{tikzpicture}[font=\small,line width=.45pt,every label/.style={font=\scriptsize}]
\path[use as bounding box] (-.80,-.2) rectangle (4.40,2.1);
\node[familyblock] (partone) at (0,.85) {};
\node[familyblock] (parttwo) at (2.45,.85) {};
\node[font=\scriptsize,anchor=south,inner sep=0pt] at (0,1.55) {$|V_{1}|=m_{1}$};
\node[font=\scriptsize,anchor=north,inner sep=0pt] at (1.22,.60) {$|V_{2}|=1$};
\node[font=\scriptsize,anchor=south,inner sep=0pt] at (2.45,1.55) {$|V_{3}|=m_{2}$};
\node[font=\scriptsize,anchor=north,inner sep=0pt] at (3.60,.60) {$|V_{4}|=1$};
\node[font=\scriptsize] at (0,.85) {$\vdots$};
\node[font=\scriptsize] at (2.45,.85) {$\vdots$};
\node[familyvertex] (a1) at (0,1.19) {};
\node[familyvertex] (a2) at (0,.51) {};
\node[familyvertex] (b) at (1.22,.85) {};
\node[familyvertex] (c1) at (2.45,1.19) {};
\node[familyvertex] (c2) at (2.45,.51) {};
\node[familyvertex] (d) at (3.60,.85) {};
\draw (a1)--(b)--(c1)--(d);
\draw (a2)--(b)--(c2)--(d);
\end{tikzpicture}%
}\par
{\small (b) $G\in\mathcal H(m_{1},1,m_{2},1)$}
\end{minipage}\hfill
\begin{minipage}[b]{.32\linewidth}
\centering
\resizebox{\linewidth}{!}{%
\begin{tikzpicture}[font=\small,line width=.45pt]
\path[use as bounding box] (-.80,-.2) rectangle (4.40,2.1);
\node[familyblock,minimum height=10mm] (zblock) at (1.8,1.14) {};
\node[font=\scriptsize,anchor=south,inner sep=0pt] at ($(zblock.north)+(0,.08)$) {$|F|=n-4$};
\node[font=\scriptsize] at (1.8,1.14) {$\vdots$};
\node[familyvertex] (z1) at (1.8,1.41) {};
\node[familyvertex] (z2) at (1.8,.87) {};
\node[familyvertex,label=below:$p$] (p) at (0,.15) {};
\node[familyvertex,label={below:$u$}] (u) at (1.15,.15) {};
\node[familyvertex,label={below:$v$}] (v) at (2.45,.15) {};
\node[familyvertex,label=below:$q$] (q) at (3.6,.15) {};
\draw (p)--(u)--(v)--(q);
\foreach \z in {z1,z2} {\draw (u)--(\z)--(v);}
\end{tikzpicture}%
}\par
{\small (c) $G\in\mathcal H_{n}$}
\end{minipage}
\caption{The graph families $\mathcal H(m_{1},m_{2},1,1)$, $\mathcal H(m_{1},1,m_{2},1)$ and $\mathcal H_{n}$.}
\label{fig:four-block-bull}
\end{figure}
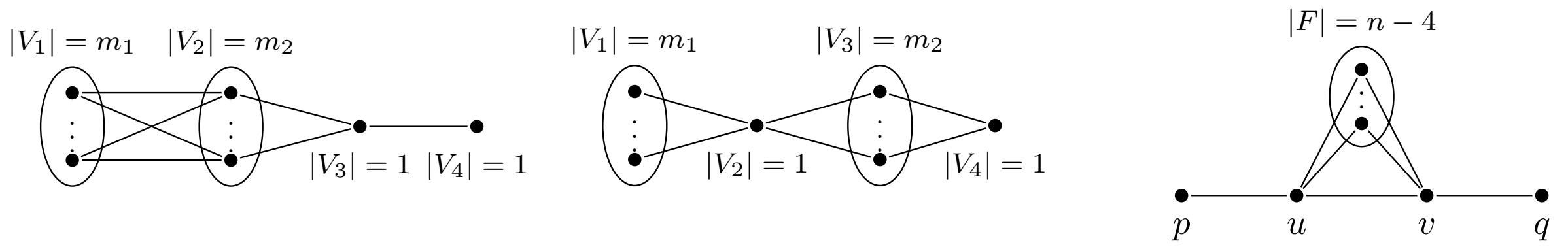

For positive integers $h_{1},\ldots,h_{4}$, let $\mathcal H(h_{1},h_{2},h_{3},h_{4})$ be the family of graphs $G$ whose vertex set has a partition $V_{1},\ldots,V_{4}$ with $|V_{i}|=h_{i}$ for $1\leq i\leq4$.
For $1\leq i<j\leq4$, $G[V_{i},V_{j}]$ is complete bipartite if $j=i+1$ and empty otherwise.
Each $G[V_{i}]$ is an arbitrary graph.
See Figure~\ref{fig:four-block-bull}(a,b) for the two special cases used below.
For $n\geq4$, it is easy to see that $\mathcal H(n-3,1,1,1)=\{\mathcal D(H):|H|=n-3\}$.

For each integer $n\geq5$, let $\mathcal H_{n}$ be the family of graphs obtained from the disjoint union of a path $puvq$ and an arbitrary graph $F$ of order $n-4$ by adding the edges $uw,vw$ for all $w\in V(F)$ (see Figure~\ref{fig:four-block-bull}(c)).
Put $\mathcal H_{4}=\varnothing$.

\begin{lemma}\label{lem:four-block}
Let $m_{1},m_{2}$ be positive integers and $n=m_{1}+m_{2}+2$.
For any $G\in\mathcal H(m_{1},m_{2},1,1)\cup\mathcal H(m_{1},1,m_{2},1)$, we have $S(G)\geq n+a_{n}-c_{n}$, with equality if and only if $m_{1}=1$ or $m_{2}=1$.
\end{lemma}

\begin{proof}

First let $G\in\mathcal H(m_{1},m_{2},1,1)$.
Then $\bar{G}\in\mathcal H(m_{2},1,m_{1},1)$.
Let $Q(G)$ and $Q(\bar{G})$ be the corresponding matrices from Lemma~\ref{lem:quotient}.
A direct calculation gives
\begin{align}
    \det(xI_{4}-Q(G))=x\Phi(x,m_{2}),\qquad
    \det(xI_{4}-Q(\bar{G}))=x\Phi(x,m_{1}), \notag
\end{align}
where $\Phi(x,s)=x^{3}-(n+1)x^{2}+(2n-2)x-s(x^{2}-nx+n)$.
Since $\Phi(0,m_{i})=-nm_{i}<0<\Phi(1,m_{i})=n-2-m_{i}$ for $i=1,2$, Lemma~\ref{lem:quotient} gives $0<\alpha(G),\alpha(\bar{G})<1$ and $\Phi(\alpha(G),m_{2})=\Phi(\alpha(\bar{G}),m_{1})=0$.

Since $x^{2}-nx+n\geq1$ on $[0,1]$, the equation $\Phi(x,s)=0$ has the unique real solution
\begin{align}
    s=\psi(x):=\frac{x^{3}-(n+1)x^{2}+(2n-2)x}{x^{2}-nx+n}
    =x-1+\frac{n-2x}{x^{2}-nx+n}. \notag
\end{align}
Hence
\begin{equation}\label{eq:comparison-relations}
\psi(\alpha(G))=m_{2},\qquad
\psi(\alpha(\bar{G}))=m_{1}.
\end{equation}

For $0\leq x\leq1$, direct calculation gives
\begin{align}
    \psi''(x)=\frac{2(n-2x)(n^{2}-nx-3n+x^{2})}{(x^{2}-nx+n)^{3}}>0,
    \qquad
    \psi'(0)=2-\frac{2}{n}>0. \notag
\end{align}
Hence $\psi'(x)>0$ and $\psi'(x)$ is strictly increasing on $(0,1)$.

Take $\mathcal G_{1}=\bigcup_{s=1}^{n-3}\mathcal H(n-2-s,s,1,1)$ and $\mathcal G_{2}=\bigcup_{s=1}^{n-3}\mathcal H(n-2-s,1,s,1)$.
By Eq.~\eqref{eq:comparison-relations} and Lemma~\ref{cor:Mn}, Lemma~\ref{lem:spectral-parameter}\textup{(i)} gives $S(G)\geq n+a_{n}-c_{n}$, with equality exactly when $m_{1}=1$ or $m_{2}=1$.

Now let $G\in\mathcal H(m_{1},1,m_{2},1)$.
Then $\bar{G}\in\mathcal H(m_{2},m_{1},1,1)$.
Since $S(G)=S(\bar{G})$, the result holds.
\end{proof}

\begin{corollary}\label{cor:bull}
Let $G\in\mathcal H_{n}$ ($n\geq5$).
Then $S(G)>n+a_{n}-c_{n}$.
\end{corollary}

\begin{proof}
Let $Q(G)$ be the matrix defined in Lemma~\ref{lem:quotient} for the partition $\{p\},\{u\},V(F),\{v\},\{q\}$.
A direct calculation gives
\begin{align}
    \det(xI_{5}-Q(G))=x(x^{2}-nx+n)(x^{2}-nx+n-2). \notag
\end{align}
Since $N_{G}(w)\setminus V(F)=\{u,v\}$ for every $w\in V(F)$, Lemma~\ref{lem:quotient} shows that $\alpha(G)\in(0,1)$ is the smaller root of $x^{2}-nx+n-2=0$.
Since $\bar{G}\in\mathcal H_{n}$ with path $vpqu$ and graph $\bar{F}$, the same calculation for $Q(\bar{G})$ gives $\alpha(\bar{G})=\alpha(G)$.
For $\lambda=\alpha(G)=\alpha(\bar{G})$, the formula for $\psi(\lambda)$ in the proof of Lemma~\ref{lem:four-block} and $\lambda^{2}-n\lambda+n=2$ give $s=\psi(\lambda)=\frac{n-2}{2}\in(1,n-3)$.
By Lemma~\ref{lem:spectral-parameter}\textup{(i)} and the proof of Lemma~\ref{lem:four-block}, we obtain $S(G)>n+a_{n}-c_{n}$.
\end{proof}

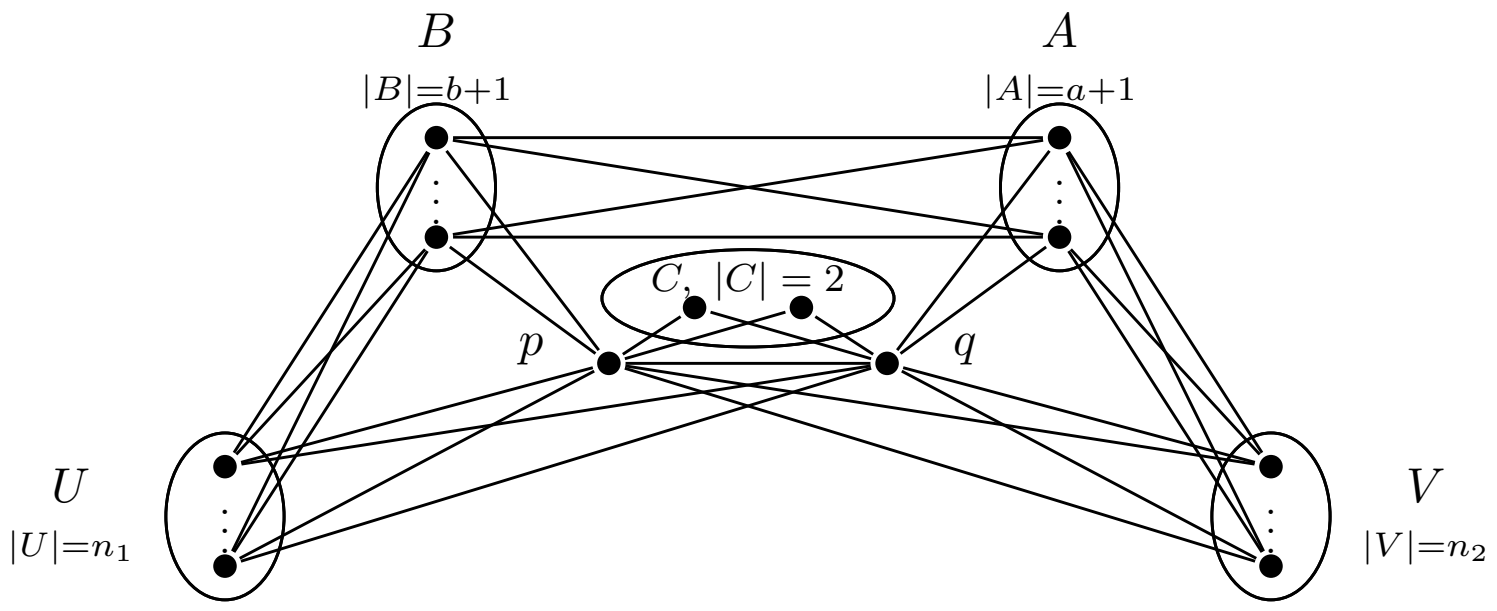
\begin{figure}[!htbp]
\centering
\begin{tikzpicture}[papergraph,x=.8cm,y=.55cm,scale=.92,transform shape]
\begin{pgfinterruptboundingbox}
\node[vertex] (p) at (-1.25,0) {};
\node[inner sep=0pt] at (-1.95,.18) {$p$};
\node[vertex] (q) at (1.25,0) {};
\node[inner sep=0pt] at (1.95,.18) {$q$};
\node[block,minimum width=8.5mm,minimum height=12mm,inner sep=0pt] (B) at (-2.8,2.3) {};
\node[vertex] (B1) at (-2.8,1.65) {};
\node[vertex] (B2) at (-2.8,2.9499999999999997) {};
\node[font=\scriptsize,inner sep=0pt] at (-2.8,2.3) {$\vdots$};
\node[align=center] at (-2.8,3.95) {$B$\\[-1pt]$\scriptstyle |B|=b+1$};
\node[block,minimum width=8.5mm,minimum height=12mm,inner sep=0pt] (A) at (2.8,2.3) {};
\node[vertex] (A1) at (2.8,1.65) {};
\node[vertex] (A2) at (2.8,2.9499999999999997) {};
\node[font=\scriptsize,inner sep=0pt] at (2.8,2.3) {$\vdots$};
\node[align=center] at (2.8,3.95) {$A$\\[-1pt]$\scriptstyle |A|=a+1$};
\node[block,minimum width=8.5mm,minimum height=12mm,inner sep=0pt] (U) at (-4.7,-2.0) {};
\node[vertex] (U1) at (-4.7,-2.65) {};
\node[vertex] (U2) at (-4.7,-1.35) {};
\node[font=\scriptsize,inner sep=0pt] at (-4.7,-2.0) {$\vdots$};
\node[anchor=east,align=center] at ($(U.west)+(-.10,0)$) {$U$\\[-1pt]$\scriptstyle |U|=n_{1}$};
\node[block,minimum width=8.5mm,minimum height=12mm,inner sep=0pt] (V) at (4.7,-2.0) {};
\node[vertex] (V1) at (4.7,-2.65) {};
\node[vertex] (V2) at (4.7,-1.35) {};
\node[font=\scriptsize,inner sep=0pt] at (4.7,-2.0) {$\vdots$};
\node[anchor=west,align=center] at ($(V.east)+(.10,0)$) {$V$\\[-1pt]$\scriptstyle |V|=n_{2}$};
\node[block,minimum width=21mm,minimum height=7mm,inner sep=0pt] (C) at (0,.85) {};
\node[vertex] (C1) at (-.48,.73) {};
\node[vertex] (C2) at (.48,.73) {};
\node[font=\scriptsize,inner sep=0pt] at (0,1.08) {$C,\ |C|=2$};
\draw (p)--(q);
\draw (p)--(B1);
\draw (p)--(B2);
\draw (p)--(C1);
\draw (p)--(C2);
\draw (p)--(U1);
\draw (p)--(U2);
\draw (q)--(A1);
\draw (q)--(A2);
\draw (q)--(C1);
\draw (q)--(C2);
\draw (q)--(V1);
\draw (q)--(V2);
\draw (U1)--(B1);
\draw (U1)--(B2);
\draw (U2)--(B1);
\draw (U2)--(B2);
\draw (V1)--(A1);
\draw (V1)--(A2);
\draw (V2)--(A1);
\draw (V2)--(A2);
\draw (B2)--(A2);
\draw (B1)--(A1);
\draw (B2)--(A1);
\draw (B1)--(A2);
\draw (p)--(V1);
\draw (p)--(V2);
\draw (q)--(U2);
\draw (q)--(U1);
\draw (-2.8,2.3) ellipse [x radius=4.25mm,y radius=6mm];
\draw (2.8,2.3) ellipse [x radius=4.25mm,y radius=6mm];
\draw (-4.7,-2.0) ellipse [x radius=4.25mm,y radius=6mm];
\draw (4.7,-2.0) ellipse [x radius=4.25mm,y radius=6mm];
\draw (0,.85) ellipse [x radius=10.5mm,y radius=3.5mm];
\foreach \v in {p,q,B1,B2,A1,A2,U1,U2,V1,V2,C1,C2} {\fill (\v) circle[radius=1.7pt];}
\end{pgfinterruptboundingbox}
\path[use as bounding box] (-6.5699,-3.2304) rectangle (6.5485,4.7936);
\end{tikzpicture}
\caption{ $F(a,b,2,n_{1},n_{2})$.}
\label{fig:Fk}
\end{figure}

\begin{lemma}\label{lem:Fk}

For nonnegative integers $a,b,c,n_{1},n_{2}$, let $F(a,b,c,n_{1},n_{2})$ be obtained from pairwise disjoint independent vertex sets $\{p\},\{q\},A,B,C,U,V$ by adding all possible edges between the following pairs of sets: $(\{p\},\{q\})$, $(\{p\},B)$, $(\{q\},A)$, $(\{p,q\},C\cup U\cup V)$, $(A,B)$, $(B,U)$ and $(A,V)$, where $|A|=a+1$, $|B|=b+1$, $|C|=c\in\{1,2\}$, $|U|=n_{1}$ and $|V|=n_{2}$(See Figure~\ref{fig:Fk} for $c=2$).
Then $\alpha(F(a,b,c,n_{1},n_{2}))>\frac{5}{4}$.

\end{lemma}

\begin{proof}
	Put $F=F(a,b,c,n_{1},n_{2})$ and suppose that $\alpha(F)\leq\frac{5}{4}$.
	Let $\mathbf{x}$ be a Fiedler vector of $F$, and choose $u_{0},v_{0}\in V(F)$ such that
	$x_{u_{0}}=\max_{v\in V(F)}x_{v}$ and $x_{v_{0}}=\min_{v\in V(F)}x_{v}$.
	Put $\mathbf{f}=\frac{\mathbf{x}-x_{v_{0}}\one}{x_{u_{0}}-x_{v_{0}}}$.
	Then $f_{u_{0}}=1$, $f_{v_{0}}=0$ and $0\leq f_{v}\leq1$ for every $v\in V(F)$.
	By Lemma~\ref{lem:normalization}\textup{(i)}, we have $u_{0}v_{0}\notin E(F)$ and $|N_{F}(u_{0})\cap N_{F}(v_{0})|\leq1$.
	
	By the definition of $F$, one of $u_{0},v_{0}$ lies in $C$ and the other in $A\cup B$.
	By replacing $\mathbf{f}$ by $\one-\mathbf{f}$ and interchanging $u_{0},v_{0}$ if necessary, assume that $u_{0}\in C$.
By interchanging $(p,A,U)$ with $(q,B,V)$ if necessary, assume that $v_{0}\in A$.
	Then $N_{F}(u_{0})=\{p,q\}$ and $N_{F}(v_{0})=\{q\}\cup B\cup V$.
	Eq.~\eqref{eq:extreme-sum} gives $\alpha(F)=2-f_{p}+\sum_{v\in B\cup V}f_{v}\geq2-f_{p}$, so $f_{p}\geq2-\alpha(F)\geq\frac{3}{4}$.
	
	For $u\in U$, Eq.~\eqref{eq:vertex-equations} and $f_{v}\geq0$ for every $v\in V$ give
\begin{align}
    (L(F)\mathbf{f})_{u}-(L(F)\mathbf{f})_{v_{0}}=\alpha(F)f_{u}
    	&=(b+3)f_{u}-f_{p}+\sum_{v\in V}f_{v}, \notag\\
    	(b+3-\alpha(F))f_{u}&=f_{p}-\sum_{v\in V}f_{v}\leq f_{p}. \notag
\end{align}
	Since $b+3-\alpha(F)>1$, we have $f_{u}\leq f_{p}$ for every $u\in U$.
	Then Eq.~\eqref{eq:vertex-equations} and $f_{v}\leq1$ for $v\in C$ give
\begin{align}
    (L(F)\mathbf{f})_{p}-(L(F)\mathbf{f})_{v_{0}}=\alpha(F)f_{p}
    	&=(b+c+n_{1}+n_{2}+2)f_{p}-\sum_{v\in C\cup U}f_{v}, \notag\\
    	(b+c+n_{1}+n_{2}+2-\alpha(F))f_{p}&=\sum_{v\in C\cup U}f_{v}
    	\leq c+n_{1}f_{p}. \notag
\end{align}
	Combining the above inequality with $f_{p}\geq\frac{3}{4}$ and $\alpha(F)\leq\frac{5}{4}$, we obtain
\begin{align}
    c\geq(b+c+n_{2}+2-\alpha(F))f_{p}
    	\geq(c+2-\alpha(F))f_{p}
    	\geq\frac34\left(c+\frac34\right). \notag
\end{align}
	Thus $c\geq\frac{9}{4}$, contrary to $c\in\{1,2\}$.
\end{proof}

\section{Auxiliary lemmas and proof of Theorem~\ref{thm:main}}

For $u,v\in V(G)$, let $P_{uv}(G)$ denote a shortest $u$--$v$ path in $G$.

\begin{lemma}\label{lem:diameter-two}
Let $G$ and $\bar{G}$ be connected graphs of order $n\geq4$.
Let $\mathbf{x},\mathbf{y}$ be Fiedler vectors of $G,\bar{G}$, respectively, and $p,q\in V(G)$ satisfy $x_{p}=\max_{v\in V(G)}x_{v}$ and $x_{q}=\min_{v\in V(G)}x_{v}$.
Then the following hold.
\begin{enumerate}
\item[\textup{(i)}] If $d_{G}(p,q)\leq2$ and $x_{p}-x_{q}\geq\max_{u,v\in V(G)}|y_{u}-y_{v}|$, then $S(G)>n+a_{n}-c_{n}$.
\item[\textup{(ii)}] If $d(G)=2$ or $d(\bar{G})=2$, then $S(G)>n+a_{n}-c_{n}$.
\end{enumerate}
\end{lemma}

\Needspace{6\baselineskip}
\begin{proof}
By Lemma~\ref{cor:Mn}, we may assume $S(G)<\frac{71}{50}$.

\textup{(i)} Put $\mathbf{z}=\frac{\mathbf{x}-x_{q}\one}{x_{p}-x_{q}}$.
Then $z_{p}=1$, $z_{q}=0$ and $0\leq z_{v}\leq1$ for every $v\in V(G)$.
By Lemma~\ref{lem:normalization}\textup{(ii)},
\begin{equation}\label{eq:diameter-two-bound}
\mathcal E_{G}(\mathbf{z})\leq\alpha(G)S(G)<\frac{71}{50}\alpha(G).
\end{equation}
Since $\alpha(G)<\frac{71}{50}$, Lemma~\ref{lem:normalization}\textup{(i)} gives $pq\notin E(G)$ and $|N_{G}(p)\cap N_{G}(q)|\leq1$.
Together with $d_{G}(p,q)\leq2$, this gives $d_{G}(p,q)=2$.
Write $N_{G}(p)\cap N_{G}(q)=\{r\}$.
Let $A=N_{G}(p)\setminus\{r\}$, $B=N_{G}(q)\setminus\{r\}$, $F=G[A,B]$, $m=|E(F)|$, and $G'$ be a spanning subgraph of $\bar{G}$ obtained by deleting $d_{\bar{G}}(r)-1$ edges incident with $r$.

If $m=0$, then any two nonadjacent vertices of $\bar{G}-r$ have a common neighbor in $\{p,q\}$.
Thus $\bar{G}-r$ is connected and $d(\bar{G}-r)\leq2$.
By Lemma~\ref{lem:diameter-bound}, $\alpha(\bar{G}-r)\geq1$.
By Lemmas~\ref{lem:subgraph} and~\ref{lem:pendant} with $c=1$, we obtain $\alpha(\bar{G})\geq\alpha(G')\geq\xi_{n}$.
By Lemma~\ref{lem:normalization}\textup{(i)}, $\alpha(G)\geq|N_{G}(p)\cap N_{G}(q)|=1$.
Hence, by Lemma~\ref{cor:Mn}, $S(G)\geq1+\xi_{n}>n+a_{n}-c_{n}$.

For $m\geq1$, Eq.~\eqref{eq:extreme-sum} gives
\begin{align}
    \frac{71}{50}>\alpha(G)=1+\sum_{v\in A}(1-z_{v})+\sum_{v\in B}z_{v}\eqqcolon 1+\varepsilon. \notag
\end{align}
For any $ab\in E(F)$ with $a\in A$ and $b\in B$, Eq.~\eqref{eq:vertex-equations} and $0\leq z_{v}\leq1$ for all $v\in V(G)$ give
\begin{align}
    \alpha(G)(z_{a}-z_{b})
    &=(L(G)\mathbf{z})_{a}-(L(G)\mathbf{z})_{b} \notag\\
    &=2(z_{a}-z_{b})+\sum_{v\in N_{G}(a)\setminus\{b\}}(z_{a}-z_{v})
    +\sum_{v\in N_{G}(b)\setminus\{a\}}(z_{v}-z_{b}) \notag\\
    &\geq2(z_{a}-z_{b})-(n-2)\bigl((1-z_{a})+z_{b}\bigr). \notag
\end{align}
By the definition of $\varepsilon$, $(1-z_{a})+z_{b}\leq\varepsilon$.
Since $\varepsilon<1$, this gives $z_{a}-z_{b}\geq1-\varepsilon>0$.

Combining the above three inequalities with $2-\alpha(G)=1-\varepsilon$, we obtain
\begin{align}
    (1-\varepsilon)^{2}
    \leq(2-\alpha(G))(z_{a}-z_{b})
    \leq(n-2)\varepsilon. \notag
\end{align}
Thus $\varepsilon>0$ and $n\varepsilon\geq1+\varepsilon^{2}>1$. Hence
\begin{equation}\label{eq:gap}
\alpha(G)=1+\varepsilon>1+\frac1n.
\end{equation}

Since $N_{G}(p)\cap N_{G}(q)=\{r\}$ and $m\geq1$, we have $d_{G-r}(p,q)=3$.
By Lemma~\ref{lem:bipartite-bound} and $\sigma(F)=\max_{ab\in E(F)}(d_{F}(a)+d_{F}(b))\leq m+1$, we have $\alpha(\bar{G}-r)\geq\frac{n-1}{n+2m}>0$.
By Lemmas~\ref{lem:subgraph} and~\ref{lem:pendant} with $c=\frac{n-1}{n+2m}$, we obtain $\alpha(\bar{G})\geq\alpha(G')\geq\theta_{m}$, where
\begin{align}
    \theta_{m}=\frac{3n+4m-1-\sqrt{(3n+4m-1)^{2}-4n(n+2m)}}{2(n+2m)}. \notag
\end{align}

For $m=1$, we have $\alpha(\bar{G})\geq\theta_{1}>\xi_{n}-\frac1n$.
Eq.~\eqref{eq:gap} and Lemma~\ref{cor:Mn} give $S(G)>1+\xi_{n}>n+a_{n}-c_{n}$.

For $m\geq2$, we have $\alpha(\bar{G})\geq\theta_{m}>\frac1{m+2}$.
Thus $S(G)<\frac{71}{50}$ gives $\alpha(G)<\frac{71}{50}-\frac{1}{m+2}<2$.
Since $z_{a}-z_{b}\geq1-\varepsilon=2-\alpha(G)>0$ for every $ab\in E(F)$, Eq.~\eqref{eq:diameter-two-bound} gives
\postdisplaypenalty=10000\relax
\begin{align}
    \frac{71}{50}\alpha(G)
    >\alpha(G)S(G)\geq\mathcal E_{G}(\mathbf{z})
    &\geq(1-z_{r})^{2}+z_{r}^{2}+\sum_{ab\in E(F)}(z_{a}-z_{b})^{2} \notag\\
    &\geq\frac12+m(2-\alpha(G))^{2}
    >\frac12+m\left(\frac{29}{50}+\frac1{m+2}\right)^{2} \notag\\
    &>\frac{71}{50}\left(\frac{71}{50}-\frac1{m+2}\right)
    >\frac{71}{50}\alpha(G), \notag
\end{align}
a contradiction.

\textup{(ii)} Without loss of generality, we need only consider the case $d(\bar{G})=2$.
Choose $a_{0},b_{0}\in V(G)$ with $y_{a_{0}}=\max_{v\in V(G)}y_{v}$ and $y_{b_{0}}=\min_{v\in V(G)}y_{v}$.

If $x_{p}-x_{q}<y_{a_{0}}-y_{b_{0}}$, then the result follows by applying \textup{(i)} to $\bar{G}$.

If $x_{p}-x_{q}\geq y_{a_{0}}-y_{b_{0}}$, then by \textup{(i)} we may assume $d_{G}(p,q)\geq3$.
Since $\alpha(\bar{G})\leq S(G)<\frac{71}{50}$, Lemma~\ref{lem:normalization}\textup{(i)} gives $a_{0}b_{0}\notin E(\bar{G})$ and $|N_{\bar{G}}(a_{0})\cap N_{\bar{G}}(b_{0})|\leq1$.
Since $d(\bar{G})=2$, we have $d_{\bar{G}}(a_{0},b_{0})=2$, and there exists a vertex $r_{0}$ such that $N_{\bar{G}}(a_{0})\cap N_{\bar{G}}(b_{0})=\{r_{0}\}$.
Put
\begin{align}
    F^\ast=\bar{G}\bigl[N_{\bar{G}}(a_{0})\setminus\{r_{0}\},\,N_{\bar{G}}(b_{0})\setminus\{r_{0}\}\bigr]. \notag
\end{align}
If $E(F^\ast)=\varnothing$, then an argument similar to that for $m=0$ in \textup{(i)} gives $S(G)\geq1+\xi_{n}>n+a_{n}-c_{n}$.
If $E(F^\ast)\ne\varnothing$, then a calculation similar to that leading to Eq.~\eqref{eq:gap} gives $\alpha(\bar{G})>1+\frac{1}{n}$.

If $\alpha(G)\geq\frac{3-\sqrt5}{2}$, then by Lemma~\ref{cor:Mn}, we have $S(G)>1+\frac1n+\frac{3-\sqrt5}{2}>1+\xi_{n}>n+a_{n}-c_{n}$.

Thus we may assume $\alpha(G)<\frac{3-\sqrt5}{2}$.
Since $V(G)\setminus(N_{G}[a_{0}]\cup N_{G}[b_{0}])=N_{\bar{G}}(a_{0})\cap N_{\bar{G}}(b_{0})=\{r_{0}\}$ and $a_{0}b_{0}\in E(G)$, we have $d(G)\leq4$ and $d_{G}(p,q)\in\{3,4\}$.
Put $\mathbf{z}=\frac{\mathbf{x}-x_{q}\one}{x_{p}-x_{q}}$.
Then $z_{p}=1$, $z_{q}=0$, and $0\leq z_{v}\leq1$ for each $v\in V(G)$.
Since $x_{p}-x_{q}\geq y_{a_{0}}-y_{b_{0}}$, Lemma~\ref{lem:normalization}\textup{(ii)} gives Eq.~\eqref{eq:diameter-two-bound}.
Since $|E(P_{pq}(G))|\leq4$ and $z_{p}-z_{q}=1$, the Cauchy--Schwarz inequality and Eq.~\eqref{eq:diameter-two-bound} give
\begin{align}
    \frac14\leq\frac{(z_{p}-z_{q})^{2}}{|E(P_{pq}(G))|}
    \leq\sum_{uv\in E(P_{pq}(G))}(z_{u}-z_{v})^{2}
    \leq\mathcal E_{G}(\mathbf{z})<\frac{71}{50}\alpha(G). \notag
\end{align}
Hence $\alpha(G)>\frac{25}{142}$, and $S(G)<\frac{71}{50}$ gives $\alpha(\bar{G})<\frac{71}{50}-\frac{25}{142}<\frac{5}{4}$.

First suppose that $d_{G}(p,q)=3$, and write $P_{pq}(G)=puvq$.
If $uv$ is a cut edge, let $X,Y$ be the vertex sets of the components of $G-uv$, with $u\in X$ and $v\in Y$.
Since $d(\bar{G})=2$, choose $v_{0}\in N_{\bar{G}}(u)\cap N_{\bar{G}}(v)$.
Without loss of generality, assume that $v_{0}\in X$.
Then $v_{0}\ne p$, $|X|\geq3$ and $|Y|\geq2$.
Since $E_{G}(X,Y)=\{uv\}$, $\bar{G}[X\setminus\{u\},Y]$ is complete bipartite and $uq\in E(\bar{G})$.
Let $H$ be the spanning subgraph of $\bar{G}$ obtained from $\bar{G}[X\setminus\{u\},Y]$ by adding $u$ and the edges $uq$ and $uv_{0}$.
Then $H\cong F(|X|-3,|Y|-2,1,0,0)$ (with the vertices labeled $p,q$ in Figure~\ref{fig:Fk} corresponding to $v_{0},q$ in $H$, respectively).
Hence Lemmas~\ref{lem:subgraph} and~\ref{lem:Fk} give $\alpha(\bar{G})\geq\alpha(H)>\frac{5}{4}$, contrary to $\alpha(\bar{G})<\frac{5}{4}$.

If $uv$ is not a cut edge, then at most four vertices of $P_{uv}(G-uv)$ lie in $N_{G}[a_{0}]\cup N_{G}[b_{0}]$, so $|E(P_{uv}(G-uv))|\leq4$ since $r_{0}$ is the unique vertex outside this union.
Since $pu,vq\in E(G-uv)$, $|E(P_{pq}(G-uv))|\leq|E(P_{uv}(G-uv))|+2\leq6$.
By Eq.~\eqref{eq:extreme-sum}, $(1-z_{u})+z_{v}\leq\alpha(G)$, so $z_{u}-z_{v}\geq1-\alpha(G)>0$.
Since $\alpha(G)<\frac{3-\sqrt5}{2}$, we have $(1-\alpha(G))^{2}>\alpha(G)$.
Since $uv\notin E(P_{pq}(G-uv))$, the Cauchy--Schwarz inequality gives
\begin{align}
    \mathcal E_{G}(\mathbf{z})&\geq(z_{u}-z_{v})^{2}+\sum_{ab\in E(P_{pq}(G-uv))}(z_{a}-z_{b})^{2} \notag\\
    &\geq(1-\alpha(G))^{2}+\frac{(z_{p}-z_{q})^{2}}{|E(P_{pq}(G-uv))|}
    >\alpha(G)+\frac16>\frac{71}{50}\alpha(G), \notag
\end{align}
contrary to Eq.~\eqref{eq:diameter-two-bound}.

Now suppose that $d_{G}(p,q)=4$.
Let $A=N_{G}(p)$, $B=N_{G}(q)$, and
\begin{align}
    C=\{c\in V(G)\setminus(\{p,q\}\cup A\cup B):
    N_{G}(c)\cap A\ne\varnothing,\ N_{G}(c)\cap B\ne\varnothing\}. \notag
\end{align}
Since $d_{G}(p,q)=4$, the middle vertex of $P_{pq}(G)$ lies in $C$, so $|C|\geq1$.
For each $c\in C$, choose $a_{c}\in N_{G}(c)\cap A$ and $b_{c}\in N_{G}(c)\cap B$.
By Eq.~\eqref{eq:extreme-sum}, $(1-z_{a_{c}})+z_{b_{c}}\leq\alpha(G)$, so $z_{a_{c}}-z_{b_{c}}\geq1-\alpha(G)>0$.
Since $\alpha(G)<\frac{3-\sqrt5}{2}$, Eq.~\eqref{eq:diameter-two-bound} and the Cauchy--Schwarz inequality give
\begin{align}
    \frac{71}{50}\alpha(G)>\mathcal E_{G}(\mathbf{z})&\geq\sum_{c\in C}\bigl((z_{a_{c}}-z_{c})^{2}+(z_{c}-z_{b_{c}})^{2}\bigr) \notag\\
    &\geq\frac12\sum_{c\in C}(z_{a_{c}}-z_{b_{c}})^{2}\geq\frac{|C|}{2}(1-\alpha(G))^{2}>\frac{|C|}{2}\alpha(G). \notag
\end{align}
Thus $|C|<\frac{71}{25}<3$, so $|C|\in\{1,2\}$.

Let $U=\{u\in V(G)\setminus(N_{G}[p]\cup N_{G}[q]\cup C):N_{G}(u)\cap B=\varnothing\}$ and $V=V(G)\setminus(N_{G}[p]\cup N_{G}[q]\cup C\cup U)$.
Then
\begin{align}
    E_{G}(N_{G}[p],N_{G}[q])=E_{G}(\{p,q\},C\cup U\cup V)=E_{G}(U,B)=E_{G}(V,A)=\varnothing. \notag
\end{align}
Then $F(|A|-1,|B|-1,|C|,|U|,|V|)$ is isomorphic to a spanning subgraph of $\bar{G}$.
Lemmas~\ref{lem:subgraph} and~\ref{lem:Fk} give $\alpha(\bar{G})>\frac{5}{4}$, contrary to $\alpha(\bar{G})<\frac{5}{4}$.
\end{proof}

\begin{lemma}\label{lem:distant-pairs}
Let $G$ be a connected graph with at least two vertices.
If all pairs $u,v\in V(G)$ with $d_{G}(u,v)\geq3$ satisfy $d_{G}(u),d_{G}(v)\geq2$ and $d_{G}(u,w),d_{G}(v,w)\leq2$ for all $w\in V(G)\setminus\{u,v\}$, then $\alpha(G)\geq1$.
\end{lemma}
\begin{proof}
Suppose $\alpha(G)<1$.
Let $\mathbf{x}$ be a Fiedler vector of $G$ and $p,q\in V(G)$ satisfy $x_{p}=\max_{v\in V(G)}x_{v}$ and $x_{q}=\min_{v\in V(G)}x_{v}$.
By Lemma~\ref{lem:normalization}\textup{(i)}, $d_{G}(p,q)\geq3$.
Then $d_{G}(p),d_{G}(q)\geq2$ and $d_{G}(p,w),d_{G}(q,w)\leq2$ for all $w\in V(G)\setminus\{p,q\}$.
Put $V_\pm=\{v:\pm x_{v}>0\}$.
Since $\sum_{v\in V(G)}x_{v}=0$, the vector $\mathbf{f}=\frac{\mathbf{x}}{\sum_{v\in V_+}x_{v}}$ satisfies $f_{p}>0>f_{q}$ and $\sum_{v\in V_+}f_{v}=-\sum_{v\in V_-}f_{v}=1$.
Since $\alpha(G)<1$, evaluating $L(G)\mathbf{f}=\alpha(G)\mathbf{f}$ at $p$ and $q$ gives
\begin{align}
    \sum_{v\in N_{G}(p)}(f_{p}-f_{v})=\alpha(G)f_{p}<f_{p},\qquad
    \sum_{v\in N_{G}(q)}(f_{v}-f_{q})=-\alpha(G)f_{q}<-f_{q}. \notag
\end{align}
Since all summands are nonnegative, we have $N_{G}(p)\subseteq V_+$ and $N_{G}(q)\subseteq V_-$.
Since $\sum_{v\in V_+}f_{v}=1$,
\begin{align}
    1\geq f_{p}+\sum_{v\in N_{G}(p)}f_{v}=(d_{G}(p)+1-\alpha(G))f_{p}>2f_{p}. \notag
\end{align}
Similarly, $1>-2f_{q}$, so $f_{p}-f_{q}<1$.
If $v\in V_+\setminus\{p\}$, then $d_{G}(v,q)=2$, and the middle vertex of $P_{vq}(G)$ belongs to $N_{G}(q)\subseteq V_-$.
Similarly, every vertex of $V_-\setminus\{q\}$ has a neighbor in $V_+$.
Summing $(L(G)\mathbf{f})_{v}=\alpha(G)f_{v}$ over $v\in V_+$ and using $\sum_{v\in V_+}f_{v}=-\sum_{v\in V_-}f_{v}=1$, we obtain
\begin{align}
    \alpha(G)
    &=\sum_{\substack{a\in V_+,\ b\notin V_+\\ab\in E(G)}}(f_{a}-f_{b})
    \geq\sum_{ab\in E_{G}(V_+,V_-)}(f_{a}-f_{b}) \notag\\
    &\geq\sum_{v\in V_+\setminus\{p\}}f_{v}-\sum_{v\in V_-\setminus\{q\}}f_{v}
    =2-f_{p}+f_{q}>1, \notag
\end{align}
a contradiction.
\end{proof}

\begin{lemma}\label{lem:two-vertices}
Let $G$ be a graph and $pq\in E(G)$ satisfy $N_{G}[p]\cup N_{G}[q]=V(G)$.
Put $X=N_{G}(p)\setminus N_{G}[q]$, $Y=N_{G}(q)\setminus N_{G}[p]$ and $C=N_{G}(p)\cap N_{G}(q)$.
Suppose that $d_{G}(c)=2$ for all $c\in C$ and that each nontrivial component of $G[X\cup Y]$ is a complete bipartite graph $K_{a_{i},b_{i}}$, with parts of sizes $a_{i}$ in $X$ and $b_{i}$ in $Y$, respectively.
For $0<\tau<1$, put $\gamma_{0}=(1-\tau)\left(1+\frac{|C|}{2-\tau}\right)$.
Then
\begin{align}
    \alpha(G)>\tau\quad\text{if and only if}\quad
    \left(\gamma_{0}+\sum_{i}\frac{a_{i} b_{i}}{1-\tau+a_{i}+b_{i}}\right)
    \left(\frac{1}{\gamma_{0}+|X|}+\frac{1}{\gamma_{0}+|Y|}\right)>\tau. \notag
\end{align}
The equivalence also holds with $\geq$ in both places.
\end{lemma}

\begin{proof}
The case $G=K_{2}$ is immediate.
We may assume $|G|\geq3$.
With the vertices ordered as $p,q,X,Y,C$, write
\begin{align}
    L(G)-\tau I=\begin{pmatrix}A&R\\R^{\mathsf T}&D\end{pmatrix}, \notag
\end{align}
where
\begin{align}
    \begin{gathered}
    A=\begin{pmatrix}
    1-\tau+|X|+|C|&-1\\
    -1&1-\tau+|Y|+|C|
    \end{pmatrix}\\
    \text{with rows and columns corresponding to $p,q$},\\[2pt]
    R=-\begin{pmatrix}
    \one_{|X|}^{\mathsf T}&0&\one_{|C|}^{\mathsf T}\\
    0&\one_{|Y|}^{\mathsf T}&\one_{|C|}^{\mathsf T}
    \end{pmatrix},\\[2pt]
    D=\operatorname{diag}\bigl((1-\tau)I+L(G[X\cup Y]),(2-\tau)I\bigr).
    \end{gathered} \notag
\end{align}
Since $1-\tau>0$, $D$ is positive definite.
By Lemma~\ref{lem:schur}, the smaller and larger eigenvalues of the Schur complement $M=A-RD^{-1}R^{\mathsf T}$ have the same signs as $-\tau$ and $\alpha(G)-\tau$, respectively.

Computing by components gives
\begin{align}
    M&=A-RD^{-1}R^{\mathsf T} \notag\\
    &=A-\frac{1}{1-\tau}\left[\operatorname{diag}(|X|,|Y|)
    -\sum_{i}\frac{a_{i} b_{i}}{1-\tau+a_{i}+b_{i}}\begin{pmatrix}1&-1\\-1&1\end{pmatrix}\right]
    -\frac{|C|}{2-\tau}\begin{pmatrix}1&1\\1&1\end{pmatrix} \notag\\
    &=\frac{1}{1-\tau}\left[
    \left(\gamma_{0}+\sum_{i}\frac{a_{i} b_{i}}{1-\tau+a_{i}+b_{i}}\right)
    \begin{pmatrix}1&-1\\-1&1\end{pmatrix}
    -\tau\begin{pmatrix}\gamma_{0}+|X|&0\\0&\gamma_{0}+|Y|\end{pmatrix}
    \right]. \notag
\end{align}

Thus $\alpha(G)>\tau$ if and only if $\det M<0$, and $\alpha(G)\geq\tau$ if and only if $\det M\leq0$.
\begin{align}
    \det M&=\frac{\tau(\gamma_{0}+|X|)(\gamma_{0}+|Y|)}{(1-\tau)^{2}} \notag\\
    &\quad\left[
    \tau-\left(\gamma_{0}+\sum_{i}\frac{a_{i} b_{i}}{1-\tau+a_{i}+b_{i}}\right)
    \left(\frac{1}{\gamma_{0}+|X|}+\frac{1}{\gamma_{0}+|Y|}\right)
    \right]. \notag
\end{align}
Thus both equivalences follow.
\end{proof}

\begin{lemma}\label{lem:delete-ecc-two}
Let $G$ be a connected graph of order $n\geq3$ and $w\in V(G)$.
If $G-w$ is connected and $d_{G}(w,v)\leq2$ for every $v\in V(G)$, then
\begin{align}
    \alpha(G)&\geq\min\left\{\alpha(G-w),\frac{n-1}{n}\right\}. \notag
\end{align}
\end{lemma}

\begin{proof}
It suffices to consider $\alpha(G)<1$.
Let $\mathbf{x}$ be a Fiedler vector of $G$ with $x_{w}\geq0$, and choose $b\in V(G)$ such that $x_{b}=\min_{v\in V(G)}x_{v}<0$.
Let $u$ be the neighbor of $w$ on $P_{wb}(G)$.
Since $d_{G}(w,b)\leq2$, either $u=b$ or $u\in N_{G}(b)$.
Since $x_{j}\geq x_{b}$ for every $j\in N_{G}(b)$, evaluating $L(G)\mathbf{x}=\alpha(G)\mathbf{x}$ at $b$ gives $x_{u}-x_{b}\leq\sum_{j\in N_{G}(b)}(x_{j}-x_{b})=-\alpha(G)x_{b}<-x_{b}$, so $x_{u}<0$.
Since $\sum_{i\ne w}x_{i}=-x_{w}$, we have $\sum_{i\ne w}(x_{i}+\frac{x_{w}}{n-1})=0$ and $\sum_{i\ne w}(x_{i}+\frac{x_{w}}{n-1})^{2}=1-\frac{n}{n-1}x_{w}^{2}\geq0$.
By Eq.~\eqref{eq:variational}, we have
\begin{align}
    \alpha(G)=\mathcal E_{G}(\mathbf{x})
    &\geq (x_{w}-x_{u})^{2}+\alpha(G-w)\sum_{i\ne w}\left(x_{i}+\frac{x_{w}}{n-1}\right)^{2} \notag\\
    &\geq x_{w}^{2}+\alpha(G-w)\left(1-\frac{n}{n-1}x_{w}^{2}\right) \notag\\
    &=\frac{nx_{w}^{2}}{n-1}\cdot\frac{n-1}{n}
    +\left(1-\frac{nx_{w}^{2}}{n-1}\right)\alpha(G-w) \notag\\
    &\geq \min\left\{\frac{n-1}{n},\,\alpha(G-w)\right\}. \notag
\end{align}
\end{proof}

\begin{lemma}\label{lem:matching-two}
	Let $G$ be a graph with at least seven vertices and $p,q\in V(G)$ satisfy $d_{G}(p,q)=3$.
	Put $U=N_{G}(p)$, $V=N_{G}(q)$ and $F=G[U,V]$.
	If $\nu(F)=2$, then the following hold.
	\begin{enumerate}
		\item[\textup{(i)}] If $F\cong 2K_{2}\cup(|F|-4)K_{1}$, then $\alpha(\bar{G})\geq1$.
		\item[\textup{(ii)}] If $F\cong P_{4}\cup(|F|-4)K_{1}$ (resp.\ $K_{2,2}\cup(|F|-4)K_{1}$), then $\alpha(\bar{G})>\frac23$ (resp.\ $\alpha(\bar{G})>\frac12$).
		\item[\textup{(iii)}] If $F$ has at least five nonisolated vertices and contains no $K_{2,3}$, then $\alpha(\bar{G})>\frac47$.
	\end{enumerate}
\end{lemma}

\begin{proof}
	Put $W=V(G)\setminus(\{p,q\}\cup U\cup V)$.
Let $H$ be the spanning subgraph of $\bar{G}$ with edge set
\begin{align}
    E(H)=E_{\bar{G}}(\{p\}\cup U,\{q\}\cup V)
    \cup\bigcup_{w\in W}\{pw,qw\}. \notag
\end{align}

In \textup{(i)}, write $E(F)=\{u_{1}v_{1},u_{2}v_{2}\}$ with $u_{i}\in U$ and $v_{i}\in V$ for $i=1,2$.
For distinct $a,b\in V(H)$ with $ab\notin E(F)$, we have $d_{H}(a,b)\leq2$.
Since $d_{H}(u_{i})=|V|\geq2$ and $d_{H}(v_{i})=|U|\geq2$ for $i=1,2$, Lemmas~\ref{lem:subgraph} and~\ref{lem:distant-pairs} give $\alpha(\bar{G})\geq\alpha(H)\geq1$.

For a subgraph of $H$ retaining $pq$ and all edges from $p,q$ to its other vertices, the sets $X,Y,C$ in Lemma~\ref{lem:two-vertices} consist of its vertices in $V,U,W$, respectively.

\textup{(ii)} If $W\ne\varnothing$, then choose $w\in W$ and let $H_{1}$ be the induced subgraph of $H$ with $V(H_{1})=\{p,q,w\}\cup\{v\in V(F):d_{F}(v)>0\}$.
For $H_{1}$, we have $|X|=|Y|=2$ and $C=\{w\}$.
The graph $H_{1}[X\cup Y]$ is isomorphic to $K_{2}\cup2K_{1}$ in the $P_{4}$ case and to $4K_{1}$ in the $K_{2,2}$ case.
In the two cases, take $\tau=\frac23$ and $\tau=\frac12$, respectively.
Then
\begin{align}
    \left(\frac7{12}+\frac37\right)\frac{24}{31}=\frac{170}{217}>\frac23,
    \qquad
    \frac56\cdot\frac{12}{17}=\frac{10}{17}>\frac12. \notag
\end{align}
By Lemma~\ref{lem:two-vertices}, $\alpha(H_{1})>\frac23$ in the $P_{4}$ case and $\alpha(H_{1})>\frac12$ in the $K_{2,2}$ case.
Lemmas~\ref{lem:subgraph} and~\ref{lem:delete-ecc-two} give $\alpha(\bar{G})\geq\alpha(H)\geq\min\{\alpha(H_{1}),\frac67\}$.
Since $\frac67>\frac23$, the required bounds follow.

If $W=\varnothing$, then $|F|=|G|-2\geq5$, so $F$ has an isolated vertex, say $w$.
Let $H_{2}$ be the subgraph of $H$ with
\begin{align}
    V(H_{2})&=\{p,q,w\}\cup\{v\in V(F):d_{F}(v)>0\}, \notag\\
    E(H_{2})&=E(H[V(H_{2})])\setminus E_{H}(U\setminus\{w\},V\setminus\{w\}). \notag
\end{align}
For $H_{2}$, we have $\{|X|,|Y|\}=\{2,3\}$ and $C=\varnothing$.
The only nontrivial component of $H_{2}[X\cup Y]$ is $K_{1,2}$.
For $\tau=\frac23$, we have $\left(\frac13+\frac35\right)\left(\frac37+\frac3{10}\right)=\frac{17}{25}>\frac23$.
By Lemma~\ref{lem:two-vertices}, $\alpha(H_{2})>\frac23$.
Lemmas~\ref{lem:subgraph} and~\ref{lem:delete-ecc-two} give $\alpha(\bar{G})\geq\alpha(H)\geq\min\{\alpha(H_{2}),\frac67\}>\frac23$.

\textup{(iii)} Since $\nu(F)=2$, Lemma~\ref{lem:konig} gives a vertex cover $\{u,v\}$ of $F$.

First suppose that $u\in U$ and $v\in V$.
Let $H_{3}$ be the subgraph of $H$ with
\begin{align}
    V(H_{3})&=\{p,q\}\cup U\cup V, \notag\\
    E(H_{3})&=E(H-W)\setminus\bigl(E_{H}(u,V)\cup E_{H}(U,v)\bigr). \notag
\end{align}
For $H_{3}$, we have $C=\varnothing$.
Since $\{u,v\}$ covers $F$, $H_{3}[X\cup Y]\cong K_{|X|-1,|Y|-1}\cup2K_{1}$.
Since $|X|,|Y|\geq2$, we have
\begin{align}
    \begin{gathered}
    \min\left\{\frac{|X|-\frac47}{|X|+\frac37},\frac{|Y|-\frac47}{|Y|+\frac37}\right\}\geq\frac{10}{17},\\
    \frac{\left(|X|-\frac47\right)\left(|Y|-\frac47\right)}{|X|+|Y|-\frac{11}{7}}
    \left(\frac{1}{|X|+\frac37}+\frac{1}{|Y|+\frac37}\right)
    \geq\frac{10}{17}\frac{|X|+|Y|-\frac87}{|X|+|Y|-\frac{11}{7}}>\frac47.
    \end{gathered} \notag
\end{align}
By Lemma~\ref{lem:two-vertices}, $\alpha(H_{3})>\frac47$.
Lemmas~\ref{lem:subgraph} and~\ref{lem:delete-ecc-two} give $\alpha(\bar{G})\geq\alpha(H)\geq\min\{\alpha(H_{3}),\frac67\}>\frac47$.

Without loss of generality, assume that $u,v\in U$.
Let $H_{4}$ be the induced subgraph of $H$ with $V(H_{4})=\{p,q\}\cup\{x\in V(F):d_{F}(x)>0\}$.
For $H_{4}$, we have $Y=\{u,v\}$ and $C=\varnothing$.
The nontrivial components of $H_{4}[X\cup Y]$ are stars centered at $u$ or $v$ and together contain $j=|\{x\in V:d_{F}(x)=1\}|$ vertices in $X$.
Since $F$ has at least five nonisolated vertices and contains no $K_{2,3}$, we have $3\leq|X|=j+|N_{F}(u)\cap N_{F}(v)|\leq j+2$, so $j\geq1$.
Then
\begin{align}
    \begin{gathered}
    \sum_{i}\frac{a_{i} b_{i}}{a_{i}+b_{i}+\frac37}
     =\sum_{i}\frac{7a_{i} b_{i}}{7a_{i} b_{i}+10}
     \geq\frac{7j}{7j+10},\\
    \left(\frac37+\sum_{i}\frac{a_{i} b_{i}}{a_{i}+b_{i}+\frac37}\right)
    \left(\frac{1}{|X|+\frac37}+\frac7{17}\right)
    \geq\frac{10}{17}\left(1+\frac{70j-68}{(7j+10)(7j+17)}\right)>\frac47.
    \end{gathered} \notag
\end{align}
By Lemma~\ref{lem:two-vertices}, $\alpha(H_{4})>\frac47$.
Lemmas~\ref{lem:subgraph} and~\ref{lem:delete-ecc-two} give $\alpha(\bar{G})\geq\alpha(H)\geq\min\{\alpha(H_{4}),\frac67\}>\frac47$.
\end{proof}

\Needspace{9\baselineskip}
\begin{lemma}\label{lem:matching-exclusion}
Let $G$ be a connected graph of order $n$ and $\mathbf{x},\mathbf{y}$ be Fiedler vectors of $G,\bar{G}$, respectively.
Choose $p,q\in V(G)$ such that $x_{p}=\max_{v\in V(G)}x_{v}$ and $x_{q}=\min_{v\in V(G)}x_{v}$.
If $d_{G}(p,q)=3$, $x_{p}-x_{q}\geq\max_{u,v\in V(G)}|y_{u}-y_{v}|$, and $S(G)<\frac{10}{7}$, then $\nu(G[N_{G}(p),N_{G}(q)])=1$.
\end{lemma}

\begin{proof}
Put $F=G[N_{G}(p),N_{G}(q)]$.
Since $d_{G}(p,q)=3$, $E(F)\ne\varnothing$ and $\bar{G}$ is connected.
Suppose that  $\nu(F)\geq2$.
Then $n\geq6$.
If $n=6$, then $\alpha(G)\geq\alpha(C_{6})=1$ by Lemma~\ref{lem:subgraph}.
Since $|E(F)|\leq4$ and $\sigma(F)=\max_{uv\in E(F)}\{d_{F}(u)+d_{F}(v)\}\leq4$, Lemma~\ref{lem:bipartite-bound} gives $\alpha(\bar{G})\geq\frac{6}{6+4+4}=\frac37$.
Thus $S(G)\geq\frac{10}{7}$, a contradiction.

Assume that $n\geq7$.
For a subgraph $F_{0}$ of $F$ with at least one edge and no isolated vertices, put $h=|F_{0}|$ and $\kappa=\kappa_{F_{0}}$.
By Lemma~\ref{lem:energy-bound}\textup{(i)}, $h-4\kappa>0$.
Define
\begin{align}
    f(\kappa,h,s)=\kappa+\frac{\max\{2\kappa+s-\frac{10}{7},0\}^{2}}{h-4\kappa}+\frac{10}{7}s. \notag
\end{align}
$f$ is increasing in $s$.
Since $S(G)<\frac{10}{7}$, we have $0<\alpha(G)<\frac{10}{7}-\alpha(\bar{G})$ and
$\alpha(G)S(G)<\frac{10}{7}\left(\frac{10}{7}-\alpha(\bar{G})\right)$.
Applying Lemma~\ref{lem:energy-bound}\textup{(ii)} if $2\kappa+\alpha(\bar{G})>\frac{10}{7}$ and Eq.~\eqref{eq:retained-energy} otherwise, we obtain
\begin{equation}\label{eq:matching-exclusion-contradiction}
f(\kappa,h,\alpha(\bar{G}))<\frac{100}{49}.
\end{equation}

If $\nu(F)\geq5$, then take $F_{0}\cong5K_{2}$.
Here $h=10$ and $\kappa=5\kappa_{K_{2}}=\frac53$ (see Remark~\ref{rem:energy-components}).
Since $\alpha(\bar{G})\geq0$, we have $f(\kappa,h,\alpha(\bar{G}))\geq f(\frac53,10,0)=\frac{135}{49}>\frac{100}{49}$, contrary to Eq.~\eqref{eq:matching-exclusion-contradiction}.

If $\nu(F)=4$, then take $F_{0}\cong4K_{2}$.
By Lemma~\ref{lem:bipartite-bound}, $\alpha(\bar{G})\geq3-2\sqrt2>\frac16$.
Then $f(\frac43,8,\frac16)=\frac{10873}{4704}>\frac{100}{49}$, contrary to Eq.~\eqref{eq:matching-exclusion-contradiction}.

If $\nu(F)=3$ and $\Delta(F)\geq5$, then $F$ contains a subgraph $F_{0}\cong K_{1,3}\cup2K_{2}$.
By Lemma~\ref{lem:bipartite-bound}, $\alpha(\bar{G})\geq\frac{5-\sqrt{21}}2>\frac15$.
Then $f(\frac{19}{15},8,\frac15)=\frac{22991}{10780}>\frac{100}{49}$, contrary to Eq.~\eqref{eq:matching-exclusion-contradiction}.

If $\nu(F)=3$, $\Delta(F)\leq4$ and $|E(F)|\geq4$, then $F$ contains a subgraph $F_{0}$ isomorphic to $K_{1,2}\cup2K_{2}$ or $P_{4}\cup K_{2}$.
Then $n\geq|F|+2\geq\Delta(F)+5$, and Lemma~\ref{lem:konig} gives $|E(F)|\leq3\Delta(F)$.
Since $F$ is bipartite, $\sigma(F)\leq|F|\leq n-2$.
Applying Lemma~\ref{lem:bipartite-bound} with the above inequalities, we obtain
\begin{align}
    \alpha(\bar{G})
    \geq\frac{n}{n+|E(F)|+\sigma(F)}
    \geq\frac{n}{2n+3\Delta(F)-2}
    \geq\frac{\Delta(F)+5}{5\Delta(F)+8}
    \geq\frac9{28}. \notag
\end{align}
For $F_{0}\cong K_{1,2}\cup2K_{2}$ and $F_{0}\cong P_{4}\cup K_{2}$, respectively, we have
\begin{align}
    f\left(\frac76,7,\frac9{28}\right)=\frac{12459}{5488}>\frac{100}{49},
    \qquad
    f\left(\frac{22}{21},6,\frac9{28}\right)=\frac{60965}{29792}>\frac{100}{49}, \notag
\end{align}
contrary to Eq.~\eqref{eq:matching-exclusion-contradiction}.

If $\nu(F)=3$, then $|E(F)|=3$, so we may take $F_{0}\cong3K_{2}$.
Since $\sigma(F)=2$, Lemma~\ref{lem:bipartite-bound} gives $\alpha(\bar{G})\geq\frac{n}{n+5}>\frac12$.
Then $f(1,6,\frac12)=\frac{897}{392}>\frac{100}{49}$, contrary to Eq.~\eqref{eq:matching-exclusion-contradiction}.

Thus $\nu(F)=2$.
If $F$ contains $K_{2,3}$, then take $F_{0}\cong K_{2,3}$.
By Lemma~\ref{lem:bipartite-bound}, $\alpha(\bar{G})\geq2-\sqrt3>\frac4{15}$.
Then $f(1,5,\frac4{15})=\frac{22969}{11025}>\frac{100}{49}$, contrary to Eq.~\eqref{eq:matching-exclusion-contradiction}.

If $F$ contains no $K_{2,3}$ and has at least five nonisolated vertices, then take $F_{0}\cong K_{1,2}\cup K_{2}$.
By Lemma~\ref{lem:matching-two}\textup{(iii)}, $\alpha(\bar{G})>\frac47$.
Then $f(\frac56,5,\frac47)=\frac{143}{70}>\frac{100}{49}$, contrary to Eq.~\eqref{eq:matching-exclusion-contradiction}.

The four nonisolated vertices of $F$ induce $F_{0}\cong2K_{2}$, $P_{4}$ or $K_{2,2}$.
By Lemma~\ref{lem:matching-two}\textup{(i)} and \textup{(ii)}, the corresponding lower bounds for $\alpha(\bar{G})$ are $1$, $\frac23$ and $\frac12$, respectively.
Then
\begin{align}
    f\left(\frac23,4,1\right)=\frac{531}{196},
    \quad f\left(\frac57,4,\frac23\right)=\frac{37}{18},
    \quad f\left(\frac45,4,\frac12\right)=\frac{1629}{784}. \notag
\end{align}
Each of these values exceeds $\frac{100}{49}$, contrary to Eq.~\eqref{eq:matching-exclusion-contradiction}.
\end{proof}

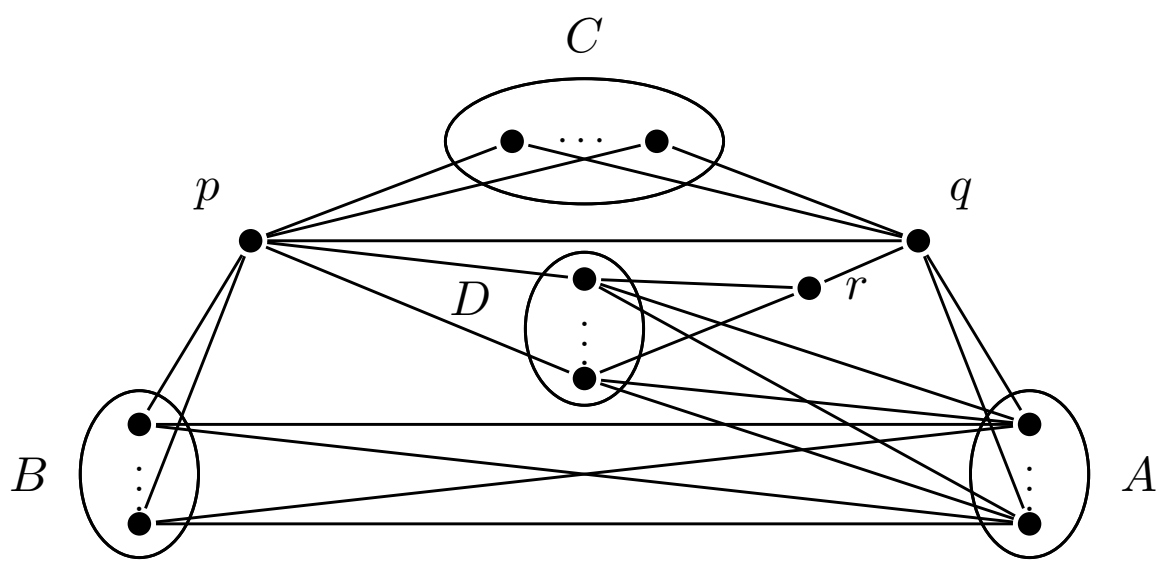
\begin{figure}[H]
\centering
\begin{tikzpicture}[papergraph,x=.8cm,y=.55cm,scale=.92,transform shape]
\begin{pgfinterruptboundingbox}
\node[vertex,label=above left:$p$] (p) at (-3,0) {};
\node[vertex,label=above right:$q$] (q) at (3,0) {};
\node[vertex,label=right:$r$] (r) at (2.02,-0.62) {};
\node[block,minimum width=20mm,minimum height=9mm,inner sep=0pt] (C) at (0,1.3) {};
\node[vertex] (C1) at (-0.65,1.3) {};
\node[vertex] (C2) at (0.65,1.3) {};
\node[font=\scriptsize,inner sep=0pt] at (0,1.3) {$\cdots$};
\node[anchor=south] at ($(C.north)+(0,.08)$) {$C$};
\node[block,minimum width=8.5mm,minimum height=11mm,inner sep=0pt] (D) at (0,-1.15) {};
\node[vertex] (D1) at (0,-1.7999999999999998) {};
\node[vertex] (D2) at (0,-0.4999999999999999) {};
\node[font=\scriptsize,inner sep=0pt] at (0,-1.15) {$\vdots$};
\node[anchor=east] at ($(D.west)+(-.12,.40)$) {$D$};
\node[block,minimum width=8.5mm,minimum height=12mm,inner sep=0pt] (B) at (-4,-3.05) {};
\node[vertex] (B1) at (-4,-3.6999999999999997) {};
\node[vertex] (B2) at (-4,-2.4) {};
\node[font=\scriptsize,inner sep=0pt] at (-4,-3.05) {$\vdots$};
\node[anchor=east] at ($(B.west)+(-.10,0)$) {$B$};
\node[block,minimum width=8.5mm,minimum height=12mm,inner sep=0pt] (A) at (4,-3.05) {};
\node[vertex] (A1) at (4,-3.6999999999999997) {};
\node[vertex] (A2) at (4,-2.4) {};
\node[font=\scriptsize,inner sep=0pt] at (4,-3.05) {$\vdots$};
\node[anchor=west] at ($(A.east)+(.10,0)$) {$A$};
\draw (p)--(q);
\draw (q)--(r);
\draw (p)--(D1);
\draw (p)--(D2);
\draw (p)--(C1);
\draw (p)--(C2);
\draw (q)--(C1);
\draw (q)--(C2);
\draw (r)--(D1);
\draw (r)--(D2);
\draw (A1)--(D1);
\draw (A1)--(D2);
\draw (A2)--(D1);
\draw (A2)--(D2);
\draw (p)--(B2);
\draw (p)--(B1);
\draw (q)--(A2);
\draw (q)--(A1);
\draw (B2)--(A2);
\draw (B1)--(A1);
\draw (B1)--(A2);
\draw (B2)--(A1);
\draw (0,1.3) ellipse [x radius=10mm,y radius=4.5mm];
\draw (0,-1.15) ellipse [x radius=4.25mm,y radius=5.5mm];
\draw (-4,-3.05) ellipse [x radius=4.25mm,y radius=6mm];
\draw (4,-3.05) ellipse [x radius=4.25mm,y radius=6mm];
\foreach \v in {p,q,r,C1,C2,D1,D2,B1,B2,A1,A2} {\fill (\v) circle[radius=1.7pt];}
\end{pgfinterruptboundingbox}
\path[use as bounding box] (-5.8388,-4.2804) rectangle (5.8445,3.9357);
\end{tikzpicture}
\caption{$H(a,b,c,d)$.}
\label{fig:H0}
\end{figure}

For nonnegative integers $a,b,c,d$, let $H(a,b,c,d)$ be obtained from pairwise disjoint independent vertex sets $\{p\},\{q\},\{r\},A,B,C,D$ by adding all possible edges between the following pairs of sets: $(\{p\}\cup A,\{q\}\cup B\cup D)$, $(\{r\},\{q\}\cup D)$ and $(\{p,q\},C)$, where $|A|=a$, $|B|=b\geq1$, $|C|=c$ and $|D|=d$ (see Figure~\ref{fig:H0}).

Define
\begin{align}
    W_{a,b,c,d}(\tau)
    &=\left(\gamma_{0}+\frac{ab}{1-\tau+a+b}+\frac{d}{2-\tau+d}\right)
       (2\gamma_{0}+a+b+d+1) \notag\\
    &\quad -\tau(\gamma_{0}+b+d)(\gamma_{0}+a+1), \notag
\end{align}
where $a,b,c,d$ are nonnegative integers, $0<\tau<1$, and $\gamma_{0}=(1-\tau)\left(1+\frac{c}{2-\tau}\right)$.

\begin{lemma}\label{lem:H0}
Let $H(a,b,c,d)$ be as defined above.
Then the following hold.
\begin{enumerate}
\item[\textup{(i)}] If $ad>0$, then $\alpha(H(a,b,c,d))\geq1$.
\item[\textup{(ii)}] If $ad=0$ and $0<\tau<1$, then $\alpha(H(a,b,c,d))>\tau$ if and only if $W_{a,b,c,d}(\tau)>0$.
In addition, $W_{a,b,c,0}(\tau)>0$ implies $W_{0,b,c,a}(\tau)>0$.
\item[\textup{(iii)}] If $c\geq1$ and $b\geq2$, then $\alpha(H(0,b,c,0))>\frac37+\frac1{2(b+1)}$.
\item[\textup{(iv)}] If $a,c\geq1$ and $b\geq5$, then $\min\{\alpha(H(a,b,c,0)),\alpha(H(0,b,c,a))\}>\frac35+\frac1{b+2}$.
\item[\textup{(v)}] If $a,c\geq1$, then
\begin{align}
    \min\{\alpha(H(a,b,c,0)),\alpha(H(0,b,c,a))\}>
    \begin{cases}
    1-\frac1{2(a+2)}, & b=1;\\
    \frac45, & b=2,\ a=c=1;\\
    \frac67, & b=2,\ a+c\geq3;\\
    \frac34, & b=3,4.
    \end{cases} \notag
\end{align}
\end{enumerate}
\end{lemma}

\begin{proof}
Put $H=H(a,b,c,d)$.

\textup{(i)}
Suppose that $\alpha(H)<1$, and let $\mathbf{x}$ be a Fiedler vector of $H$.
Let $u_{0},v_{0}\in V(H)$ satisfy $x_{u_{0}}=\max_{u\in V(H)}x_{u}$ and $x_{v_{0}}=\min_{u\in V(H)}x_{u}$.
By Lemma~\ref{lem:normalization}\textup{(i)}, $d_{H}(u_{0},v_{0})\geq3$, so by the construction of $H$, $\{u_{0},v_{0}\}=\{r,v\}$ for some $v\in B$.
Replacing $\mathbf{x}$ by $-\mathbf{x}$ if necessary, we may assume that $x_{r}>x_{v}$.
Put $\mathbf{f}=\frac{\mathbf{x}-x_{v}\one}{x_{r}-x_{v}}$.
Then $f_{r}=1$, $f_{v}=0$ and $0\leq f_{u}\leq1$ for all $u\in V(H)$.
For $w\in C$, Eq.~\eqref{eq:vertex-equations} gives
\begin{align}
    \alpha(H)f_{w}=(L(H)\mathbf{f})_{w}-(L(H)\mathbf{f})_{v}=2f_{w}-f_{q}+\sum_{u\in A}f_{u}. \notag
\end{align}
Since $f_{u}\geq0$ for $u\in A$ and $2-\alpha(H)>1$, we have $f_{w}\leq f_{q}$ for $w\in C$.
By Eq.~\eqref{eq:vertex-equations}, we have
\begin{align}
    \alpha(H)f_{q}
    &=(L(H)\mathbf{f})_{q}-(L(H)\mathbf{f})_{v} \notag\\
    &=(a+c+2)f_{q}-1-\sum_{w\in C}f_{w}\geq(a+2)f_{q}-1; \notag\\
    \alpha(H)f_{u}
    &=(L(H)\mathbf{f})_{u}-(L(H)\mathbf{f})_{v}=(a+2)f_{u}-1,\qquad u\in D. \notag
\end{align}
Hence $f_{u}\leq\frac{1}{a+2-\alpha(H)}<\frac{1}{2}$ for $u\in\{q\}\cup D$.
Combining the above inequality with Eq.~\eqref{eq:extreme-sum} gives
\begin{align}
    \alpha(H)
    =\sum_{u\in\{q\}\cup D}(1-f_{u})+\sum_{u\in\{p\}\cup A}f_{u}
    >\frac{d+1}{2}\geq1, \notag
\end{align}
a contradiction.

\textup{(ii)} It is easy to see that $H$ satisfies the hypotheses of Lemma~\ref{lem:two-vertices} with $X=B\cup D$, $Y=A\cup\{r\}$ and $N_{H}(p)\cap N_{H}(q)=C$.
Since $ad=0$, $H[X\cup Y]$ is edgeless if $a=d=0$ and otherwise has a unique nontrivial component, $K_{b,a}$ if $a>0$ and $K_{d,1}$ if $d>0$.
By Lemma~\ref{lem:two-vertices} and $(\gamma_{0}+b+d)(\gamma_{0}+a+1)>0$, we have $W_{a,b,c,d}(\tau)>0$ if and only if $\alpha(H)>\tau$.

Direct calculation gives $W_{0,b,c,a}(\tau)\geq\frac{(2-\tau)(1-\tau+a+b)}{(2-\tau+a)(1-\tau+b)}W_{a,b,c,0}(\tau)$, so $W_{a,b,c,0}(\tau)>0$ implies $W_{0,b,c,a}(\tau)>0$.

\textup{(iii)} Take $\tau=\frac37+\frac1{2(b+1)}$.
Since $b\geq2$, we have $\tau\leq\frac{25}{42}<\frac35$.
Since $c\geq1$, we have $\gamma_{0}\geq\frac{(1-\tau)(3-\tau)}{2-\tau}>\frac{3-\tau}{3+\tau}>\frac23$.
Then
\begin{align}
    \frac{W_{0,b,c,0}(\tau)}{(\gamma_{0}+b)(\gamma_{0}+1)}
    =\frac{\gamma_{0}}{\gamma_{0}+1}+\frac{\gamma_{0}}{\gamma_{0}+b}-\tau
    >\frac{3-\tau}{6}+\frac{2}{3(b+1)}-\tau
    =\frac1{12(b+1)}>0. \notag
\end{align}
Since $(\gamma_{0}+b)(\gamma_{0}+1)>0$, \textup{(ii)} gives $\alpha(H(0,b,c,0))>\tau=\frac37+\frac1{2(b+1)}$.

\textup{(iv)} For fixed $b$ and $\tau\in(0,1)$, $W_{a,b,c,0}(\tau)$ is increasing in both $a$ and $c$.
Take $\tau=\frac35+\frac1{b+2}$.
By \textup{(ii)}, the desired bound holds if and only if $W_{a,b,c,0}(\tau)>0$ and $W_{0,b,c,a}(\tau)>0$, and the first inequality implies the second.
Since $a,c\geq1$, it suffices to prove that $W_{1,b,1,0}(\tau)>0$.
If $b=5$, then $\tau=\frac{26}{35}$ and $W_{1,5,1,0}(\tau)=\frac{392597}{206571750}>0$.
If $b\geq6$, then $\tau=\frac35+\frac1{b+2}\leq\frac{29}{40}$, so $(1-\tau)(2-\tau)(5-\tau)-\frac75>0$.
Then
\begin{align}
    \frac{(2-\tau)(b+2-\tau)}{b+1-\tau}W_{1,b,1,0}(\tau)
    >\frac{13b-69}{125}
    +(1-\tau)\left[(1-\tau)(2-\tau)(5-\tau)-\frac75\right]>0. \notag
\end{align}

\textup{(v)} By \textup{(ii)}, each stated bound holds if and only if $W_{a,b,c,0}(\tau)>0$ and $W_{0,b,c,a}(\tau)>0$, and the first inequality implies the second.
For $b=1$, take $\tau=1-\frac1{2(a+2)}$.
Then $W_{1,1,1,0}(\frac56)=\frac{5311}{19656}>0$.
If $a\geq2$, then $(1-\tau)(a+2-\tau)(a+3-\tau)>\frac{a+1}{2}>2-\tau$, so
\begin{align}
    (a+2-\tau)W_{a,1,0,0}(\tau)
    =(2-\tau)\bigl[(1-\tau)(a+2-\tau)(a+3-\tau)-(2-\tau)\bigr]>0. \notag
\end{align}
Since $W_{a,1,c,0}(\tau)$ is increasing in $c$ and $c\geq1$, we have $W_{a,1,c,0}(\tau)>0$ in both cases.

For $b\in\{2,3,4\}$, it suffices to check $W_{a,b,c,0}(\tau)>0$ at the minimal pairs in the table below, since $W_{a,b,c,0}(\tau)$ is increasing in $a$ and $c$.
\begin{align}
    \begin{array}{@{}lccc@{}}
    \toprule
    \text{Case} & \tau & \text{Minimal pair }(a,c) & W_{a,b,c,0}(\tau)\\
    \midrule
    b=2,\ a=c=1 & 4/5 & (1,1) & 213/1000\\
    b=2,\ a+c\geq3 & 6/7 & (1,2) & 67/3773\\
    b=2,\ a+c\geq3 & 6/7 & (2,1) & 18911/39788\\
    b=3 & 3/4 & (1,1) & 2613/5440\\
    b=4 & 3/4 & (1,1) & 83/448\\
    \bottomrule
    \end{array} \notag
\end{align}
Thus $W_{a,b,c,0}(\tau)>0$ in each case, so the stated bounds follow from \textup{(ii)}.
\end{proof}

\begin{lemma}\label{lem:diameter-three}
Let $G$ be a graph of order $n$ with $d(\bar{G})=3$.
Let $\mathbf{x},\mathbf{y}$ be Fiedler vectors of $G,\bar{G}$, respectively, and $p,q\in V(G)$ satisfy $x_{p}=\max_{v\in V(G)}x_{v}$ and $x_{q}=\min_{v\in V(G)}x_{v}$.
Suppose that $d_{G}(p,q)=3$ and $x_{p}-x_{q}\geq\max_{u,v\in V(G)}|y_{u}-y_{v}|$.
If either \textup{(i)} $S(G)<\min\{\frac{71}{50},1+\xi_{n}\}$ or \textup{(ii)} $S(G)<\frac{10}{7}$ and $\alpha(\bar{G})<1$, then $G$ or $\bar{G}$ belongs to $\mathcal H(m_{1},m_{2},1,1)\cup\mathcal H_{n}$ ( $m_{1}+m_{2}=n-2$).
\end{lemma}

\begin{proof}
Since $d(\bar{G})=3$, $G$ is connected.
Both \textup{(i)} and \textup{(ii)} imply $S(G)<\frac{10}{7}$, so Lemma~\ref{lem:matching-exclusion} shows that $F=G[N_{G}(p),N_{G}(q)]$ consists of a star and possibly some isolated vertices.
Since $d_{G}(p,q)=3$, there exists a path $prsq$ of length three, where $r\in N_{G}(p)$ and $s\in N_{G}(q)$.
Every dominating edge of $G$ belongs to $E(F)$.

For $\mathbf{z}=\frac{\mathbf{x}-x_{q}\one}{x_{p}-x_{q}}$, we have $0=z_{q}\leq z_{v}\leq z_{p}=1$ for all $v\in V(G)$.
Lemma~\ref{lem:normalization}\textup{(ii)} gives
\begin{equation}\label{eq:diam-three-energy}
\mathcal E_{G}(\mathbf{z})\leq\alpha(G)S(G)<\frac{10}{7}\alpha(G).
\end{equation}

Without loss of generality, we consider the following two cases.

\emph{Case 1.} $|E(F)|=1$.
Then $E(F)=\{rs\}$, so $prsq$ is the unique $p$--$q$ path of length three and $rs$ is the only dominating edge of $G$.
Put $A_{1}=N_{G}(r)\setminus N_{G}[s]$, $A_{2}=N_{G}(s)\setminus N_{G}[r]$ and $Z=N_{G}(r)\cap N_{G}(s)$.
If either of $A_{1},A_{2}$ is a singleton, then we may assume without loss of generality that $A_{2}=\{q\}$.

\emph{Subcase 1.1.} $rs$ is a cut edge.
Then $Z=\varnothing$ and $E_{G}(A_{1},A_{2})=\varnothing$.
If $A_{2}=\{q\}$, then $G\in\mathcal H(n-3,1,1,1)$ with ordered blocks $(A_{1},\{r\},\{s\},\{q\})$.
Now suppose that $|A_{1}|,|A_{2}|\geq2$.
Since $rs$ is the only edge between $A_{1}\cup\{r\}$ and $A_{2}\cup\{s\}$, $\bar{G}$ contains a spanning subgraph $H=K_{|A_{1}|+1,\,|A_{2}|+1}-rs$.
By Lemma~\ref{lem:subgraph}, $\alpha(H)\leq\alpha(\bar{G})<\frac{10}{7}$.
Let $\mathbf{w}$ be a Fiedler vector of $H$ and $u_{0},v_{0}\in V(H)$ satisfy $w_{u_{0}}=\max_{v\in V(H)}w_{v}$ and $w_{v_{0}}=\min_{v\in V(H)}w_{v}$.
Since $\alpha(H)<2$, Lemma~\ref{lem:normalization}\textup{(i)} gives $u_{0}v_{0}\notin E(H)$ and $|N_{H}(u_{0})\cap N_{H}(v_{0})|\leq1$.
Then $\{u_{0},v_{0}\}=\{r,s\}$.
Put $\mathbf{f}=\frac{\mathbf{w}-w_{s}\one}{w_{r}-w_{s}}$, so $f_{r}=\max_{v\in V(H)}f_{v}=1$ and $f_{s}=\min_{v\in V(H)}f_{v}=0$.
For $u\in A_{1}$ and $v\in A_{2}$, Eq.~\eqref{eq:vertex-equations} gives
\begin{align}
    \alpha(H)(1-f_{u})&=(L(H)\mathbf{f})_{r}-(L(H)\mathbf{f})_{u}=(|A_{2}|+1)(1-f_{u})-1, \notag\\
    \alpha(H)f_{v}&=(L(H)\mathbf{f})_{v}-(L(H)\mathbf{f})_{s}=(|A_{1}|+1)f_{v}-1. \notag
\end{align}
Since $|A_{1}|,|A_{2}|\geq2$ and $\alpha(H)<\frac{10}{7}$, we have $1-f_{u}<\frac7{11}$ and $f_{v}<\frac7{11}$.
Combining the above inequalities with Eq.~\eqref{eq:extreme-sum} gives
\begin{align}
    \alpha(H)=\sum_{v\in A_{2}}(1-f_{v})+\sum_{u\in A_{1}}f_{u}
    >\frac4{11}(|A_{1}|+|A_{2}|)
    >\frac{10}{7}, \notag
\end{align}
a contradiction.

\emph{Subcase 1.2.} $rs$ is not a cut edge and $A_{2}=\{q\}$.
The uniqueness of $prsq$ gives $N_{G}(q)=\{s\}$ and $N_{G}(p)\cap Z=\varnothing$.
Since $rs$ lies on a cycle and $d_{G}(q)=1$, $Z\ne\varnothing$.
If $A_{1}=\{p\}$, then $G\in\mathcal H_{n}$ with path $prsq$ and graph $G[Z]$.
Now suppose that $|A_{1}|\geq2$.
Then $G$ contains a spanning copy of $H(0,|A_{1}|,|Z|,0)$ with $(p,q,r,B,C)$ corresponding to $(r,s,q,A_{1},Z)$ of $G$ (see Figure~\ref{fig:H0}).
Also, $\bar{G}$ contains a spanning copy of $H(|A_{1}|-1,1,|Z|,0)$ with $(p,q,r,A,B,C)$ corresponding to $(p,q,r,A_{1}\setminus\{p\},\{s\},Z)$ of $\bar{G}$ (see Figure~\ref{fig:H0}).

Applying Lemmas~\ref{lem:subgraph} and~\ref{lem:H0}\textup{(iii), (v)}, we have
\begin{align}
    \alpha(G)>\frac37+\frac1{2(|A_{1}|+1)},\qquad
    \alpha(\bar{G})>1-\frac1{2(|A_{1}|+1)}, \notag
\end{align}
so $S(G)>\frac{10}{7}$, a contradiction.

\emph{Subcase 1.3.} $rs$ is not a cut edge and $|A_{1}|,|A_{2}|\geq2$.
The unique pair at distance three in $\bar{G}$ is $\{r,s\}$.
As $d_{\bar{G}}(r)=|A_{2}|\geq2$ and $d_{\bar{G}}(s)=|A_{1}|\geq2$, Lemma~\ref{lem:distant-pairs} gives $\alpha(\bar{G})\geq1$.
Thus \textup{(i)} holds and $0<\alpha(G)\leq S(G)-1<\min\{\frac{21}{50},\xi_{n}\}$.
Put $h=(1-z_{r})+z_{s}$, so $z_{r}-z_{s}=1-h$ and, by Eq.~\eqref{eq:extreme-sum}, $0\leq h\leq\alpha(G)$.

If $Z=\varnothing$, choose $a\in A_{1}$ and $b\in A_{2}$ with $ab\in E(G)$, for $rs$ is not a cut edge.
If $Z\ne\varnothing$, choose $a=b\in Z$.
Then $ra,bs\in E(G)$, and the uniqueness of $prsq$ gives $a,b\notin\{p,q\}$.
Since $z_{r}-z_{s}=1-h$, we have
\begin{equation}\label{eq:alternate-path-energy}
\begin{aligned}
(z_{r}-z_{a})^{2}+(z_{a}-z_{b})^{2}+(z_{b}-z_{s})^{2}
&=\frac{(1-h)^{2}}{3}
+\sum_{(i,j)\in\{(r,a),(a,b),(b,s)\}}\left(z_{i}-z_{j}-\frac{1-h}{3}\right)^{2}\\
&\geq\frac{(1-h)^{2}}{3}+\left(z_{a}-z_{b}-\frac{1-h}{3}\right)^{2}.
\end{aligned}
\end{equation}

Using Eq.~\eqref{eq:vertex-equations} at $r,s$ together with $z_{p}=1$, $z_{q}=0$ and $0\leq h\leq\alpha(G)<\frac{21}{50}$, we obtain
\begin{align}
    &\begin{aligned}[t]
    \alpha(G)(z_{r}-z_{s})&=(L(G)\mathbf{z})_{r}-(L(G)\mathbf{z})_{s}\\
    &=\left(2+1+|Z\setminus\{a,b\}|\right)(z_{r}-z_{s})-\bigl[(z_{p}-z_{r})+(z_{s}-z_{q})\bigr]\\
    &-(z_{a}-z_{b})-\sum_{v\in A_{1}\setminus\{p,a\}}(z_{v}-z_{r})
     -\sum_{v\in A_{2}\setminus\{q,b\}}(z_{s}-z_{v}).
    \end{aligned} \notag\\[2pt]
    &\text{Then}\quad
    \begin{aligned}[t]
    &\sum_{v\in A_{1}\setminus\{p,a\}}(z_{v}-z_{r})
     +\sum_{v\in A_{2}\setminus\{q,b\}}(z_{s}-z_{v})
     +z_{a}-z_{b}-\frac{1-h}{3}\\
    &=\left(\frac83-\alpha(G)+|Z\setminus\{a,b\}|\right)(1-h)-h\\
    &\geq\left(\frac83-\alpha(G)\right)(1-\alpha(G))-\alpha(G)>\frac78.
    \end{aligned} \notag
\end{align}
Since $S(G)-1<\xi_{n}$, we have
$5S(G)-S(G)^{2}-5\leq3\xi_{n}-\xi_{n}^{2}-1=\frac1{n-1}$.
Combining the preceding inequalities with Eq.~\eqref{eq:alternate-path-energy},
we apply the Cauchy--Schwarz inequality and use $0\leq h\leq S(G)-1<\frac{21}{50}$ to obtain
\begin{align}
    \mathcal E_{G}(\mathbf{z})
    &\geq\sum_{uv\in\{pr,rs,sq\}}(z_{u}-z_{v})^{2}
     +(z_{r}-z_{a})^{2}+(z_{a}-z_{b})^{2}+(z_{b}-z_{s})^{2} \notag\\
    &\quad+\sum_{v\in A_{1}\setminus\{p,a\}}(z_{v}-z_{r})^{2}
     +\sum_{v\in A_{2}\setminus\{q,b\}}(z_{s}-z_{v})^{2} \notag\\
    &\geq (1-z_{r})^{2}+(1-h)^{2}+z_{s}^{2}+\frac{(1-h)^{2}}{3}+\left(z_{a}-z_{b}-\frac{1-h}{3}\right)^{2} \notag\\
    &\quad+\sum_{v\in A_{1}\setminus\{p,a\}}(z_{v}-z_{r})^{2}
     +\sum_{v\in A_{2}\setminus\{q,b\}}(z_{s}-z_{v})^{2}
      \notag\\
    &>\frac{h^{2}}{2}+\frac43(1-h)^{2}+\frac{\bigl(\frac{7}{8}\bigr)^{2}}{n-1} \notag\\
    &>\frac{h^{2}}{2}+\frac43(1-h)^{2}+\frac{3}{4(n-1)} \notag\\
    &\geq\frac{(S(G)-1)^{2}}{2}+\frac43(2-S(G))^{2}+\frac34\bigl(5S(G)-S(G)^{2}-5\bigr). \notag
\end{align}
Combining the above inequality with $\alpha(G)\leq S(G)-1$ and $S(G)<\frac{71}{50}$, we obtain
\begin{align}
    \mathcal E_{G}(\mathbf{z})-\alpha(G)S(G)>\frac{S(G)^{2}-19S(G)+25}{12}\geq\frac{91}{30000}>0, \notag
\end{align}
a contradiction to Eq.~\eqref{eq:diam-three-energy}.

\emph{Case 2.} $|E(F)|\geq2$.
First suppose that $r$ is the center of the star in $F$.
Put $B=\{v:rv\in E(F)\}\subseteq N_{G}(q)$ and $b=|B|\geq2$.
Let $A=N_{G}(p)\setminus\{r\}$, $C=V(G)\setminus(N_{G}[p]\cup N_{G}[q])$ and $D=N_{G}(q)\setminus B$, with $(a,c,d)=(|A|,|C|,|D|)$.
Thus $n=a+b+c+d+3$.
Then $H(a,b,c,d)$ is a spanning subgraph of $\bar{G}$ (see Figure~\ref{fig:H0}).

If $a=c=d=0$, then $\bar{G}\in\mathcal H(n-3,1,1,1)$ with ordered blocks $(B,\{p\},\{q\},\{r\})$.
Assume that $a+c+d>0$.

Eq.~\eqref{eq:extreme-sum} gives $\sum_{v\in B}z_{v}\leq z_{r}+\alpha(G)-1$ and hence $z_{r}\geq1-\alpha(G)$.
Since $b\geq2$, we have
\begin{align}
    \sum_{v\in B}(z_{r}-2z_{v})=bz_{r}-2\sum_{v\in B}z_{v}
    \geq(b-2)z_{r}+2(1-\alpha(G))
    \geq b(1-\alpha(G)). \notag
\end{align}
If $\alpha(G)\leq1$, then $b(1-\alpha(G))\geq0$.
Combining the above inequalities with the Cauchy--Schwarz inequality and using $z_{p}=1$ and $z_{q}=0$, we obtain
\begin{equation}\label{eq:star-energy}
\begin{aligned}
\mathcal E_{G}(\mathbf{z})
&\geq(z_{p}-z_{r})^{2}+\sum_{v\in B}\bigl((z_{r}-z_{v})^{2}+(z_{v}-z_{q})^{2}\bigr)\\
&=(1-z_{r})^{2}+\sum_{v\in B}\bigl((z_{r}-z_{v})^{2}+z_{v}^{2}\bigr)\\
&=(1-z_{r})^{2}+\frac b2z_{r}^{2}+\frac12\sum_{v\in B}(z_{r}-2z_{v})^{2}\\
&\geq\frac{\bigl((1-z_{r})+z_{r}\bigr)^{2}}{1+\frac{2}{b}}
 +\frac{\bigl(\sum_{v\in B}(z_{r}-2z_{v})\bigr)^{2}}{2b}\\
&\geq\frac{b}{b+2}+\frac b2(1-\alpha(G))^{2}\qquad(\alpha(G)\leq1).
\end{aligned}
\end{equation}

If $ad>0$, then Lemmas~\ref{lem:subgraph} and~\ref{lem:H0}\textup{(i)} give $\alpha(\bar{G})\geq1$.
So $\alpha(G)<\frac37$.
Eqs.~\eqref{eq:diam-three-energy} and~\eqref{eq:star-energy} for $b\geq2$ give
\begin{align}
    \frac{81}{98}=\frac12+\left(\frac47\right)^{2}<\mathcal E_{G}(\mathbf{z})<\frac{30}{49}, \notag
\end{align}
a contradiction.
Thus $ad=0$.

Since $d(\bar{G})=3$, $G$ has a dominating edge in $E(F)$, say $rv_\ast$ with $v_\ast\in B$.
Then $A\subseteq N_{G}(r)$, $D\subseteq N_{G}(v_\ast)$ and $C\subseteq N_{G}(r)\cup N_{G}(v_\ast)$.
Put $C_{r}=C\cap N_{G}(r)$ (see Figure~\ref{fig:star}).
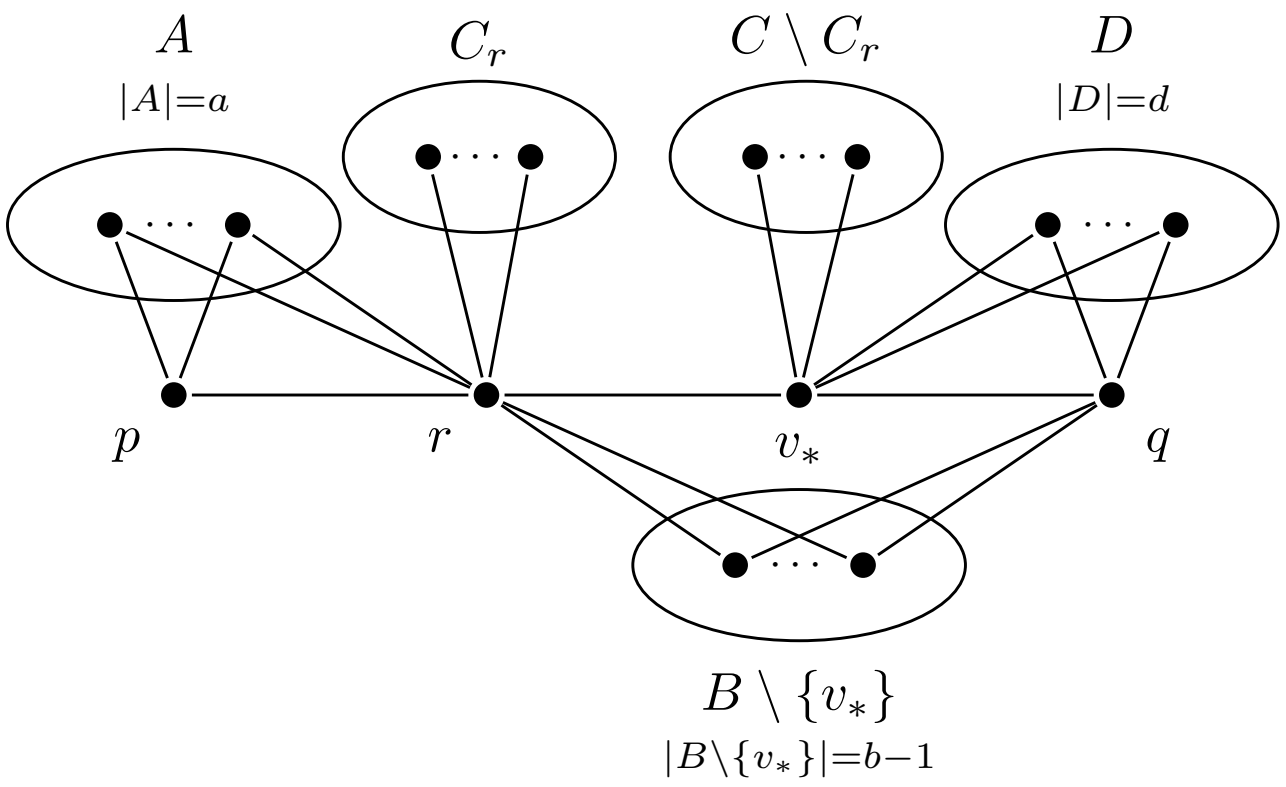
\begin{figure}[H]
\centering
\begin{tikzpicture}[papergraph,xscale=.94,yscale=.75]
\node[vertex,label=below left:$p$] (p) at (-3.3,0) {};
\node[vertex,label=below left:$r$] (r) at (-1.1,0) {};
\node[vertex,label=below:$v_\ast$] (vstar) at (1.1,0) {};
\node[vertex,label=below right:$q$] (q) at (3.3,0) {};
\node[block,minimum width=22mm,minimum height=10mm,inner sep=0pt] (A) at (-3.3,1.5) {};
\node[vertex] (A1) at (-3.75,1.5) {};
\node[vertex] (A2) at (-2.85,1.5) {};
\node[font=\scriptsize,inner sep=0pt] at (-3.3,1.5) {$\cdots$};
\node[anchor=south,align=center] at ($(A.north)+(0,.06)$) {$A$\\[-1pt]$\scriptstyle |A|=a$};
\node[block,minimum width=22mm,minimum height=10mm,inner sep=0pt] (D) at (3.3,1.5) {};
\node[vertex] (D1) at (2.85,1.5) {};
\node[vertex] (D2) at (3.75,1.5) {};
\node[font=\scriptsize,inner sep=0pt] at (3.3,1.5) {$\cdots$};
\node[anchor=south,align=center] at ($(D.north)+(0,.06)$) {$D$\\[-1pt]$\scriptstyle |D|=d$};
\node[block,minimum width=18mm,minimum height=10mm,inner sep=0pt] (Cr) at (-1.15,2.1) {};
\node[vertex] (Cr1) at (-1.51,2.1) {};
\node[vertex] (Cr2) at (-.79,2.1) {};
\node[font=\scriptsize,inner sep=0pt] at (-1.15,2.1) {$\cdots$};
\node[anchor=south,inner sep=1pt] at ($(Cr.north)+(0,.06)$) {$C_{r}$};
\node[block,minimum width=18mm,minimum height=10mm,inner sep=0pt] (Crest) at (1.15,2.1) {};
\node[vertex] (Crest1) at (.79,2.1) {};
\node[vertex] (Crest2) at (1.51,2.1) {};
\node[font=\scriptsize,inner sep=0pt] at (1.15,2.1) {$\cdots$};
\node[anchor=south,inner sep=1pt] at ($(Crest.north)+(0,.06)$) {$C\setminus C_{r}$};
\node[block,minimum width=22mm,minimum height=10mm,inner sep=0pt] (B) at (1.1,-1.5) {};
\node[vertex] (B1) at (0.65,-1.5) {};
\node[vertex] (B2) at (1.55,-1.5) {};
\node[font=\scriptsize,inner sep=0pt] at (1.1,-1.5) {$\cdots$};
\node[anchor=north,align=center] at ($(B.south)+(0,-.06)$) {$B\setminus\{v_\ast\}$\\[-1pt]$\scriptstyle |B\setminus\{v_\ast\}|=b-1$};
\draw (p)--(r);
\draw (r)--(vstar);
\draw (vstar)--(q);
\draw (p)--(A1);
\draw (p)--(A2);
\draw (r)--(A1);
\draw (r)--(A2);
\draw (q)--(D1);
\draw (q)--(D2);
\draw (vstar)--(D1);
\draw (vstar)--(D2);
\draw (r)--(B1);
\draw (r)--(B2);
\draw (q)--(B1);
\draw (q)--(B2);
\draw (r)--(Cr1);
\draw (r)--(Cr2);
\draw (vstar)--(Crest1);
\draw (vstar)--(Crest2);
\foreach \v in {p,r,vstar,q,A1,A2,D1,D2,Cr1,Cr2,Crest1,Crest2,B1,B2} {\fill (\v) circle[radius=1.7pt];}
\end{tikzpicture}
\caption{A spanning subgraph of $G$.}
\label{fig:star}
\end{figure}

Since $z_{q}=0$ and $N_{G}(q)=B\cup D$, Eq.~\eqref{eq:vertex-equations} gives
\begin{equation}\label{eq:case2-vertex-equations}
(L(G)\mathbf{z})_{u}=(L(G)\mathbf{z})_{q}+\alpha(G)z_{u}=\alpha(G)z_{u}-\sum_{v\in B\cup D}z_{v}
\qquad(u\in V(G)).
\end{equation}
Evaluating Eq.~\eqref{eq:case2-vertex-equations} at $p$, where $z_{p}=1$ and $N_{G}(p)=\{r\}\cup A$, gives
\begin{equation}\label{eq:case2-p-identity}
z_{r}-\sum_{v\in B\cup D}z_{v}=1-\alpha(G)+\sum_{u\in A}(1-z_{u}).
\end{equation}

If $a=0$, then $\alpha(G)<1$.
Indeed, if $a=0$, then $N_{G}(p)=\{r\}$ and $rq\notin E(G)$.
So $d_{G}(r)\leq n-2$.
Take $g_{p}=n-2$, $g_{r}=0$ and $g_{v}=-1$ for $v\notin\{p,r\}$.
Since $\mathbf{g}\perp\one$, Eq.~\eqref{eq:variational} gives
\begin{align}
    \alpha(G)\leq\frac{\mathcal E_{G}(\mathbf{g})}{\|\mathbf{g}\|^{2}}
    =\frac{(n-2)^{2}+d_{G}(r)-1}{(n-2)^{2}+n-2}<1. \notag
\end{align}

It suffices to consider $a+c+d>0$ and $ad=0$.
Without loss of generality, we distinguish the following four subcases.

\emph{Subcase 2.1.} $a=0$ and $C_{r}=\varnothing$.
Since $r\in N_{G}(v)$ for $v\in B$ and $z_{u}\geq0$ for $u\in V(G)$, Eq.~\eqref{eq:case2-vertex-equations} at $v$ and Eq.~\eqref{eq:case2-p-identity} with $A=\varnothing$ give
\begin{align}
    (d_{G}(v)-\alpha(G))z_{v}
    =\sum_{u\in N_{G}(v)}z_{u}-\sum_{u\in B\cup D}z_{u}
    \geq z_{r}-\sum_{u\in B\cup D}z_{u}=1-\alpha(G)>0. \notag
\end{align}
Hence $z_{v}>0$ for every $v\in B$.

Since $N_{G}(p)=\{r\}$ and $N_{G}(r)=\{p\}\cup B$, {we have $N_{G}(v)\setminus(D\cup C)\subseteq B\cup\{q\}$ for every $v\in D\cup C$.}
Each vertex of $D\cup C$ is adjacent to $v_\ast\in B$ since $rv_\ast$ dominates $G$.
Choose any ${t_{0}}\in D\cup C$ with $z_{t_{0}}=\max_{{t}\in D\cup C}z_{t}$.
By Eq.~\eqref{eq:case2-vertex-equations} at ${t_{0}}$, with $z_{q}=0$ and $\alpha(G)<1$, we have
\begin{align}
    0\leq(1-\alpha(G))z_{t_{0}}
    &\leq\sum_{v\in N_{G}({t_{0}})\cap(D\cup C)}(z_{t_{0}}-z_{v})+\bigl(|N_{G}({t_{0}})\setminus(D\cup C)|-\alpha(G)\bigr)z_{t_{0}} \notag\\
    &=(L(G)\mathbf{z})_{t_{0}}-\alpha(G)z_{t_{0}}+\sum_{v\in N_{G}(t_{0})\cap B}z_{v} \notag\\
    &=\sum_{v\in N_{G}({t_{0}})\cap B}z_{v}-\sum_{v\in B\cup D}z_{v} \notag\\
    &\leq-\sum_{v\in B\setminus N_{G}({t_{0}})}z_{v}\leq0. \notag
\end{align}
Thus $0\leq z_{t}\leq z_{t_{0}}=0$ for every ${t}\in D\cup C$.
The same inequalities with ${t_{0}=t}$ give $\sum_{v\in B\setminus N_{G}({t})}z_{v}=0$, so $B\subseteq N_{G}({t})$ since $z_{v}>0$ for $v\in B$.
Together with $B\subseteq N_{G}(q)$, this gives $G\in\mathcal H(c+d+1,b,1,1)$ with ordered blocks $(\{q\}\cup D\cup C,B,\{r\},\{p\})$.

In the remaining subcases, Eq.~\eqref{eq:case2-vertex-equations} at $r$ and Eq.~\eqref{eq:case2-p-identity} give
\begin{align}
    (L(G)\mathbf{z})_{r}=\alpha(G)z_{r}-\sum_{v\in B\cup D}z_{v}=(1-\alpha(G))(1-z_{r})+\sum_{u\in A}(1-z_{u}). \notag
\end{align}
Eq.~\eqref{eq:vertex-equations} at $r$ gives $n(L(G)\mathbf{z})_{r}=\alpha(G)\sum_{v\in V(G)}(z_{r}-z_{v})$.
Combining the preceding two equalities with $(L(G)\mathbf{z})_{r}=\sum_{v\in N_{G}(r)}(z_{r}-z_{v})$, we obtain
\begin{equation}\label{eq:star-center-identity}
\begin{aligned}
(n-\alpha(G))\left((1-\alpha(G))(1-z_{r})+\sum_{u\in A}(1-z_{u})\right)
&=(n-\alpha(G))(L(G)\mathbf{z})_{r}\\
&=\alpha(G)\sum_{v\in V(G)\setminus N_{G}[r]}(z_{r}-z_{v}).
\end{aligned}
\end{equation}

Put $m_{1}=a+|C_{r}|$ and $m_{2}=d+|C\setminus C_{r}|$, so $n=b+m_{1}+m_{2}+3$.
Since $A\cup C_{r}\ne\varnothing$ in the remaining subcases, we have $m_{1}\geq1$.
Using $z_{p}=1$, $N_{G}(r)=\{p\}\cup A\cup B\cup C_{r}$ and $0\leq z_{v}\leq1$ for $v\in V(G)$ in Eq.~\eqref{eq:case2-vertex-equations} at $r$, we obtain
\begin{align}
    (b-\alpha(G))z_{r}=1-z_{r}+\sum_{u\in A\cup C_{r}}(z_{u}-z_{r})-\sum_{v\in D}z_{v}
    \leq(m_{1}+1)(1-z_{r}). \notag
\end{align}

For $\alpha(G)<1$, Eq.~\eqref{eq:case2-p-identity} and the above inequality give $0<1-\alpha(G)\leq z_{r}<1$, since $b\geq2$.
Eq.~\eqref{eq:star-center-identity} and $|V(G)\setminus N_{G}[r]|=m_{2}+1$ then give
\begin{align}
    (1-\alpha(G))(n-\alpha(G))(1-z_{r})
    \leq\alpha(G)\sum_{v\in V(G)\setminus N_{G}[r]}(z_{r}-z_{v})
    \leq\alpha(G)(m_{2}+1)z_{r}.\notag
\end{align}
Multiplying the two preceding inequalities and dividing by $z_{r}(1-z_{r})>0$, we obtain
\begin{equation}\label{eq:star-product}
(b-\alpha(G))(1-\alpha(G))(n-\alpha(G))\leq\alpha(G)(m_{1}+1)(m_{2}+1)\qquad(\alpha(G)<1).
\end{equation}

\noindent\hypertarget{claim:large-b}{\textbf{Claim.}}
If $S(G)<\frac{10}{7}$ and $\bar{G}$ contains a spanning copy of either $H(a_{1},b,c_{1},0)$ or $H(0,b,c_{1},a_{1})$, with $a_{1},c_{1}\geq1$,
then $b\in\{2,3,4\}$ and $\alpha(G)<\ell_{b}$, where $\ell_{2}=\frac47$ and $\ell_{3}=\ell_{4}=\frac{19}{28}$.

\begin{proof}[Proof of the claim]
For $b\geq5$, Lemmas~\ref{lem:subgraph} and~\ref{lem:H0}\textup{(iv)} give $\alpha(\bar{G})>\frac35+\frac1{b+2}$.
Together with $S(G)<\frac{10}{7}$, this gives $\alpha(G)<\frac{29}{35}-\frac1{b+2}$.
If $b\geq6$, then Eq.~\eqref{eq:star-energy} gives
\begin{align}
    \mathcal E_{G}(\mathbf{z})
    >\frac{b}{b+2}+\frac b2\left(\frac6{35}+\frac1{b+2}\right)^{2}
    \geq\frac{47b}{35(b+2)}>\frac{10}{7}\left(\frac{29}{35}-\frac1{b+2}\right)
    >\frac{10}{7}\alpha(G), \notag
\end{align}
a contradiction to Eq.~\eqref{eq:diam-three-energy}.

If $b=5$ and $a_{1}+c_{1}\geq3$, then it suffices to check $W_{a_{1},5,c_{1},0}(\frac34)>0$ at $(a_{1},c_{1})=(1,2)$ and $(2,1)$, since $W_{a_{1},5,c_{1},0}(\frac34)$ is increasing in $a_{1}$ and $c_{1}$.
Since $W_{1,5,2,0}(\frac34)=\frac{1289}{1600}>0$ and $W_{2,5,1,0}(\frac34)=\frac{20221}{9280}>0$, Lemmas~\ref{lem:subgraph} and~\ref{lem:H0}\textup{(ii)} give $\alpha(\bar{G})>\frac34$.
Thus $\alpha(G)<\frac{19}{28}$, so Eq.~\eqref{eq:star-energy} gives
\begin{align}
    \mathcal E_{G}(\mathbf{z})>\frac57+\frac52\left(\frac9{28}\right)^{2}
    =\frac{95}{98}+\frac{5}{1568}>\frac{10}{7}\alpha(G), \notag
\end{align}
a contradiction to Eq.~\eqref{eq:diam-three-energy}.

Suppose that $b\in\{2,5\}$ and $a_{1}=c_{1}=1$.
Then $n=b+5$, so $m_{1}+m_{2}=2$ and $(m_{1}+1)(m_{2}+1)\leq4$.
If $b=2$, then Lemmas~\ref{lem:subgraph} and~\ref{lem:H0}\textup{(v)} give $\alpha(\bar{G})>\frac45$, so $\alpha(G)<\frac{22}{35}<\frac23$.
If $b=5$, then $S(G)<\frac{10}{7}$ and $\alpha(\bar{G})>\frac{26}{35}$ give $\alpha(G)<\frac{24}{35}<\frac56$.
Since $0<\alpha(G)<\frac{b}{b+1}$ and $n-\alpha(G)>6$, Eq.~\eqref{eq:star-product} gives
\begin{align}
    4\leq\frac{6b}{b+1}<\frac{(b-\alpha(G))(1-\alpha(G))(n-\alpha(G))}{\alpha(G)}
    \leq(m_{1}+1)(m_{2}+1)\leq4, \notag
\end{align}
a contradiction.

Thus $b\in\{2,3,4\}$, with $a_{1}+c_{1}\geq3$ if $b=2$.
Lemmas~\ref{lem:subgraph} and~\ref{lem:H0}\textup{(v)} give $\alpha(\bar{G})>\frac67$ for $b=2$ and $\alpha(\bar{G})>\frac34$ for $b\in\{3,4\}$.
Together with $S(G)<\frac{10}{7}$, these inequalities give $\alpha(G)<\ell_{b}$.
\end{proof}

If the claim applies, then $b\in\{2,3,4\}$ and $0<\alpha(G)<\ell_{b}<1$.
By Eq.~\eqref{eq:star-product},
\begin{equation}\label{eq:star-parameter-bound}
\begin{aligned}
\frac{15}{14}(n-\ell_{b})
&\leq\frac{(b-\ell_{b})(1-\ell_{b})(n-\ell_{b})}{\ell_{b}}\\
&<\frac{(b-\alpha(G))(1-\alpha(G))(n-\alpha(G))}{\alpha(G)}
\leq(m_{1}+1)(m_{2}+1).
\end{aligned}
\end{equation}

\emph{Subcase 2.2.} $a=0$, $C_{r}\ne\varnothing$, and $m_{2}=d+|C\setminus C_{r}|\geq1$.
Since $C\setminus C_{r}\subseteq N_{\bar{G}}(p)\cap N_{\bar{G}}(r)$, $\bar{G}$ contains a spanning copy of $H(0,b,m_{1},m_{2})$ with $(C,D)$ replaced by $(C_{r},D\cup(C\setminus C_{r}))$.
By claim, $b\in\{2,3,4\}$ and $\alpha(G)<\ell_{b}$.

If $\min\{m_{1},m_{2}\}=1$, then $(m_{1}+1)(m_{2}+1)=2(m_{1}+m_{2})$, so Eq.~\eqref{eq:star-parameter-bound} gives $13(m_{1}+m_{2})>15(b+3-\ell_{b})\geq\frac{465}{7}>65$.
Thus $m_{1}\geq k_{1}$ and $m_{2}\geq k_{2}$ for some $(k_{1},k_{2})\in\{(1,5),(2,2),(5,1)\}$.

Take $\tau=1-\frac{b}{24}$.
Then $\gamma_{0}=\frac{b}{24}+\frac{bm_{1}}{24+b}$ and
\begin{align}
    \frac{24(2-\tau+m_{2})}{b}W_{0,b,m_{1},m_{2}}(\tau)
    &=\gamma_{0}\left[(2-\tau+m_{2})\left(\frac{25b}{24}+m_{1}+m_{2}\right)+2(2-\tau)\right] \notag\\
    &\quad+\frac{b(b+27m_{2})}{24}+m_{2}(2m_{1}+m_{2}+2)-24. \notag
\end{align}
The right-hand side increases with $b$, $m_{1}$, and $m_{2}$.
For $b=2$ and $(m_{1},m_{2})\in\{(1,5),(2,2),(5,1)\}$, direct calculation shows that the right-hand side is positive.
So $W_{0,b,m_{1},m_{2}}(\tau)>0$, and Lemmas~\ref{lem:subgraph} and~\ref{lem:H0}\textup{(ii)} give $\alpha(\bar{G})>\tau$.

Since $S(G)<\frac{10}{7}$, we have $\alpha(G)<\frac{10}{7}-\tau=\frac37+\frac{b}{24}<1$.
By Eqs.~\eqref{eq:diam-three-energy} and~\eqref{eq:star-energy},
\begin{align}
    0&>\mathcal E_{G}(\mathbf{z})-\frac{10}{7}\alpha(G) \notag\\
    &>\frac{b}{b+2}+\frac b2\left(\frac47-\frac b{24}\right)^{2}
    -\frac{10}{7}\left(\frac37+\frac b{24}\right)=f(b). \notag
\end{align}
For $2\leq b\leq4$,
\begin{align}
    f'(b)=\frac{2}{(b+2)^{2}}+\frac12\left(\frac47-\frac b{24}\right)\left(\frac47-\frac b8\right)-\frac5{84}
    \geq\frac{37}{3528}>0. \notag
\end{align}
Thus $f(b)\geq f(2)=\frac1{144}$, a contradiction.

\emph{Subcase 2.3.} $a=0$, $C_{r}\ne\varnothing$, and $m_{2}=d+|C\setminus C_{r}|=0$.
Then $d=0$ and $C_{r}=C$. 
Thus $N_{G}(r)=\{p\}\cup B\cup C$.
If $E_{G}(C,B)=\varnothing$, then $G\in\mathcal H(c+1,1,b,1)$ with ordered blocks $(\{p\}\cup C,\{r\},B,\{q\})$.

Assume that $E_{G}(C,B)\ne\varnothing$.
For $v\in C$, Eq.~\eqref{eq:case2-vertex-equations} gives
\begin{align}
    \bigl(1+|N_{G}(v)\cap B|-\alpha(G)\bigr)z_{v}
    +\sum_{u\in N_{G}(v)\cap C}(z_{v}-z_{u})
    =z_{r}-\sum_{u\in B\setminus N_{G}(v)}z_{u}. \notag
\end{align}
If $N_{G}(v)\cap B\ne\varnothing$ for all $v\in C$, then summing over $v\in C$ gives
\begin{align}
    (2-\alpha(G))\sum_{v\in C}z_{v}
    \leq\sum_{v\in C}\bigl(1+|N_{G}(v)\cap B|-\alpha(G)\bigr)z_{v}
    =cz_{r}-\sum_{v\in C}\sum_{u\in B\setminus N_{G}(v)}z_{u}\leq cz_{r}. \notag
\end{align}
Since $a=0$, we have $\alpha(G)<1$, and Eq.~\eqref{eq:case2-p-identity} gives $z_{r}\geq1-\alpha(G)>0$.
So $\sum_{v\in C}z_{v}<cz_{r}$.
By Eq.~\eqref{eq:case2-vertex-equations} at $r$ and $z_{p}=1$,
\begin{align}
    1=(b+1-\alpha(G))z_{r}+\sum_{v\in C}(z_{r}-z_{v})>bz_{r}. \notag
\end{align}
Lemmas~\ref{lem:subgraph} and~\ref{lem:H0}\textup{(iii)} applied to the spanning copy of $H(0,b,c,0)$ in $\bar{G}$ give $\alpha(\bar{G})>\frac37+\frac1{2(b+1)}$.
Since $S(G)<\frac{10}{7}$, we have $\alpha(G)<1-\frac1{2(b+1)}$.
Here $A=\varnothing$, $V(G)\setminus N_{G}[r]=\{q\}$ and $z_{q}=0$.
Using $bz_{r}<1$ and $n-\alpha(G)=b+c+3-\alpha(G)>b+3$ in Eq.~\eqref{eq:star-center-identity}, we obtain
\begin{align}
    2b+1>\frac{\alpha(G)}{1-\alpha(G)}
    =(n-\alpha(G))\frac{1-z_{r}}{z_{r}}
    >(b+3)(b-1)\geq2b+1, \notag
\end{align}
a contradiction.

Otherwise, choose $u\in C$ with $N_{G}(u)\cap B=\varnothing$.
Then $B\cup\{q\}\subseteq N_{\bar{G}}(u)$, and $c\geq2$ since $E_{G}(C,B)\ne\varnothing$.
Thus replacing $(A,C)$ by $(\{u\},C\setminus\{u\})$ gives a spanning copy of $H(1,b,c-1,0)$ in $\bar{G}$.
By the claim and Eq.~\eqref{eq:star-parameter-bound} with $(m_{1},m_{2})=(c,0)$,
\begin{align}
    \frac{15}{14}<\frac{c+1}{b+c+3-\ell_{b}}<1, \notag
\end{align}
a contradiction.

\emph{Subcase 2.4.} $a\geq1$ and $d=0$.
If $c=0$, then $G\in\mathcal H(a+1,1,b,1)$ with ordered blocks $(\{p\}\cup A,\{r\},B,\{q\})$.
For $c\geq1$, applying claim to $H(a,b,c,0)\subseteq\bar{G}$ gives $b\in\{2,3,4\}$ and $\alpha(G)<\ell_{b}$.

Choose $u\in A\subseteq N_{G}(r)$ and let $F_{0}$ be the star with edge set $\{ru\}\cup\{rv:v\in B\}$.
Remark~\ref{rem:energy-components}\textup{(ii)} gives $\kappa_{F_{0}}=\frac{3b}{2(b+3)}$.

For $b=2,3,4$, we have $|F_{0}|=b+2$ and $\alpha(G)<\ell_{b}<2\kappa_{F_{0}}$.
By Lemma~\ref{lem:energy-bound}\textup{(ii)} and Eq.~\eqref{eq:diam-three-energy},
\begin{align}
    \frac{10}{7}\ell_{b}
    >\mathcal E_{G}(\mathbf{z})
    >\kappa_{F_{0}}+\frac{(2\kappa_{F_{0}}-\ell_{b})^{2}}{b+2-4\kappa_{F_{0}}}
    >\frac{10}{7}\ell_{b}. \notag
\end{align}

Finally, suppose that $s$ is the center of the star in $F$.
Put $\mathbf{x}'=-\mathbf{x}$.
Then $x_{q}'=\max_{v} x_{v}'$ and $x_{p}'=\min_{v} x_{v}'$, with $x_{q}'-x_{p}'=x_{p}-x_{q}$.
We have $\frac{\mathbf{x}'-x_{p}'\one}{x_{q}'-x_{p}'}=\one-\mathbf{z}$ and $\mathcal E_{G}(\one-\mathbf{z})=\mathcal E_{G}(\mathbf{z})$.
Thus the preceding argument with $\mathbf{x}'$ and the path $qsrp$ gives the result.
\end{proof}

\begin{proof}[Proof of Theorem~\ref{thm:main}]
We may assume that $S(G)\leq n+a_{n}-c_{n}$.
By Lemma~\ref{lem:diameter-two}\textup{(ii)}, we have $d(G),d(\bar{G})\geq3$.
Thus each graph has a dominating edge, so $d(G)=d(\bar{G})=3$.
Let $\mathbf{x},\mathbf{y}$ be Fiedler vectors of $G,\bar{G}$, respectively, and choose $p,q\in V(G)$ with $x_{p}=\max_{v\in V(G)}x_{v}$ and $x_{q}=\min_{v\in V(G)}x_{v}$.
Without loss of generality, assume that $x_{p}-x_{q}\geq\max_{u,v\in V(G)}|y_{u}-y_{v}|$.
Then Lemma~\ref{lem:diameter-two}\textup{(i)} gives $d_{G}(p,q)=3$.

Lemmas~\ref{cor:Mn}, \ref{lem:four-block} and~\ref{lem:diameter-three}\textup{(i)}, together with Corollary~\ref{cor:bull}, give $S(G)\geq n+a_{n}-c_{n}$, with equality if and only if $G$ or $\bar{G}$ is isomorphic to a graph $\mathcal D(H)$, where $H$ is an arbitrary graph of order $n-3$.

Let $n\geq6$.
Then $p_{n}(1)>0>p_{n}(n-1)$.
For $n=6$, direct calculation gives $p_{6}(\frac{12}{25})<0<p_{6}(\frac{51}{10})$, so $a_{6}>\frac{12}{25}$ and $c_{6}<\frac{51}{10}$.
Hence $6+a_{6}-c_{6}>1.38$.
For $n\geq7$, put $d=\frac{1}{(n-3)^{2}}$.
Direct calculation gives $p_{n}(\frac{19}{50}+d)<0<p_{n}(n-1+d)$, so $a_{n}>\frac{19}{50}+d$ and $c_{n}<n-1+d$.
Hence $n+a_{n}-c_{n}>1.38$.
\end{proof}
\FloatBarrier

\section{Applications}\label{sec:applications}

\subsection{An alternative proof of \texorpdfstring{$S(G)\geq1$}{S(G) >= 1}}\label{sec:sum-application}

\begin{lemma}[{\cite[Lemma~3]{EK2021}}]\label{lem:disconnected-complement}
Let $G$ be a disconnected graph of order $n\geq2$.
Then $\alpha(\bar{G})\geq1$, with equality if and only if $\bar{G}$ is isomorphic to a graph $K_{1}\vee H$, where $H$ is any disconnected graph of order $n-1$.
\end{lemma}

\begin{corollary}[{\cite[Theorem~1]{EK2021}}]\label{cor:sum-all-graphs}
For every graph $G$ of order $n\geq2$, we have $S(G)\geq1$.
Equality holds if and only if $G$ or $\bar{G}$ is isomorphic to a graph $K_{1}\vee H$, where $H$ is any disconnected graph of order $n-1$.
\end{corollary}

\begin{proof}
If one of $G$ and $\bar{G}$ is disconnected, then we may assume without loss of generality that $G$ is disconnected.
Then $S(G)=\alpha(\bar{G})\geq1$ by Lemma~\ref{lem:disconnected-complement}, which also gives the equality condition.

If both graphs are connected, then $2(n-1)\leq\binom{n}{2}$, so $n\geq4$.
By Theorem~\ref{thm:main}, $S(G)\geq n+a_{n}-c_{n}>1$, where the last inequality for $n=4,5$ follows by direct calculation.
\end{proof}

\subsection{The product of the Laplacian spectral radii}\label{sec:laplacian-product}

\begin{lemma}\label{lem:laplacian-product-threshold}
Let $H$ be a graph of order $n-3$, where $n\geq5$.
Then
\begin{align}
    \mu_{1}(\mathcal D(H))\mu_{1}(\overline{\mathcal D(H)})
    =c_{n}(n-a_{n})>\left(n-\frac57\right)^{2}, \notag
\end{align}
where $a_{n}$ and $c_{n}$ are as in Lemma~\ref{lem:dandelion-spectrum}.
\end{lemma}

\begin{proof}
By Lemmas~\ref{lem:dandelion-spectrum} and~\ref{lem:complement-spectrum}, $\mu_{1}(\mathcal D(H))=c_{n}$ and $\mu_{1}(\overline{\mathcal D(H)})=n-a_{n}$.

For $5\leq n\leq10$, direct calculation from $p_{n}$ gives $c_{n}(n-a_{n})>(n-\frac57)^{2}$.

Now suppose that $n\geq11$.
Since $p_{n}(0)<0<p_{n}(1)$ and $p_{n}(n-1)<0<p_{n}(n)$, we have $a_{n}<1$ and $c_{n}>n-1$.
By Lemma~\ref{cor:Mn}, $n+a_{n}-c_{n}<\frac{71}{50}$.
Thus
\begin{align}
    c_{n}(n-a_{n})
    &=(n-1)(c_{n}-a_{n}+1)+(c_{n}-n+1)(1-a_{n}) \notag\\
    &>(n-1)(c_{n}-a_{n}+1)>(n-1)(n-\frac{21}{50}) \notag\\
    &=(n-\frac{5}{7})^{2}+\frac{21n-221}{2450}>(n-\frac{5}{7})^{2}. \notag
\end{align}
\end{proof}

\begin{proof}[Proof of Theorem~\ref{thm:laplacian-product}]
For $n=4$, we have $G\cong P_{4}=\mathcal D(K_{1})$ since both $G$ and $\bar{G}$ are connected, so Lemmas~\ref{lem:complement-spectrum} and~\ref{lem:dandelion-spectrum} give $\mu_{1}(G)\mu_{1}(\bar{G})=c_{4}(4-a_{4})$.

Suppose that $n\geq5$ and $\mu_{1}(G)\mu_{1}(\bar{G})\geq c_{n}(n-a_{n})$.
By Lemma~\ref{lem:complement-spectrum}, $\mu_{1}(G)=n-\alpha(\bar{G})$ and $\mu_{1}(\bar{G})=n-\alpha(G)$.
Lemma~\ref{lem:laplacian-product-threshold} gives
\begin{align}
    n-\frac57<\sqrt{c_{n}(n-a_{n})}\leq\sqrt{\mu_{1}(G)\mu_{1}(\bar{G})}
    \leq\frac{\mu_{1}(G)+\mu_{1}(\bar{G})}{2}=n-\frac12 S(G). \notag
\end{align}
So $S(G)<\frac{10}{7}$.

By Theorem~\ref{thm:main}, $S(G)\geq n+a_{n}-c_{n}$.
Since $n-\alpha(G)>0$, we have
\begin{align}
    0&\leq\mu_{1}(G)\mu_{1}(\bar{G})-c_{n}(n-a_{n}) \notag\\
    &=(n-\alpha(G))(n-S(G)+\alpha(G))-c_{n}(n-a_{n}) \notag\\
    &\leq(n-\alpha(G))(c_{n}-a_{n}+\alpha(G))-c_{n}(n-a_{n}) \notag\\
    &=-(\alpha(G)-a_{n})(\alpha(G)-n+c_{n}). \notag
\end{align}
By symmetry, this gives $\max\{\alpha(G),\alpha(\bar{G})\}\leq\max\{a_{n},n-c_{n}\}$.
Since $p_{n}(0)<0<p_{n}(1)$ and $p_{n}(n-1)<0<p_{n}(n)$, we have $\max\{a_{n},n-c_{n}\}<1$.

Let $\mathbf{x},\mathbf{y}$ be Fiedler vectors of $G,\bar{G}$, respectively, and choose $p,q\in V(G)$ such that $x_{p}=\max_{v\in V(G)}x_{v}$ and $x_{q}=\min_{v\in V(G)}x_{v}$.
Without loss of generality, assume that $x_{p}-x_{q}\geq\max_{u,v\in V(G)}|y_{u}-y_{v}|$.
Since $\alpha(G),\alpha(\bar{G})<1$, Lemma~\ref{lem:normalization}\textup{(i)} gives $d_{G}(p,q)\geq3$ and $d(\bar{G})\geq3$.
Thus  $d(G)=d(\bar{G})=d_{G}(p,q)=3$.
Since $S(G)<\frac{10}{7}$ and $\alpha(\bar{G})<1$, Lemma~\ref{lem:diameter-three}\textup{(ii)} implies that $G$ or $\bar{G}$ belongs to $\mathcal H(m_{1},m_{2},1,1)\cup\mathcal H_{n}$, where $m_{1}+m_{2}=n-2$.
Without loss of generality, assume that $G$ belongs to this union.

By the proof of Lemma~\ref{lem:four-block},
\begin{align}
    \psi(\lambda)&=\lambda-1+\frac{n-2\lambda}{\lambda^{2}-n\lambda+n}, \notag
\end{align}
with $\psi'(\lambda)>0$ on $(0,1)$ and $\psi((0,1))=(0,n-2)$.
By Eq.~\eqref{eq:comparison-relations} and the proof of Corollary~\ref{cor:bull}, we have $\psi(\alpha(G))=s$ and $\psi(\alpha(\bar{G}))=n-2-s$ for some $s\in[1,n-3]$, with $s=\frac{n-2}{2}$ if $G\in\mathcal H_{n}$.
Then $\mu_{1}(G)\mu_{1}(\bar{G})=P(s)$, where
\begin{align}
    P(x)=\bigl(n-\psi^{-1}(x)\bigr)\bigl(n-\psi^{-1}(n-2-x)\bigr)\qquad (1\leq x\leq n-3). \notag
\end{align}

Put $q_{n}(\lambda)=\lambda^{2}-n\lambda+n$.
Since $0<q_{n}(\lambda)\leq n-2\lambda<n$ for $0<\lambda<1$, we have
\begin{align}
    \frac{d}{d\lambda}\bigl((n-\lambda)\psi'(\lambda)\bigr)
    &=-1+\frac{2}{q_{n}(\lambda)}+\frac{n(n-2\lambda)}{q_{n}(\lambda)^{2}}
    +\frac{2n(n-4)(n-\lambda)(n-2\lambda)}{q_{n}(\lambda)^{3}}>0. \notag
\end{align}

By Lemmas~\ref{lem:spectral-parameter}\textup{(ii)} and~\ref{lem:laplacian-product-threshold}, we have $P(s)\leq P(1)=P(n-3)=c_{n}(n-a_{n})$, with equality if and only if $s=1$ or $s=n-3$.
Then $G$ or $\bar{G}$ is isomorphic to a graph $\mathcal D(H)$, where $H$ is an arbitrary graph of order $n-3$.
\end{proof}

\subsection{The algebraic connectivity of complements of trees}\label{sec:tree-complements}

\begin{lemma}[{\cite[4.4, p.~304]{Fiedler1973}, \cite[Corollaries~4.4 and~4.5]{GMS1990}}]\label{lem:path-connectivity}
\leavevmode
\begin{enumerate}
\item[\textup{(i)}] For $n\geq2$, $\alpha(P_{n})=2\bigl(1-\cos\frac{\pi}{n}\bigr)$.
For $n\geq3$, $\alpha(K_{1,n-1})=1$.

\item[\textup{(ii)}] If $T$ is a tree of diameter $d\geq1$, then $\alpha(T)\leq\alpha(P_{d+1})$.
If $T$ is a nonstar tree of order $n\geq6$, then $\alpha(T)<0.49$.
\end{enumerate}
\end{lemma}

For trees $T$ of order $n$ and fixed maximum degree $\Delta(T)\geq3$, Kolokolnikov~\cite{Kolokolnikov2015} proved that $\alpha(T)\leq \frac{2(\Delta(T)-2)}{n}+O\!\left(\frac{\ln n}{n^{2}}\right)$ as $n\to\infty$.
The following lemma gives an explicit upper bound on $\alpha(T)$ for every nonstar tree of order $n\geq4$.

\begin{lemma}\label{lem:tree-max-degree}
Let $T$ be a nonstar tree of order $n\geq4$.
Then
\begin{align}
    \alpha(T)\leq\frac{\Delta(T)-(3\sqrt{2}-4)}{n-1}. \notag
\end{align}
Equality holds if and only if $T\cong P_{4}$.
\end{lemma}

\begin{proof}
Put $\Delta=\Delta(T)$.
For $n=4,5$, the result follows by direct calculation.
Now suppose that $n\geq6$.
Since $3\sqrt{2}-4<\frac14$, it suffices to prove $\Delta-(n-1)\alpha(T)>\frac14$.

If $\Delta=2$, then $T\cong P_{n}$, so Lemma~\ref{lem:path-connectivity}\textup{(i)} gives $\alpha(T)<\frac{\pi^{2}}{n^{2}}$.
Thus $\Delta-(n-1)\alpha(T)>2-\frac{10(n-1)}{n^{2}}\geq\frac25$.

If $2\Delta\geq n$, then Lemma~\ref{lem:path-connectivity}\textup{(ii)} gives $\alpha(T)<\frac12$.
Thus $\Delta-(n-1)\alpha(T)>\Delta-\frac{n-1}{2}\geq\frac12$.

It remains to consider $3\leq\Delta\leq\frac{n-1}{2}$.
Let $v\in V(T)$ minimize $\sum_{u\in V(T)}d_{T}(u,v)$, $B$ be the vertex set of a largest component of $T-v$, and $w$ be the neighbor of $v$ in $B$.
By the choice of $v$,
\begin{align}
    0\leq\sum_{u\in V(T)}\bigl(d_{T}(u,w)-d_{T}(u,v)\bigr)=n-2|B|. \notag
\end{align}
Since $n-1\leq d_{T}(v)|B|\leq\Delta|B|$, we have $\frac{n-1}{\Delta}\leq|B|\leq\frac n2$.

Define $\mathbf{y}$ by $y_{v}=\frac{n-1-\Delta|B|}{n}$ and
\begin{align}
    y_{u}-y_{v}=
    \begin{cases}
        \Delta-1, & u\in B,\\
        -1, & u\in V(T)\setminus(B\cup\{v\}).
    \end{cases} \notag
\end{align}
Then $\sum_{u\in V(T)}y_{u}=ny_{v}+\Delta|B|-(n-1)=0$.
Since $\frac{n-1}{\Delta}\leq |B|\leq\frac n2$, we have $1-\frac1n-\frac{\Delta}{2}\leq y_{v}\leq0$.
Moreover, $1-\frac1n-\frac{\Delta}{2}>2-\Delta$, so $y_{v}>-(\Delta-2)$.
Since $\mathbf{y}$ is constant on each component of $T-v$ and $\mathbf{y}\perp\one$, we have
\begin{equation}\label{eq:tree-y-bounds}
\begin{aligned}
    \mathcal E_{T}(\mathbf{y})
    &=(\Delta-1)^{2}+d_{T}(v)-1\leq\Delta(\Delta-1),\\
    \|\mathbf{y}\|^{2}
    &=\sum_{u\in V(T)}(y_{u}-y_{v})^{2}-ny_{v}^{2}\\
    &=(\Delta-1)^{2}|B|+n-1-|B|-ny_{v}^{2}\\
    &=(\Delta-1)(n-1)-ny_{v}(\Delta-2+y_{v})\geq(\Delta-1)(n-1).
\end{aligned}
\end{equation}

For $u\sim v$, $(L(T)\mathbf{y})_{u}=y_{u}-y_{v}$, and for $u\notin N_{T}[v]$, $(L(T)\mathbf{y})_{u}=0$.
Using $L(T)\mathbf{y}\perp\one$ and $\mathcal E_{T}(\mathbf{y})=\sum_{u\sim v}(y_{u}-y_{v})^{2}$, we have
\begin{align}
    \sum_{u\sim v}(y_{u}-y_{v})&=-(L(T)\mathbf{y})_{v},\qquad
    (L(T)\mathbf{y})_{v}^{2}=\|L(T)\mathbf{y}\|^{2}-\mathcal E_{T}(\mathbf{y}). \notag
\end{align}
Since $N_{T}(v)$ is independent, we obtain
\begin{align}
    \mathcal E_{T}(L(T)\mathbf{y})
    &=\sum_{u\sim v}\bigl(y_{u}-y_{v}-(L(T)\mathbf{y})_{v}\bigr)^{2}+\sum_{u\sim v}(d_{T}(u)-1)(y_{u}-y_{v})^{2} \notag\\
    &=\sum_{u\sim v}d_{T}(u)(y_{u}-y_{v})^{2}+(d_{T}(v)+2)(L(T)\mathbf{y})_{v}^{2} \notag\\
    &\leq\Delta\mathcal E_{T}(\mathbf{y})+4\Delta(L(T)\mathbf{y})_{v}^{2} \notag\\
    &=4\Delta\|L(T)\mathbf{y}\|^{2}-3\Delta\mathcal E_{T}(\mathbf{y}). \label{eq:tree-laplacian-energy}
\end{align}

Put $\mathbf{z}=\mathbf{y}-\frac{1}{2\Delta}L(T)\mathbf{y}$, so $\mathbf{z}\perp\one$.
Using Eqs.~\eqref{eq:tree-y-bounds} and~\eqref{eq:tree-laplacian-energy}, we obtain
\begin{align}
    \mathcal E_{T}(\mathbf{z})
    &=\mathcal E_{T}(\mathbf{y})-\frac{\|L(T)\mathbf{y}\|^{2}}{\Delta}
    +\frac{\mathcal E_{T}(L(T)\mathbf{y})}{4\Delta^{2}}             \notag\\
    &\leq\left(1-\frac{3}{4\Delta}\right)\mathcal E_{T}(\mathbf{y})
    \leq\left(\Delta-\frac34\right)(\Delta-1),        \notag\\
    \|\mathbf{z}\|^{2}
    &=\|\mathbf{y}\|^{2}-\frac{\mathcal E_{T}(\mathbf{y})}{\Delta}+\frac{\|L(T)\mathbf{y}\|^{2}}{4\Delta^{2}} \notag\\
    &\geq\|\mathbf{y}\|^{2}-\frac{\mathcal E_{T}(\mathbf{y})}{\Delta}\geq(\Delta-1)(n-2)>0.      \notag
\end{align}
By Eq.~\eqref{eq:variational}, $\alpha(T)\leq\frac{\mathcal E_{T}(\mathbf{z})}{\|\mathbf{z}\|^{2}}\leq\frac{\Delta-\frac34}{n-2}$.
Since $2\Delta\leq n-1$, we obtain
\begin{align}
    \Delta-(n-1)\alpha(T)\geq\frac{\frac34(n-1)-\Delta}{n-2}\geq\frac{n-1}{4(n-2)}>\frac14. \notag
\end{align}
\end{proof}

\begin{remark}
Combining Lemmas~\ref{lem:path-connectivity}\textup{(i)} and~\ref{lem:tree-max-degree}, we obtain $\alpha(T)\leq\frac{\Delta(T)}{n-1}$ for every tree $T$ of order $n\geq3$, with equality if and only if $T\cong K_{1,n-1}$.
\end{remark}

\begin{lemma}[{\cite[Proposition~1, pp.~82--83]{Fan2008}}]\label{lem:dandelion}
Let $T$ be the double star obtained from an edge $uv$ by attaching $a$ pendant vertices to $u$ and $n-a-2$ to $v$, where $n\geq4$ and $1\leq a\leq n-3$.
Then $\alpha(T)$ and $\mu_{1}(T)$ are, respectively, the smallest and largest roots of the equation
\begin{align}
    (x-n)(x-1)^{2}+a(n-a-2)x=0. \notag
\end{align}
\end{lemma}

\begin{proof}[Proof of Theorem~\ref{thm:tree-complement}]
Since $d(T)\geq3$, $\bar{T}$ has a dominating edge and is connected.
Put $m=2n-2-\Delta$.
If $\Delta=n-2$, then $T\cong\mathcal D((n-3)K_{1})$ and $m=n$.
Lemmas~\ref{lem:complement-spectrum} and~\ref{lem:dandelion-spectrum} give $\alpha(T)=a_{n}$ and $\alpha(\bar{T})=n-c_{n}$.
By Lemma~\ref{lem:tree-max-degree}, $\frac{\Delta-(3\sqrt{2}-4)}{n-1}\geq a_{n}$, so equality holds.

Now suppose that $\Delta\leq n-3$, so $m>n$ and $n\geq5$.
Since $p_{m}(0)<0<p_{m}(1)$, we have $0<a_{m}<1$.
If $d(T)\geq4$, then Lemma~\ref{lem:path-connectivity}\textup{(i)} and \textup{(ii)} give $\alpha(T)\leq\alpha(P_{5})=\frac{3-\sqrt{5}}{2}$.
Since $p_{m}(a_{m})=0$, we have $3a_{m}-a_{m}^{2}-1=\frac{1}{m-1-a_{m}}>0$, so $\alpha(T)<a_{m}$.

If $d(T)=3$, then $T$ is a double star and $\Delta\geq3$.
By Lemma~\ref{lem:dandelion} with $a=\Delta-1$, $\alpha(T)$ is the smallest root of the polynomial
\begin{align}
    (x-n)(x-1)^{2}+(\Delta-1)(n-\Delta-1)x
    =p_{m}(x)+(m-n)\bigl((1-x)^{2}+(\Delta-3)x\bigr). \notag
\end{align}
This polynomial is negative at $0$ and positive at $a_{m}$, so $\alpha(T)<a_{m}$.

Since $T\not\cong P_{4}$, Lemma~\ref{lem:tree-max-degree} gives $\alpha(T)<\frac{\Delta-(3\sqrt{2}-4)}{n-1}$.
Thus Theorem~\ref{thm:main} gives
\begin{align}
    \alpha(\bar{T})=S(T)-\alpha(T)
    &\geq n+a_{n}-c_{n}-\alpha(T) \notag\\
    &>n+a_{n}-c_{n}-\min\left\{\frac{\Delta-(3\sqrt{2}-4)}{n-1},\,a_{m}\right\}. \notag
\end{align}
\end{proof}

Since $p_{n+1}(0)<0<p_{n+1}(a_{n})=\frac{1}{n-1-a_{n}}$, we have $a_{n+1}<a_{n}$.
Hence Theorem~\ref{thm:tree-complement} refines the lower bound $n-c_{n}$ of Jiang et al.~\cite{JYC2015}.

\section{A Nordhaus--Gaddum problem on Laplacian eigenvalue sums}\label{sec:laplacian-sum-problem}

For a graph $G$ of order $n\geq2$, put $S_{k}(G)=\mu_{k}(G)+\mu_{k}(\bar{G})$ for $1\leq k\leq n-1$.
For $n\geq3$, Grij\'o et al.~\cite[Conjecture~3]{Grijo2019} conjectured that $S_{2}(G)\leq2n-2$.
By Lemma~\ref{lem:complement-spectrum}, $S_{n-2}(G)=2n-S_{2}(G)$, so this upper bound is equivalent to $S_{n-2}(G)\geq2$.
Li and Guo~\cite[Conjecture~1]{LiGuo2024} conjectured that $S_{k}(G)\geq n-k$ for $1\leq k\leq n-1$.
They also conjectured that equality holds if and only if $G$ or $\bar{G}$ is isomorphic to $K_{n-k}\vee H$, where $H$ has order $k$ and at least $n-k+1$ connected components.
For $1\leq k\leq \frac{n}{2}$, Li and Guo~\cite[Theorem~18]{LiGuo2024} proved the stronger bound $S_{k}(G)\geq n$.

The join $K_{n-k}\vee H$ has disconnected complement, so at least one of $G$ and $\bar{G}$ is disconnected in the conjectured equality cases.
When both graphs are connected, Theorem~\ref{thm:main} gives the sharp bound $S_{n-1}(G)\geq n+a_{n}-c_{n}$.

For an integer $k$ with $\frac{n}{2}<k\leq n-2$, let $\mathcal A_{n,k}$ consist of the graphs obtained from $K_{n-k}\vee H$ by attaching a new pendant vertex to each vertex of $K_{n-k}$, where $H$ is an arbitrary graph of order $2k-n$.
These graphs are thin spiders~\cite[p.~3]{CH2024}.
By construction, both $G$ and $\bar{G}$ are connected for every $G\in\mathcal A_{n,k}$.

\begin{lemma}[{\cite[Section~2]{DJ2018}}]\label{lem:thin-spider-spectrum}
Let $G_{0}\in\mathcal A_{n,k}$ be obtained by taking $H$ to be edgeless.
The Laplacian spectrum of $G_{0}$ is the multiset
\begin{align}
    \left\{
    0,\ (n-k)^{[2k-n-1]},\
    \frac{k+2\pm\sqrt{k^{2}-4(n-k-1)}}{2},\
    \left(\frac{k+2\pm\sqrt{k^{2}+4}}{2}\right)^{[n-k-1]}
    \right\}, \notag
\end{align}
where $x^{[s]}$ denotes $s$ copies of $x$ and both signs are taken in each $\pm$ expression.
\end{lemma}

\begin{corollary}\label{cor:pendant-core-sum}
If $G\in\mathcal A_{n,k}$, then
\begin{align}
    S_{k}(G)=n-\sqrt{k^{2}-4(n-k-1)}. \notag
\end{align}
\end{corollary}

\begin{proof}
Let $G_{0}$ be obtained from $G$ by deleting the edges of $H$.
By Lemmas~\ref{lem:quotient} and~\ref{lem:complement-spectrum}, adding the edges of $H$ only replaces one copy of $n-k$ by $n-k+\mu_{j}(H)\in[n-k,k]$ for each $1\leq j\leq2k-n-1$.

Put $\lambda_\pm=\frac{k+2\pm\sqrt{k^{2}-4(n-k-1)}}{2}$.
Since $2\leq n-k<k$, we have
\begin{align}
    (k-2)^{2}<k^{2}-4(n-k-1)<k^{2}<k^{2}+4, \notag
\end{align}
so $0<\lambda_-<2\leq n-k<k<\lambda_+$.
Thus Lemma~\ref{lem:thin-spider-spectrum} gives
\begin{align}
    \mu_{n-k}(G)&=\mu_{n-k}(G_{0})=\lambda_+,\qquad
    \mu_{k}(G)=\mu_{k}(G_{0})=\lambda_-. \notag
\end{align}
By Lemma~\ref{lem:complement-spectrum},
\begin{align}
    S_{k}(G)=n+\mu_{k}(G)-\mu_{n-k}(G)=n-\sqrt{k^{2}-4(n-k-1)}. \notag
\end{align}
\end{proof}

The same formula holds when $\bar{G}\in\mathcal A_{n,k}$.

Let $t\geq3$ be an integer, and put $h=2^{t}-1$ and $d=\frac{h+1}{2}$.
Let $C$ and $U$ be disjoint sets of size $h$, each indexed by the nonzero vectors of $\mathbb F_{2}^{t}$.
Form a bipartite graph with parts $C$ and $U$ by joining $\mathbf{c}\in C$ to $\mathbf{u}\in U$ exactly when $\mathbf{c}^{\mathsf T}\mathbf{u}=1$ over $\mathbb F_{2}$.
This is the incidence graph of a symmetric $2$-$\bigl(h,d,\frac{d}{2}\bigr)$ design.

For each integer $n\geq2h$, add a set $W$ of $n-2h$ new vertices and all edges within $C\cup W$ to obtain $G_{n,h}$~\cite[Sections~3--4]{Anonymous2026}.
When $n=2h$, $W=\varnothing$ and the construction is a special case of that in~\cite[Example~5]{MTR2011}.
By construction, both $G_{n,h}$ and $\overline{G_{n,h}}$ are connected.

\begin{lemma}[{\cite[pp.~4--5]{Anonymous2026}}]\label{lem:design-spectrum}
Let $G=G_{n,h}$ with $n\geq2h$.
Then
\begin{align}
    \mu_{h-1}(G)=\lambda_{+},\qquad \mu_{n-h+1}(G)=\lambda_{-},
    \quad\text{where}\quad
    \lambda_{\pm}=\frac{n-h+2d\pm\sqrt{(n-h)^{2}+2d}}{2}. \notag
\end{align}
\end{lemma}

\begin{corollary}\label{cor:design-laplacian-sum}
Let $G=G_{n,h}$ with $2h\leq n\leq\frac{7h-9}{2}$ and put $k=n-h+1$.
Then
\begin{align}
    \frac{n}{2}<k<\frac{5n}{7},\qquad
    S_{k}(G)<n-\sqrt{k^{2}-4(n-k-1)}. \notag
\end{align}
\end{corollary}

\begin{proof}
From $2h\leq n\leq\frac{7h-9}{2}$, we obtain $2k-n=n-2h+2\geq2$ and $5n-7k=7h-2n-7\geq2$, so $\frac{n}{2}<k<\frac{5n}{7}$.
Since $n-k=h-1$, Lemmas~\ref{lem:design-spectrum} and~\ref{lem:complement-spectrum} give
\begin{align}
    S_{k}(G)
    &=n-\sqrt{(k-1)^{2}+h+1} \notag\\
    &=n-\sqrt{k^{2}-4(n-k-1)+5n-7k-1} \notag\\
    &<n-\sqrt{k^{2}-4(n-k-1)}. \notag
\end{align}
\end{proof}

\begin{remark}
Since $h$ is odd, we may take $n=\frac{7h-9}{2}$ in Corollary~\ref{cor:design-laplacian-sum}.
Then $5n=7k+2$, so $S_{k}(G_{n,h})<n-\sqrt{k^{2}-4(n-k-1)}$
for $k=\lceil\frac{5n}{7}\rceil-1$ and arbitrarily large $n$.
\end{remark}

\begin{problem}\label{prob:connected-laplacian-sum}
For each $n\geq7$ and each integer $k$ with $\lceil\frac{5n}{7}\rceil\leq k\leq n-2$, characterize all graphs $G$ of order $n$ such that both $G$ and $\bar{G}$ are connected and
\begin{align}
    S_{k}(G)\leq n-\sqrt{k^{2}-4(n-k-1)}. \notag
\end{align}
\end{problem}

Using a computer, we check that for $7\leq n\leq10$, the graphs satisfying the conditions in Problem~\ref{prob:connected-laplacian-sum} are precisely those in $\mathcal A_{n,n-2}$ and their complements.

\Needspace{5\baselineskip}


\begin{thebibliography}{99}
\setlength{\itemsep}{0pt}
\bibitem{Afshari2023}
B. Afshari, Algebraic connectivity of the second power of a graph, \emph{Journal of Graph Theory}, \textbf{104}(2) (2023), 275--281.

\bibitem{Anonymous2026}
Anonymous, Nordhaus--Gaddum products of Laplacian eigenvalues: a counterexample and the sharp limiting profile, Zenodo preprint (2026), version 0.1.0-candidate, \href{https://doi.org/10.5281/zenodo.23942304}{doi:10.5281/zenodo.23942304}.

\bibitem{ATR2014}
F. Ashraf and B. Tayfeh-Rezaie, Nordhaus-Gaddum Type Inequalities for Laplacian and Signless Laplacian Eigenvalues, \emph{The Electronic Journal of Combinatorics}, \textbf{21}(3) (2014), \#P3.6.

\bibitem{BTF2009}
Y.-H. Bao, Y.-Y. Tan, and Y.-Z. Fan, The Laplacian spread of unicyclic graphs, \emph{Applied Mathematics Letters}, \textbf{22}(7) (2009), 1011--1015.

\bibitem{BEHK2022}
W. Barrett, E. Evans, H. T. Hall, and M. Kempton, New conjectures on algebraic connectivity and the Laplacian spread of graphs, \emph{Linear Algebra and its Applications}, \textbf{648} (2022), 104--132.

\bibitem{BH2012}
A. E. Brouwer and W. H. Haemers, \emph{Spectra of Graphs}, Springer, New York, 2012.

\bibitem{CGLW2025}
Q. Chen, J.-M. Guo, W.-J. Li, and Z. Wang, A conjecture on the Nordhaus-Gaddum product type inequality for Laplacian eigenvalues of a graph, \emph{The Electronic Journal of Combinatorics}, \textbf{32}(4) (2025), \#P4.37.

\bibitem{CD2016}
X. Chen and K. Ch. Das, Some results on the Laplacian spread of a graph, \emph{Linear Algebra and its Applications}, \textbf{505} (2016), 245--260.

\bibitem{CW2009}
Y. Chen and L. Wang, The Laplacian spread of tricyclic graphs, \emph{The Electronic Journal of Combinatorics}, \textbf{16}(1) (2009), \#R80.

\bibitem{CH2024}
F. E. Contreras-Mendoza and C. Hern\'andez-Cruz, $2$-polarity and algorithmic aspects of polarity variants on cograph superclasses, \emph{Discrete Mathematics \& Theoretical Computer Science}, \textbf{26}(3) (2024), \#11.

\bibitem{DJ2018}
R. R. Del-Vecchio and \'{A}.~A. Jones, Laplacian integrality in $P_{4}$-sparse and $P_{4}$-extendible graphs, \emph{Applied Mathematics and Computation}, \textbf{330} (2018), 307--315.

\bibitem{EK2021}
M. Einollahzadeh and M. M. Karkhaneei, On the lower bound of the sum of the algebraic connectivity of a graph and its complement, \emph{Journal of Combinatorial Theory, Series B}, \textbf{151} (2021), 235--249.

\bibitem{Fan2008}
Y.-Z. Fan, J. Xu, Y. Wang, and D. Liang, The Laplacian spread of a tree, \emph{Discrete Mathematics \& Theoretical Computer Science}, \textbf{10}(1) (2008), 79--86.

\bibitem{FLT2010}
Y.-Z. Fan, S.-D. Li, and Y.-Y. Tan, The Laplacian spread of bicyclic graphs, \emph{Journal of Mathematical Research and Exposition}, \textbf{30}(1) (2010), 17--28.

\bibitem{Fiedler1973}
M. Fiedler, Algebraic connectivity of graphs, \emph{Czechoslovak Mathematical Journal}, \textbf{23}(2) (1973), 298--305.

\bibitem{Grijo2019}
R. Grij\'o, L. S. de Lima, C. S. Oliveira, G. Porto, and V. Trevisan, Nordhaus--Gaddum type inequalities for the two largest Laplacian eigenvalues, \emph{Discrete Applied Mathematics}, \textbf{267} (2019), 176--183.

\bibitem{GMS1990}
R. Grone, R. Merris, and V. S. Sunder, The Laplacian spectrum of a graph, \emph{SIAM Journal on Matrix Analysis and Applications}, \textbf{11}(2) (1990), 218--238.

\bibitem{Haynsworth1968}
E. V. Haynsworth, Determination of the inertia of a partitioned Hermitian matrix, \emph{Linear Algebra and its Applications}, \textbf{1}(1) (1968), 73--81.

\bibitem{JYC2015}
G. Jiang, G. Yu, and J. Cao, The least algebraic connectivity of graphs, \emph{Discrete Dynamics in Nature and Society}, \textbf{2015} (2015), Article ID 756960.

\bibitem{Kirkland2010}
S. Kirkland, Algebraic connectivity for vertex-deleted subgraphs, and a notion of vertex centrality, \emph{Discrete Mathematics}, \textbf{310}(4) (2010), 911--921.

\bibitem{Kolokolnikov2015}
T. Kolokolnikov, Maximizing algebraic connectivity for certain families of graphs, \emph{Linear Algebra and its Applications}, \textbf{471} (2015), 122--140.

\bibitem{Konig1931}
D. K\H{o}nig, Graphok \'es matrixok, \emph{Matematikai \'es Fizikai Lapok}, \textbf{38} (1931), 116--119 (in Hungarian).

\bibitem{LiGuo2024}
W.-J. Li and J.-M. Guo, Nordhaus-Gaddum Type Inequalities for the $k$th Largest Laplacian Eigenvalues, \emph{The Electronic Journal of Combinatorics}, \textbf{31}(1) (2024), \#P1.31.

\bibitem{Liu2010}
Y. Liu, The Laplacian spread of cactuses, \emph{Discrete Mathematics \& Theoretical Computer Science}, \textbf{12}(3) (2010), 35--40.

\bibitem{Mathil2022}
P. Mathil, B. Baloda, and J. Kumar, On the cozero-divisor graphs associated to rings, \emph{AKCE International Journal of Graphs and Combinatorics}, \textbf{19}(3) (2022), 238--248.

\bibitem{MTR2011}
A. Mohammadian and B. Tayfeh-Rezaie, Graphs with four distinct Laplacian eigenvalues, \emph{Journal of Algebraic Combinatorics}, \textbf{34}(4) (2011), 671--682.

\bibitem{NG1956}
E. A. Nordhaus and J. W. Gaddum, On complementary graphs, \emph{The American Mathematical Monthly}, \textbf{63}(3) (1956), 175--177.

\bibitem{XM2011}
Y. Xu and J. Meng, The Laplacian spread of quasi-tree graphs, \emph{Linear Algebra and its Applications}, \textbf{435}(1) (2011), 60--66.

\bibitem{YL2012}
Z. You and B. Liu, The Laplacian spread of graphs, \emph{Czechoslovak Mathematical Journal}, \textbf{62}(1) (2012), 155--168.

\bibitem{ZSH2011}
M. Zhai, J. Shu, and Y. Hong, On the Laplacian spread of graphs, \emph{Applied Mathematics Letters}, \textbf{24}(12) (2011), 2097--2101.

\end{thebibliography}
\end{document}